\documentclass[12pt]{iopart}
\usepackage[utf8]{inputenc}
\usepackage{setstack,amsthm,iopams,bm,dirtytalk,float, subcaption,graphicx} 
\newtheorem{thm}{Theorem}

\newtheorem{result}[thm]{Result}

\begin{document}

\title[Shrinking without rotation in  high dimensional GLM regression]{Penalization-induced shrinking without rotation in  high dimensional GLM regression: a cavity analysis}

\author{E. Massa$^\dag$, M.A. Jonker$^\ddag$, A.C.C. Coolen$^{\dag,\S}$}

\address{$\dag$ Biophysics Department, Donders Institute, Radboud University,\\ \hspace*{2mm} 6525AJ  Nijmegen, The Netherlands
\\
$\ddag$ Department for Health Evidence, Radboud UMC, Geert Grooteplein 21,\\ \hspace*{2mm}  6525EZ Nijmegen, The Netherlands
\\
$\S$ Saddle Point Science Europe, Universitair Bedrijven Centrum, Toernooiveld 100,\\ \hspace*{2mm}  6525EC Nijmegen, The Netherlands}

\ead{emanuele.massa@donders.ru.nl, marianne.jonker@radboudumc.nl, a.coolen@science.ru.nl}

\begin{abstract}
    In high dimensional regression, where  the number of covariates  is of the order of the number of observations, ridge penalization is often used as  a  remedy against overfitting. Unfortunately, for correlated covariates such regularisation typically induces in generalized linear models not only shrinking of the estimated parameter vector, but also an unwanted \emph{rotation} relative to the true vector. 
    We show analytically how this problem can be removed by using a generalization of ridge penalization, and we analyse the asymptotic properties of the corresponding estimators in the high dimensional regime, using the cavity method. Our results also  provide a quantitative rationale for tuning the parameter that controlling the amount of shrinking.
    We compare our theoretical predictions with simulated data and find excellent agreement.
\end{abstract}

%

\section{Introduction}
\subsection{Setting and motivation}
In statistical regression the goal is to understand the relationship between a response $T\in \mathbb{R}$ and a set of covariates $\mathbf{X}\in \mathbb{R}^p$, given $n$ previously recorded observations. 
In generalized linear models (GLM) our hypotheses are encoded in a statistical model 
\begin{equation}
    p(T|\mathbf{X}\cdot\bm{\beta},\bm{\sigma})
    \label{model}
\end{equation} 
where $\bbeta\in \mathbb{R}^p$ are the regression parameters and $\bm{\sigma}\in \mathbb{R}^d$ are the nuisance parameters.
For simplicity we assume that the covariates follow a multivariate normal distribution
\begin{equation}
    \mathbf{X}\sim \mathcal{N}(\bm{0},\mathbf{A}^{-1}_0)
\end{equation}
and that our model is correctly specified, so our $n$ observations $\{(T_1,\mathbf{X}_1),\dots,(T_n,\mathbf{X}_n)\}$ are actually generated according to (\ref{model}), i.e.  $p(T_i|\mathbf{X}_i\cdot\bm{\beta}_0,\bm{\sigma}_0)$ 
for some \say{true} unknown values $\bbeta_0$ and $\bsigma_0$ that we would like to estimate.
The parameter estimators  $\hat{\bbeta}_n,\hat{\bsigma}_n$ are found by optimizing an objective function constructed from the observations available
\begin{eqnarray}
    \hat{\bbeta}_n,\hat{\bsigma}_n &:=& \underset{\bbeta,\bsigma}{\arg\max} \Big\{ l_n(\bbeta,\bsigma) \Big\}
    \label{PMLE}\\
    l_n(\bbeta,\bsigma)&=&\sum_{i=1}^n \log p (T_i|\mathbf{X}_i\cdot\bm{\beta},\bm{\sigma}) +   \pi(\bm{\beta},\bm{\sigma})
    \label{PLLE}
\end{eqnarray}
In (\ref{PLLE}) the first term on the right hand side is the log-likelihood of the observations under the model (\ref{model}) and $\pi$ is the penalization function.
This term is usually added in order to enforce a constraint by Lagrange multipliers.

In the high dimensional regime, the number of covariates $p$  in the model is large and proportional to the number of observations $n$, so that as $n,p\rightarrow \infty$ the ratio $\zeta := p/n$ is fixed. 
The true regression parameter vector $\bbeta_0$ is assumed to have a modulus $S_0 := \|\bbeta_0\|= O(1)$. A common choice for $\pi(\bm{\beta},\bm{\sigma})$ is the so-called ridge penalty
\begin{equation}
   \pi_{\rm ridge}(\bm{\beta},\bm{\sigma}|\Delta)=-\frac{1}{2}\Delta \bm{\beta}\cdot\bm{\beta}
 \label{ridge}
\end{equation}
which forces the estimator $\hat{\bm{\beta}}_n$ to take values closer to the origin. The effect is more pronounced for larger values of  $\Delta$. Ridge penalization is often employed  to prevent overfitting, and  thereby improve the prediction performance of the fitted model \cite{harrell2}. Overfitting causes the Maximum  Likelihood (ML) estimator $\hat{\bbeta}^{ML}_n$ of GLMs to lie typically in the same direction of $\bbeta_0$, but with modulus larger than that of $\bbeta_0$ \cite{GLM,massa}, and  ridge penalization counter-acts this phenomenon.
When the covariates are correlated, however, the effect of the penalty (\ref{ridge}) is to return an estimator that is both a shrunken and rotated version of the ML estimator:  $\hat{\bm{\beta}}_n$ differs on average both in  magnitude and in {\em direction} from $\bm{\beta}_0$ \cite{GLM,sheik}. In fact, ridge penalization was originally introduced for the purpose of obtaining a lower variance of $\hat{\bbeta}_n$ by \say{trading} variance for bias \cite{Hoerl}. While the amount of rotation induced by the penalization is in principle predictable \cite{GLM}, computing the latter is in practice hard to do, because it requires the computation of averages over the spectrum of the population covariance matrix $\mathbf{A}_0$, which is unknown and not easy to estimate \cite{ledoit_estimation,karoui_estimation}. 
Furthermore, while the Mean Squared Error of the estimator might be lower than the Maximum Likelihood one, it is evident that one would still report a value of $\hat{\bbeta}_n$ that is on average not even in the direction of $\bbeta_0$. This represents a major problem in inference if one wants to use penalization as a remedy against overfitting.
These considerations leads us to look for a generalization of the ridge penalty such that the resulting estimator: 1) is geometrically unbiased, that is, it lies typically in the  direction of $\bbeta_0$, and 2) has an amplitude tuned by a penalization parameter, so that by selecting the proper value of the latter we can make $\hat{\bbeta}_n$ unbiased.

\subsection{Covariant penalties}

In absence of the ridge penalization (i.e.\ in ML regression) one can formally {map} the problem for  correlated covariates to one where the transformed covariates are uncorrelated. This is possible because the ML estimator has the following covariant property: if we rotate and rescale our covariates according to the linear transformation $\bm{\mathcal{X}}_i:=\mathbf{M}^{-1}\mathbf{X}_i$, then   $\hat{\bbeta}_n\big(\bbeta_0,\{T_i,\mathbf{X}_i\}_{i=1}^n\big)$ satisfies
\begin{eqnarray} 
    \hat{\bbeta}_n\Big(\bbeta_0,\{T_i,\mathbf{X}_i\}_{i=1}^n\Big)&:=&  \underset{\bbeta}{\arg\max} \ \bigg(\underset{\bsigma}{\max}\Big\{ \sum_{i=1}^n \log p (T_i|\mathbf{X}_i\cdot\bm{\beta},\bm{\sigma}) \Big\}\bigg)\nonumber \\
    &~{=}&\mathbf{M}^{-1}\hat{\bbeta}_n \Big(\mathbf{M}\bbeta_0,\{T_i,\bm{\mathcal{X}}_i\}_{i=1}^n\Big)
    \label{property_ml}
\end{eqnarray}
The ML estimator $\hat{\bbeta}^{ML}_n$ transforms under rotation and rescaling of the covariates in the exact opposite way, so that the scalar product $\mathbf{X}_i\cdot\hat{\bbeta}_n$ does not change. Because the quantity $\bbeta_0$ that we want to estimate is a fixed vector in $\mathbb{R}^p$, a linear transformation executed on the underlying basis will change its components, hence  we must have $\mathbf{M}\bbeta_0$ in the right hand side. 
Upon transforming $\bm{\mathcal{X}}_i:=\mathbf{A}_0^{-1/2}\mathbf{X}_i$, where $\mathbf{A}_0$ the covariance matrix of $\mathbf{X}_i$, the transformed covariates will be  uncorrelated, $\bm{\mathcal{X}}_i\sim \mathcal{N}(\bm{0},\bm{I})$, and 
\begin{eqnarray} 
    \hat{\bbeta}_n\Big(\bbeta_0,\{T_i,\mathbf{X}_i\}_{i=1}^n\Big)&\overset{d}{=}&\mathbf{A}_0^{-1/2}\hat{\bbeta}_n \Big(\mathbf{A}_0^{1/2}\bbeta_0,\{T_i,\bm{\mathcal{X}}_i\}_{i=1}^n\Big) 
    \label{ml_uncorr}
\end{eqnarray}
Since (\ref{property_ml}) does not hold in general for the ridge penalized estimator (only under rotations), the correlations cannot be \say{transformed out} as in (\ref{ml_uncorr}).  
This motivates our interest in  \emph{covariant} generalizations of ridge penalization, such that the resulting estimator $\hat{\bbeta}_n$ would satisfy (\ref{property_ml}). 
For simplicity and analytical tractability, we consider quadratic penalties. A first candidate is
\begin{equation}
    \pi_{\rm O}(\bm{\beta}|\eta,\mathbf{A}_0)=-\frac{1}{2} \eta \bm{\beta}\cdot\mathbf{A}_0\bm{\beta}
    \label{oracle_pen}
\end{equation}
which we call \say{oracle} penalty, as it requires the knowledge of the population covariance matrix $\mathbf{A}_0$. For uncorrelated covariates,  (\ref{oracle_pen}) reduces to the ridge penalty. For correlated covariates, (\ref{oracle_pen}) shrinks differently in different directions, depending on the covariate correlations. 
While in some cases the population under investigation is sufficiently well characterized, in most cases the matrix $\mathbf{A}_0$ is not known. A possible \say{quick and dirty} alternative is to use the empirical correlation matrix as an estimate of $\mathbf{A}_0$, giving
\begin{equation}
    \pi_{\rm E}(\bm{\beta},\bm{\sigma}|\tau)=-\frac{1}{2} \tau \bm{\beta}\cdot  \Big(\frac{1}{n}\sum_{i=1}^n\mathbf{X}_{i}\mathbf{X}_{i}\Big)\bm{\beta}
    \label{empirical_pen}
\end{equation}
We refer to (\ref{empirical_pen}) as the \say{empirical} penalty. 
We know from Random Matrix Theory (RMT) that the sample covariance matrix is not a {good} estimator of the population covariance matrix in the high dimensional regime \cite{rmt_vivo}, so there we have no guarantee that the estimators obtained using (\ref{oracle_pen}) and (\ref{empirical_pen}) will have similar properties. Furthermore, the sample covariance matrix develops eigenvalues approaching zero \cite{rmt_vivo,couillet} as $\zeta \rightarrow 1$,  and (\ref{empirical_pen}) will not define a well posed optimization problem. On the other hand, if $\mathbf{A}_0$ is non-singular, the {oracle} penalty defines a well posed optimization problem  for all  $\zeta>0$.

\subsection{Previous work}

The idea of biasing an estimator to decrease its Mean Squared Error, i.e.\ shrinking, goes back to \cite{stein}. 
Ridge penalization is the simplest implementation of this idea, equivalent to a constrained optimization where the length of the maximizer is fixed \cite{Hoerl}. 
Our goal here differs from the original application of shrinking in that we seek to make $\hat{\bbeta}_n$ unbiased.
The geometrical properties of the objective function play a role in determining the distribution of the estimator $\hat{\bbeta}_n$. We will show that our penalization functions (\ref{oracle_pen},\ref{empirical_pen}) give shrunken and {geometrically} unbiased estimators. Their origin can be traced  back to the Bayesian interpretation of uniform shrinkage \cite{Copas,copas2}; here we provide another interpretation of uniform shrinkage priors. We will show that the asymptotic properties of the resulting estimators can be computed analytically, under simple assumptions about the data generating process, for arbitrary generalized linear models (GLM). The covariant penalties (\ref{oracle_pen},\ref{empirical_pen}) provide a natural  route to asymptotically unbiased estimators in GLM inference, together with our asymptotic analytical results to select an appropriate value of the penalization parameter, even when the sample size is modest.

Several papers in literature deal with the asymptotic behaviour of M-estimators, i.e.\ those derived by optimizing a sample average of an objective function (see \cite{el_karoui1,GLM,Asymptotic}). While the techniques used vary, e.g.\ Random Matrix Theory, Cavity Method and  Replica Method, the conclusions are always expressed via self consistent equations that give the bias and variance of the estimators for the regression parameters. 
The replica approach  \cite{GLM,PH,sheik,massa} appears the most {flexible}, as it does not assume any specific form of the utility function and (unlike other approaches) allows one to find also the asymptotic expected estimators of nuisance parameters. 
On the other hand, the replica approach tends to obscure some assumptions needed for its results to hold. We show in \ref{appendix:derivation} that the Replica Symmetric (RS) equations can be obtained also via  another method inspired by the statistical physics analogy with optimization \cite{mezard}: the cavity method \cite{virasoro,Mezard_89}. We extend the results obtained with the replica method (\cite{massa,GLM}) to incorporate the penalties (\ref{oracle_pen},\ref{empirical_pen}), and identify the key assumptions needed to reach the result obtained using replicas.
Our results confirm again that {the cavity method is complementary to the replica method} \cite{Mezard_89}. However, \say{it is always easier to descend a mountain than to climb it}, and the replica method is always our first step. Alternative heuristic approaches, based on Random Matrix Theory, already appeared in the statistical literature too \cite{el_karoui1}. 
Rigorous proofs that confirm the RS equations have been carried out so far in restricted settings. In \cite{el_karoui_rigorous}, the authors proved rigorous results for robust linear M-estimators. In  \cite{barbier} the asymptotic behaviour of the estimators in a linear model with normally distributed errors and fixed variance was establish. 
Ultimately, all methods (including ours) rely on the concentration of measure phenomenon, i.e.\ that specific random variables fluctuate with very low probability around a deterministic value \cite{talagrand1,el_karoui_rigorous,barbier,Boucheron13}, and hence we support our assumptions with simulations. 
Our present results generalize previous ones by allowing for nonlinear regression models and for the presence of nuisance parameters. 

\subsection{Aim and structure of the paper}

In this paper we study the asymptotic properties of the Penalized Maximum Likelihood (PML) estimator, mathematically equivalent to the Maximum A  Posteriori (MAP) estimator, obtained by maximizing the objective function 
\begin{equation}
    \fl l_n(\bm{\beta},\bm{\sigma}|\eta',\tau',\mathbf{A}_0) = \sum_{i=1}^n \log p (T_i|\mathbf{X}_i\cdot\bm{\beta},\bm{\sigma}) -\frac{1}{2}p \bm{\beta}\cdot\Big( \tau'  \frac{1}{n}\sum_{i=1}^n \mathbf{X}_i\mathbf{X}_i +\eta' \mathbf{A}_0\Big)\bm{\beta} 
\label{pll}
\end{equation}
The penalization combines (\ref{oracle_pen}) and (\ref{empirical_pen}), with re-scaled penalization parameters $\eta= p\eta'$ and $\tau = p\tau'$ such that penalization and log-likelihood are both $O(n)$ in the high dimensional  limit \cite{sheik}. We aim to compare the asymptotic properties of the estimators derived by using the two different penalizations (\ref{oracle_pen}, \ref{empirical_pen}), so we will only be interested in the cases where either $\eta$ or $\tau$ equal zero. In particular we want to answer the following questions: do the {covariant} penalizations (\ref{oracle_pen}, \ref{empirical_pen}) lead to a {geometrically} unbiased PML estimator for $\bbeta_0$, what are the properties of the PML estimator in the high dimensional regime, and how do these depend on the penalization parameters?

The article is divided in four sections. In section \ref{section:theory} we present our results concerning the asymptotic properties of the estimator in the high dimensional regime, and 
show that the PML estimator $\hat{\bbeta}_n$ is typically oriented in the same direction as $\bbeta_0$. 
In section \ref{section:applications} we compare our theoretical  predictions with simulated data for two prototypical models: the Logit regression model for binary data, and the Weibull Proportional Hazard model for skewed data. 
We conclude in section \ref{section:conclusion} with a discussion of our results and their applicability, and we identify future directions of investigation.

\section{The Covariantly Penalized Maximum Likelihood  estimator}
\label{section:theory}

In this section we give an explicit expression for the PML/MAP  estimator, and derive its properties in the high dimensional limit $n,p \rightarrow \infty$, with fixed ratio $p/n = \zeta>0$.

\subsection{Representation of $\hat{\bbeta}_n$}
The estimator $\hat{\bbeta}_n$ obtained by maximizing (\ref{pll}) satisfies (\ref{property_ml}) and (\ref{ml_uncorr}). Hence we can simply study the properties of the following estimator, because the PML/MAP estimator with correlated covariates is simply a rotated and re-scaled version of it:
\begin{eqnarray}
\label{simple_est}
    \hat{\bbeta}_n^{\sim}:=\hat{\bbeta}_n(\bbeta_0^{\sim},\{T_i,\bm{\mathcal{X}}_i\}) \qquad \mbox{with} \quad \bm{\mathcal{X}}_i \sim \mathcal{N}(\bm{0},\bm{I}), \quad \bbeta_0^{\sim} &:=& \mathbf{A}_0^{1/2}\bbeta_0
\end{eqnarray} 
Following a similar strategy as in \cite{el_karoui1}, we argue in \ref{appendix:representation} that this estimator admits the following representation
\begin{equation}
\label{uncorr_repr}
   \hat{\bbeta}_n^{\sim} \overset{d}{=}    K_n \bbeta_0^{\sim} + V_n \mathbf{U}
\end{equation}
where 
\begin{eqnarray}
\label{overlaps}
    K_n &:=& R^{(n)}_{0,1}/\|\bbeta^{\sim}_0\|  \qquad R^{(n)}_{0,1} := (\bbeta_0^{\sim}\cdot\hat{\bbeta}^{\sim}_n)/\|\bbeta^{\sim}_0\|\\
    V^2_n &:=& R^{(n)}_{1,1}-K_n^2  \qquad R^{(n)}_{1,1} := \hat{\bbeta}^{\sim}_n\cdot\hat{\bbeta}^{\sim}_n
\end{eqnarray}
with $\mathbf{U}$ independent from $K_n$ and $V_n$, uniformly distributed on the unit sphere in the $p-1$ dimensional subspace of $\mathbb{R}^p$ that is orthogonal to $\bbeta_0^{\sim}$. 
Using (\ref{ml_uncorr}) and (\ref{uncorr_repr}) we have 
\begin{equation}
   \hat{\bbeta}_n\overset{d}{=}   K_n \bbeta_0 + V_n\mathbf{A}_0^{-1/2}\mathbf{U}
   \label{representation_PMLE}
\end{equation}
The PML/MAP estimator consists of two geometrically independent contributions: a signal component in the direction of $\bbeta_0$, and an error component  which is independent from and  orthogonal to the former. The expected value of the error component over the possible data-sets is zero.
Representation (\ref{representation_PMLE}) shows explicitly that the covariance matrix $\mathbf{A}_0$ multiplies the error term $\mathbf{U}$ and enters otherwise in the factor $ \|\bbeta_0^{\sim}\|:= \|\mathbf{A}_0^{1/2}\bbeta_0\|$ only.
Furthermore, (\ref{representation_PMLE}) tells us that the distribution of the PML/MAP estimator is fully determined by $R_{0,1}^{(n)}$ and $R_{1,1}^{(n)}$.
These are called \say{overlaps} in the statistical physics of disordered systems \cite{virasoro,parisi,mezard}. It is appealing that these quantities,  familiar from replica calculations, also have a geometric interpretation for finite $n$. The {generalized overlaps} $K_n$ and  $V_n$ determine the amplitude of the two independent components of $\hat{\bbeta}_n$: $K_n$ gives the (random) inflation factor with respect to the true value, and $V_n$ is the (random) factor controlling the size of the error term. 

With (\ref{representation_PMLE}) we can compute and understand the expected properties of the estimator $\hat{\bbeta}_n$.
For instance, the Mean Square Error (MSE) of the estimator $\hat{\bbeta}_n$ takes the form
\begin{eqnarray}
\label{risk}
    \mathbb{E}\big[\|\hat{\bbeta}_n-\bbeta_0\|^2\big] &=&  \mathbb{E}\big[\|\hat{\bbeta}_n-\mathbb{E}[\hat{\bbeta}_n]\|^2\big] + \|\mathbb{E}[\hat{\bbeta}_n]-\bbeta_0\|^2
\end{eqnarray}
It consists of two distinct terms, namely the squared bias
\begin{equation}
\label{bias_beta}
    \|\mathbb{E}[\hat{\bbeta}_n]-\bbeta_0\|^2 = (\mathbb{E}[K_n]-1)^2\|\bbeta_0\|^2
\end{equation}
and the variance 
\begin{equation}
\label{var_beta}
    \mathbb{E}\big[\|\hat{\bbeta}_n-\mathbb{E}[\hat{\bbeta}_n]\|^2\big] =  \big(\mathbb{E}[K_n^2]-\mathbb{E}[K_n]^2\big) \|\bbeta_0\|^2 + \mathbb{E}[V_n^2] \mathbb{E}[A_p^2]
\end{equation}
where we defined 
\begin{equation}
\label{a_p}
    A_p^2 := \mathbf{U}\cdot\mathbf{A}_0^{-1}\mathbf{U}
\end{equation}
and used the fact that $\mathbf{U}$ is orthogonal to $\bbeta_0$ and independent of $V_n$ (see \ref{appendix:representation}).
According to (\ref{var_beta}), the moments of the generalized overlaps $K_n,V_n$ weigh the contribution of the true signal strength $\|\bbeta_0\|$ and of the fluctuations $\mathbb{E}[A_p^2]$ to the MSE.

\subsection{The PML/MAP estimator in the high dimensional limit} 

The next step is to assume that $K_n,V_n$ and the estimators for the nuisance parameters $\hat{\bsigma}_n$ (as well as $A_p$), should concentrate, as $n,p\rightarrow \infty$ with $\zeta := p/n$ fixed, around deterministic values. These latter values are expected to satisfy a set of self consistent equations.
In this section we derive these self-consistent asymptotic equations. The derivation based on cavity arguments is presented in \ref{appendix:derivation}.
We introduce some variations of the standard cavity argument, leading to a faster derivation and allowing us to tackle also non-linear objective functions and nuisance parameters, thus we believe this derivation to be of interest on its own. We have obtained the same results by means of the replica method, as an independent verification.

\begin{result}[Replica Symmetric (RS) equations]
\label{result:mainresult}
Suppose we are given $n$ i.i.d.\ observations $\{(T_i,\mathbf{X}_i)\}_{i=1}^n$, generated as $\mathbf{X}_i\sim \mathcal{N}(\bm{0},\mathbf{A}_0)$ and $T_i|\mathbf{X}_i\sim p(T_i|\mathbf{X}_i\cdot\bbeta_0,\bsigma_0)$.
As $n,p\rightarrow \infty$ with fixed ratio $\zeta =p/n$, under regularity conditions on the objective function $l_n(\bbeta,\bsigma)$ as in (\ref{pll}), the random variables $K_n$, $V_n$ and $\hat{\bsigma}_n$ concentrate around deterministic values $w_{\star}/S$, $v_{\star}$ and $\bsigma_{\star}$ respectively. 
These deterministic values, together with $u_{\star}$ are obtained by solving the so-called Replica Symmetric (RS) equations
\begin{eqnarray}
    &&\zeta v^2 = \mathbb{E}_{T,Q,Z_0}\Big[\big(\xi-vQ-wZ_0\big)^2\Big]
    \label{RS1}\\
    &&v\Big(1 + \zeta(u^2\eta'-1)\Big) = \mathbb{E}_{T,Q,Z_0}\Big[\frac{\partial}{\partial Q}\xi\Big]
    \label{RS2}\\
    &&w \zeta =\mathbb{E}_{T,Q,Z_0}\Big[ \xi\frac{\partial}{\partial Z_0}  \log p(T|SZ_0,\bm{\sigma}_0)\Big] 
    \label{RS3}\\
    &&0=\mathbb{E}_{T,Q,Z_0}\Big[\nabla_{\bsigma}\rho(T|\xi,\bm{\sigma})\Big] 
    \label{RS4}
\end{eqnarray}
where $S^2 := \lim_{p\rightarrow \infty}\bbeta_0\cdot\mathbf{A}_0\bbeta_0=\lim_{p\rightarrow \infty}\|\bbeta_0^{\sim}\|^2$. In the equations above  we have $Z_0,Q \sim \mathcal{N}(0,1), \ Z_0\!\perp\! Q,$ and $T|Y_0 \sim p(.|SZ_0,\bm{\sigma}_0)$.
We also used the short hand $\xi:=\xi(vQ+wZ_0,u,\bm{\sigma},T,\tau')$ to denote the quantity
\begin{equation}
    \xi_{\star}(x,u,\bsigma,T,\tau^\prime) := \underset{y}{\arg\min} \Big\{ \frac{1}{2}\Big(\frac{y-x}{u}\Big)^2  -  \rho(T|y,\bm{\sigma}) \Big\} 
    \label{RS_eqs}
\end{equation}
known as the proximal mapping of the function $-\rho(T|.,\bsigma)$, where
\begin{equation}
    \rho(T|y,\bm{\sigma}) = \log p(T|y,\bm{\sigma}) - \frac{1}{2}\zeta\tau'y^2
\end{equation}
\end{result}

As a consequence of the above result, the asymptotic properties of the estimators $\hat{\bbeta}_n$ and $\hat{\bsigma}_n$ depend on the RS order parameters $v_{\star},w_{\star},\bsigma_{\star}$.
For instance, the asymptotic squared bias of $\hat{\bbeta}_n$ reads
\begin{equation}
\label{asymptotic_bias_beta}
     \lim_{n,p\rightarrow\infty \:\ p/n=\zeta} \|\mathbb{E}[\hat{\bbeta}_n]-\bbeta_0\|^2 = (w_{\star}/S-1)^2 S_0^2 
\end{equation}   
where we have defined 
    $S^2_0:= \lim_{p\rightarrow \infty} \|\bbeta_0\|^2$. 
The variability of $\hat{\bbeta}_n$ around its expected value takes asymptotically the following simple form: 
\begin{equation}
\label{asymptotic_var_beta}
    \lim_{n,p\rightarrow \infty \:\ p/n=\zeta}\mathbb{E}\big[\|\hat{\bbeta}_n-\mathbb{E}[\hat{\bbeta}_n]\|^2\big]= \lim_{n,p\rightarrow \infty \:\ p/n=\zeta} \mathbb{E}[V_n^2]\mathbb{E}[A^2_p]
    = v_{\star}^2\alpha^2 
\end{equation}
where 
\begin{equation}
\label{alpha_def}
    \alpha^2 := \lim_{p\rightarrow\infty} \mathbb{E}[A^2_p] = \lim_{p\rightarrow\infty} \frac{1}{p}\Tr(\mathbf{A}_0^{-1})
\end{equation}
The derivation of  (\ref{alpha_def}) in found in \ref{appendix:alpha}.
One observes that the fluctuations of $\hat{\bbeta}_n$ are asymptotically controlled only by the trace of the inverse of the covariance matrix.
Furthermore, equation (\ref{asymptotic_var_beta}) gives the statistical meaning of the order parameter $v_{\star}$. The latter is zero in the low-dimensional regime, i.e. for $\zeta=0$, where the number of components of $\bm{\beta}$ is finite and the {typical size} of the fluctuations of each component is $O(n^{-1/2})$ by standard asymptotic theory \cite{Asymptotic}. In the high dimensional regime $\zeta>0$ this is no longer true  \cite{el_karoui1,el_karoui1,PH,massa,GLM,sheik}.
Combining (\ref{asymptotic_bias_beta}) with (\ref{asymptotic_var_beta}), we obtain an expression for the (asymptotic) MSE of $\hat{\bbeta}_n$ 
\begin{eqnarray}
    \fl \lim_{n,p\rightarrow \infty \:\ p/n=\zeta}\mathbb{E}\big[\|\hat{\bbeta}_n-\bbeta_0\|^2\big] &=& \lim_{n,p\rightarrow \infty \:\ p/n=\zeta}\mathbb{E}\big[\|\hat{\bbeta}_n-\mathbb{E}[\hat{\bbeta}_n]\|^2\big] + \|\mathbb{E}[\hat{\bbeta}_n]-\bm{\beta}_0\|^2 \nonumber \\
    &=& v_{\star}^2\alpha^2 + S_0^2\big(w_{\star}/S- 1\big)^2
    \label{mse}
\end{eqnarray}
It is important to note that the data set enters only in three scalar quantities, namely $\alpha$, $S$ and $S_0$, but only $S$ determines the solution of the RS equations. 

Since $(u_{\star},v_{\star},w_{\star},\bm{\sigma}_{\star})$ depend implicitly on the regularization parameters $\tau'$ and $\eta'$, this is also true for the asymptotic bias and variance (and thus the MSE) of $\hat{\bm{\beta}}_n$, and for the asymptotic properties of the nuisance parameter estimators $\hat{\bsigma}_n$. One can therefore, in principle, define the {optimal} value of the regularization parameters by imposing an additional condition. 
In inference the focus is usually on obtaining an unbiased estimator for the regression parameters $\hat{\bbeta}_n$, possibly with minimum variance. 
In this case one can require that the asymptotic bias be zero, or equivalently that $ w_{\star}/S = 1$ and solve for the value of the penalization parameter that satisfies this condition.

\section{Application to selected regression models}
\label{section:applications}

In this section we compare the solution of the RS equations against simulated data for two popular regression models, namely the Logit regression model for binary data and the Weibull proportional hazard model for time to event data.
This confirms that asymptotically the overlaps fluctuate increasingly tightly around the values predicted by the theory. We also show that the theory describes accurately how the properties of the PML/MAP estimator depend on the regularization parameters, and that the RS equations can be used to select the parameter values that gives unbiased estimators. 

\subsection{Simulation protocol}

The data  used in this section were generated with the open source programming language R \cite{R}. The covariates in the simulations are always taken to be  $\mathbf{X}\sim\mathcal{N}(\bm{0},\mathbf{A}^{-1}_0)$.
In each simulation, with given $p$, the covariance matrix $\mathbf{A}_0$ is generated randomly as follows. First we generate a random orthogonal matrix $\mathbf{O}$ in $\mathbb{R}^{p\times p}$ with the R-package \texttt{pracma} \cite{R_pracma}. We then sample the $p$ eigenvalues of $\mathbf{A}_0$ from a uniform distribution with support in $[0.1,10.0]$ and store them in a diagonal matrix $\bm{\Lambda}_0$. Finally $\mathbf{A}_0 = \mathbf{O}\cdot\bm{\Lambda}_0\mathbf{O}$.
The true value of the regression parameters $\bbeta_0$ is fixed to 
$\bbeta_0 = \mathbf{e}_1$, the first 
basis vector of $\mathbb{R}^p$, for simplicity. 
All elements of the covariance matrix $\mathbf{A}_0$ are rescaled by a common factor in order to enforce that $S=\|\mathbf{A}_0^{1/2}\bbeta_0\|=1$.
The present definitions clearly differ from the simpler orthonormal design case $\mathbf{X} = \mathcal{N}(\bm{0},\bm{I})$. 
We assume that there is no model mis-specification, so the responses are generated according to the same regression model that is used to construct the likelihood.
The generalized overlaps are computed from their definitions. For $K_n$ we build on  (\ref{representation_PMLE}), so that asymptotically
\begin{equation}
    K_n =\hat{\bbeta}_n\cdot \bbeta_0/\|\bbeta_0\|
\end{equation}
We compute $V_n$ approximately using (\ref{asymptotic_var_beta}) 
\begin{equation}
    V_n = \hat{\bbeta}_n\Big(\mathbf{I}-\frac{\bbeta_0\bbeta_0}{\|\bbeta_0\|^2}\Big)\cdot\hat{\bbeta}_n/\alpha
\end{equation}
Although the theory was derived in the asymptotic limit, we carried out our simulations at the rather modest sample size of $n=200$. This is done on purpose, in order to show that the high dimensional asymptotic theory is suitable in practice to understand even the typical properties of estimators obtained from small sample size datasets.

\begin{figure}[t]
\begin{subfigure}{.5\textwidth}
  \includegraphics[width=\linewidth]{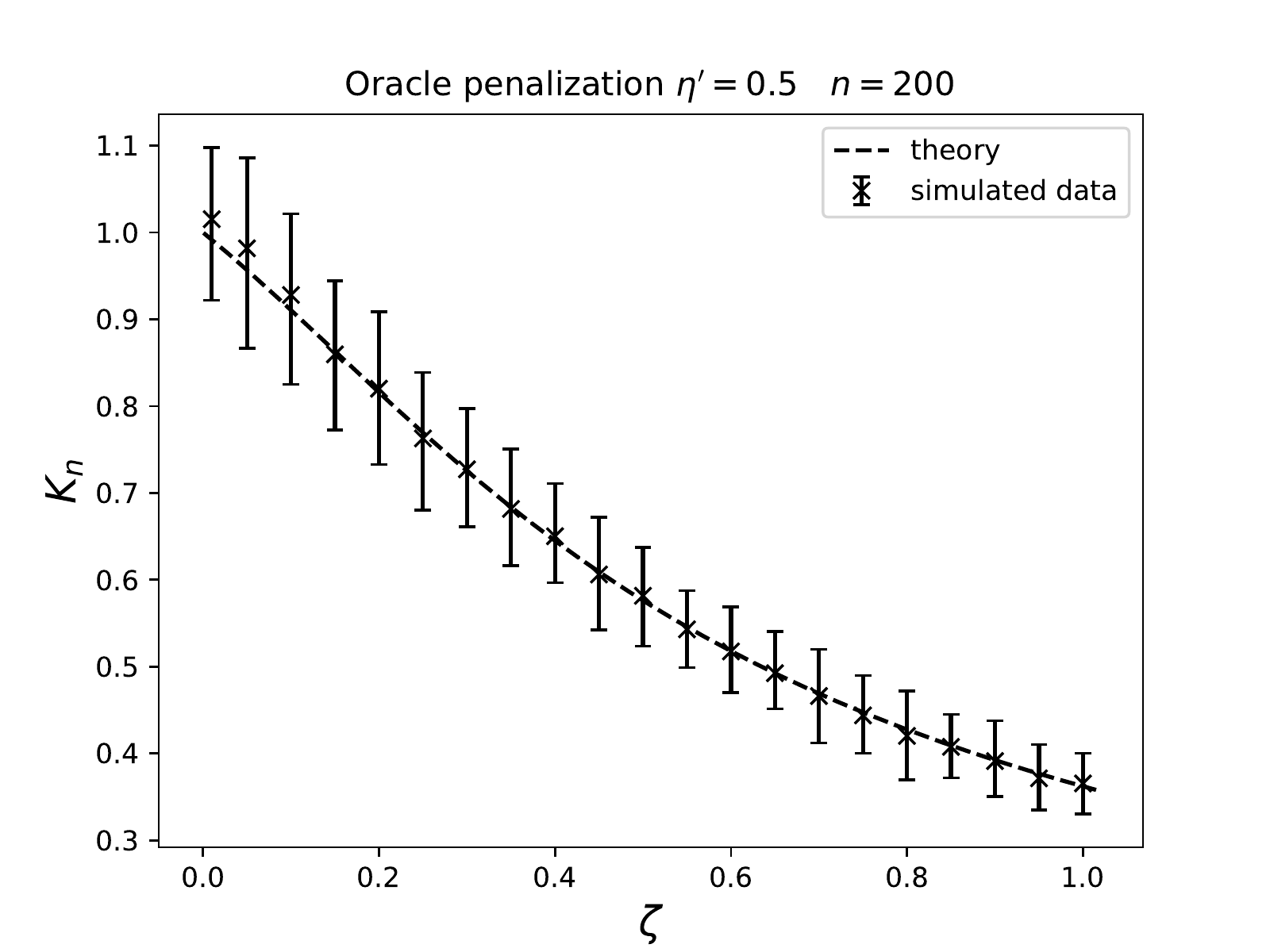}
  \caption{}
  \label{fig:1a}
\end{subfigure}%
\hfill
\begin{subfigure}{.5\textwidth}
  \includegraphics[width=\linewidth]{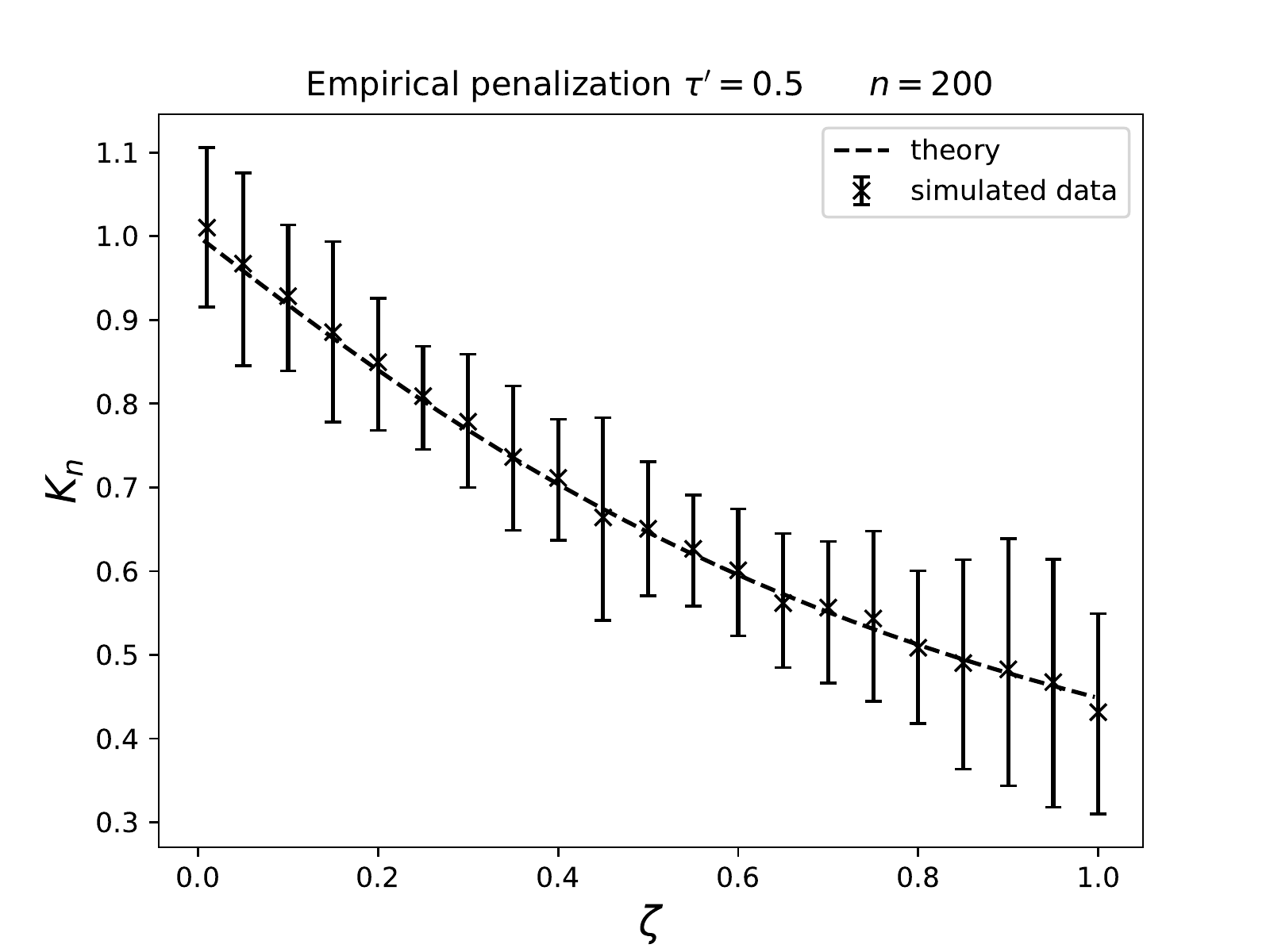}
  \caption{}
  \label{fig:1b}
\end{subfigure}
\medskip 
\begin{subfigure}{.5\textwidth}
  \includegraphics[width=\linewidth]{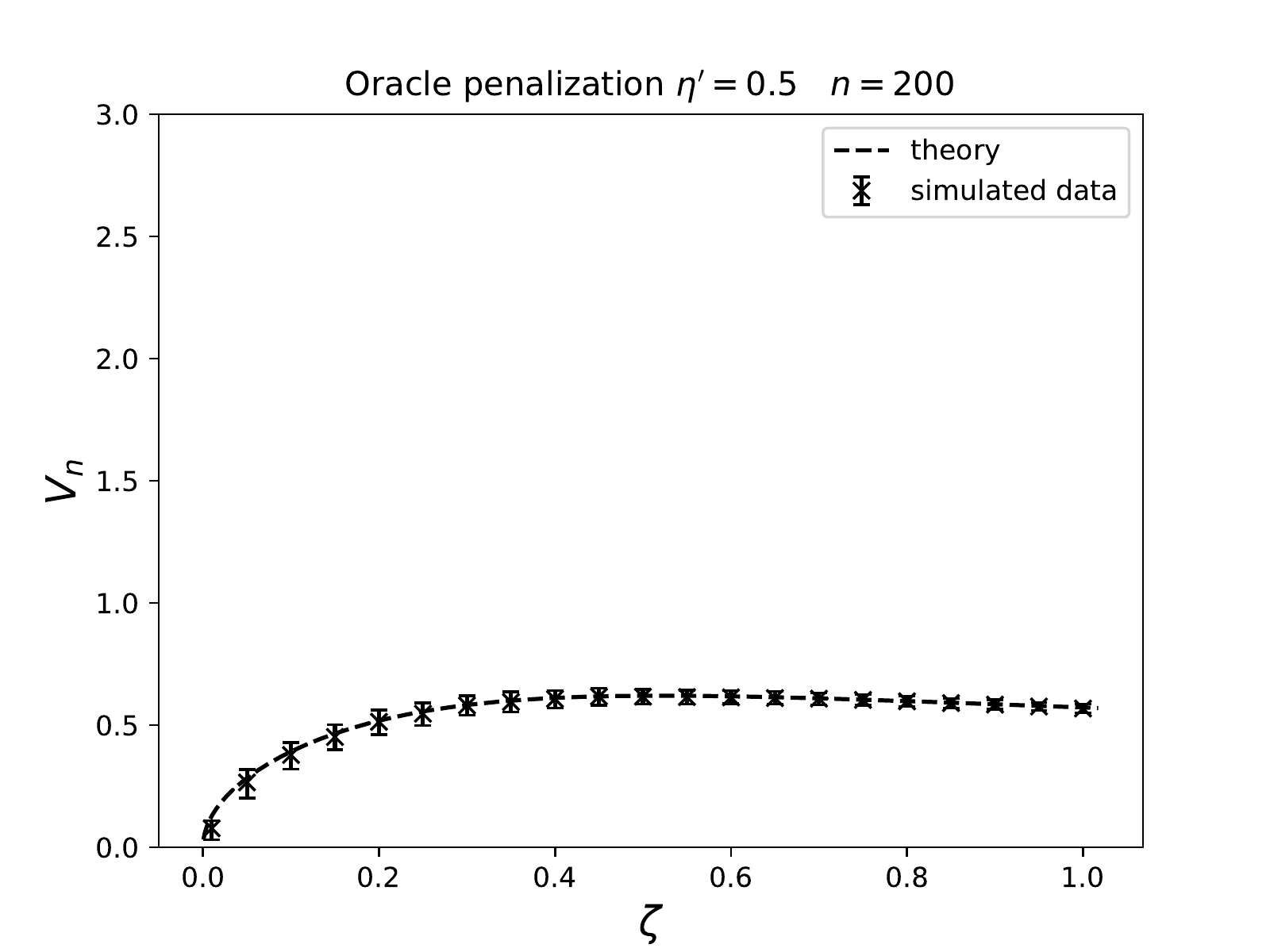}
  \caption{}
  \label{fig:1c}
\end{subfigure}%
\hfill
\begin{subfigure}{.5\textwidth}
  \includegraphics[width=\linewidth]{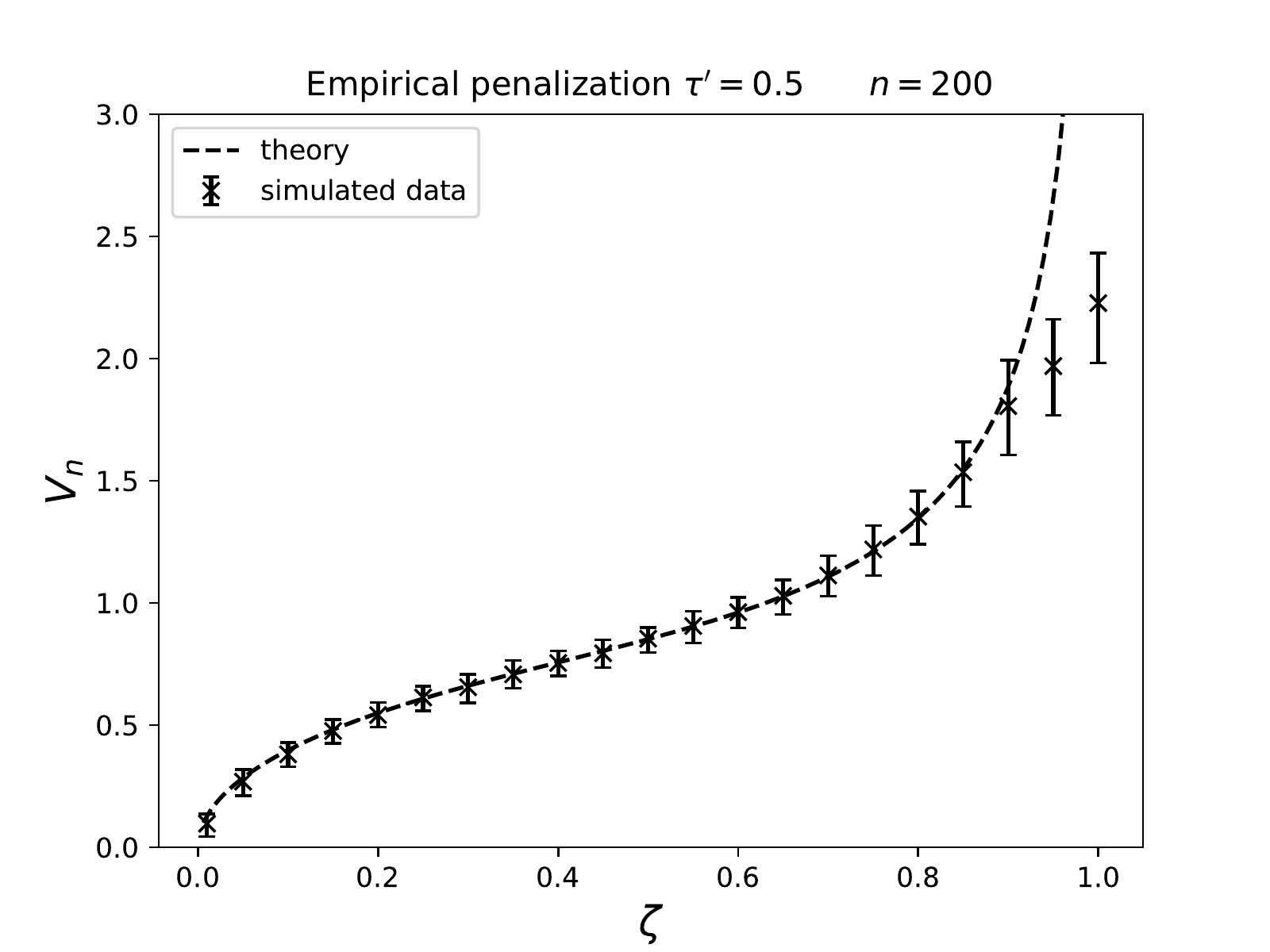}
  \caption{}
  \label{fig:1d}
\end{subfigure}
\caption{Simulated data for the Logit regression model, tested against the solution of the RS equations ($n=200$, $S=1.0$). The crosses represent sample averages, and the  error bars represent the first and third sample quartiles based on $500$ realizations of the PML/MAP estimator. Dashed curves: solution of the RS equations. Left: results for the covariant `oracle' regularizer. Right: results for the covariant `empirical' regularizer.  }
\label{fig:info_overlaps_logit}
\end{figure}

\subsection{Logit regression model}

The Logit regression model for a binary response $T \in \{-1,1\}$ is defined as
\begin{equation}
   p(T|\mathbf{X}\cdot\bbeta) = \frac{\rme^{T(\mathbf{X}\cdot\bm{\beta})}}{2\cosh(\mathbf{X}\cdot\bm{\beta})}
   \label{eq:logit}
\end{equation}
It is the simplest non-linear regression model $p(T|\mathbf{X}\cdot\bbeta,\bsigma)$ that cannot be rewritten in the form $p(T-\mathbf{X}\cdot\bbeta)$, and probably the most widely used model in applications.

It is known that for the model (\ref{eq:logit}) the MLE is biased for $\zeta>0$. Furthermore, this bias diverges at a critical value of $\zeta$, as was shown in previous studies, with a sharp phase transition  \cite{GLM,sheik,HD_logit} beyond which the MLE does not exist. 
Figure \ref{fig:info_overlaps_logit} shows that in PML regression for both covariant regularizers (\ref{oracle_pen}) (oracle) and (\ref{empirical_pen}) (empirical) the overlaps fluctuate around average values that are correctly predicted by the theory. The theoretical curves represented by dashed lines are the solution of the RS equations (see \ref{appendix:logit} for the derivation), reading 
\begin{eqnarray}
    \label{logit_RS1}
    \frac{\zeta\nu^2}{(1-\tau \mu^2)^2} &=& \mu^4~\mathbb{E}_{T,Q,Z_0}\Bigg[\bigg( T-\tau  \frac{\nu  Q+\omega Z_0}{1-\tau \mu^2} -  \tanh(x_{\star})\bigg)^2\Bigg]\\
    \label{logit_RS2}
    \zeta\Big(1-\frac{\mu^2\eta'}{1-\tau'\zeta\mu^2}\Big) &=& 1- (1-\tau'\zeta \mu^2) \mathbb{E}_{T,Q,Z_0}\Bigg[\frac{\cosh^2(x_{\star})}{u^2 + \cosh^2(x_{\star})}\Bigg] \\
    \label{logit_RS3}
    \zeta \omega/S &=& (1-\tau\mu^2) \mathbb{E}_{T,Q,Z_0}\Big[x\Big(T-\tanh(SZ_0)\Big)\Big]
\end{eqnarray}
where $Z_0,Q\sim\mathcal{N}(0,1)$, $Q\perp Z_0$, 
\begin{equation}
    T|Z_0\sim \frac{\rme^{T S Z_0}}{2\cosh(S Z_0)}
\end{equation}
and where $x_{\star}$ is the solution of the transcendental equation 
    $x = a - b\tanh(x)$, in which
\begin{eqnarray}
&&    a = \nu Q+\omega Z_0 +\mu^2 T,~~~~~~~
    b = \mu^2
\\
&&
\label{transformation}
    \mu^2 = u^2/(1\!+\!\tau'\zeta u^2),~~~~~~ \nu = v/(1\!+\!\tau'\zeta u^2), ~~~~~~ \omega = w/(1\!+\!\tau'\zeta u^2)
\end{eqnarray}
The above coupled equations can be solved numerically by fixed point iteration, once $S$  and the regularization parameter ($\tau'$ or $\eta'$) are specified. 

For $\zeta$ close to one, the sample covariance matrix becomes singular and the empirical penalization fails to regularize the optimization problem, leading to discrepancy between theory and simulated data (Fig \ref{fig:1d}).
The estimator obtained with the oracle penalization is slightly more biased than the one obtained with the empirical penalization, but has a smaller (and bounded) variance. Note that, away from $\zeta=1$, the fluctuations of the estimators obtained with the empirical penalty (Fig \ref{fig:1d}) are of the same order as those found with  oracle penalization (Fig \ref{fig:1c}).

\begin{figure}[t]
\begin{subfigure}{.47\textwidth}
  \includegraphics[width=\linewidth]{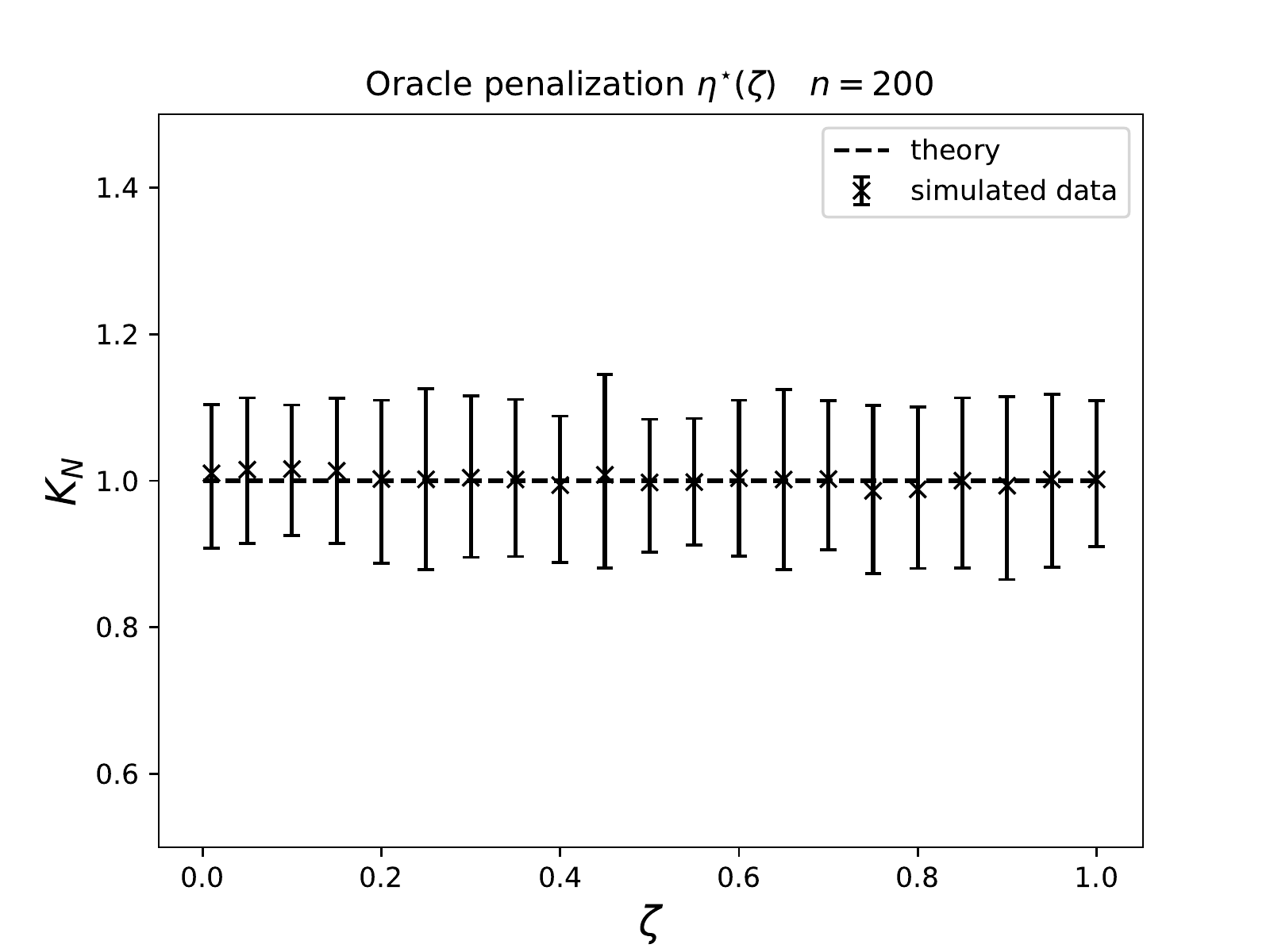}
  \caption{}
  \label{fig:2a}
\end{subfigure}%
\hfill
\begin{subfigure}{.47\textwidth}
  \includegraphics[width=\linewidth]{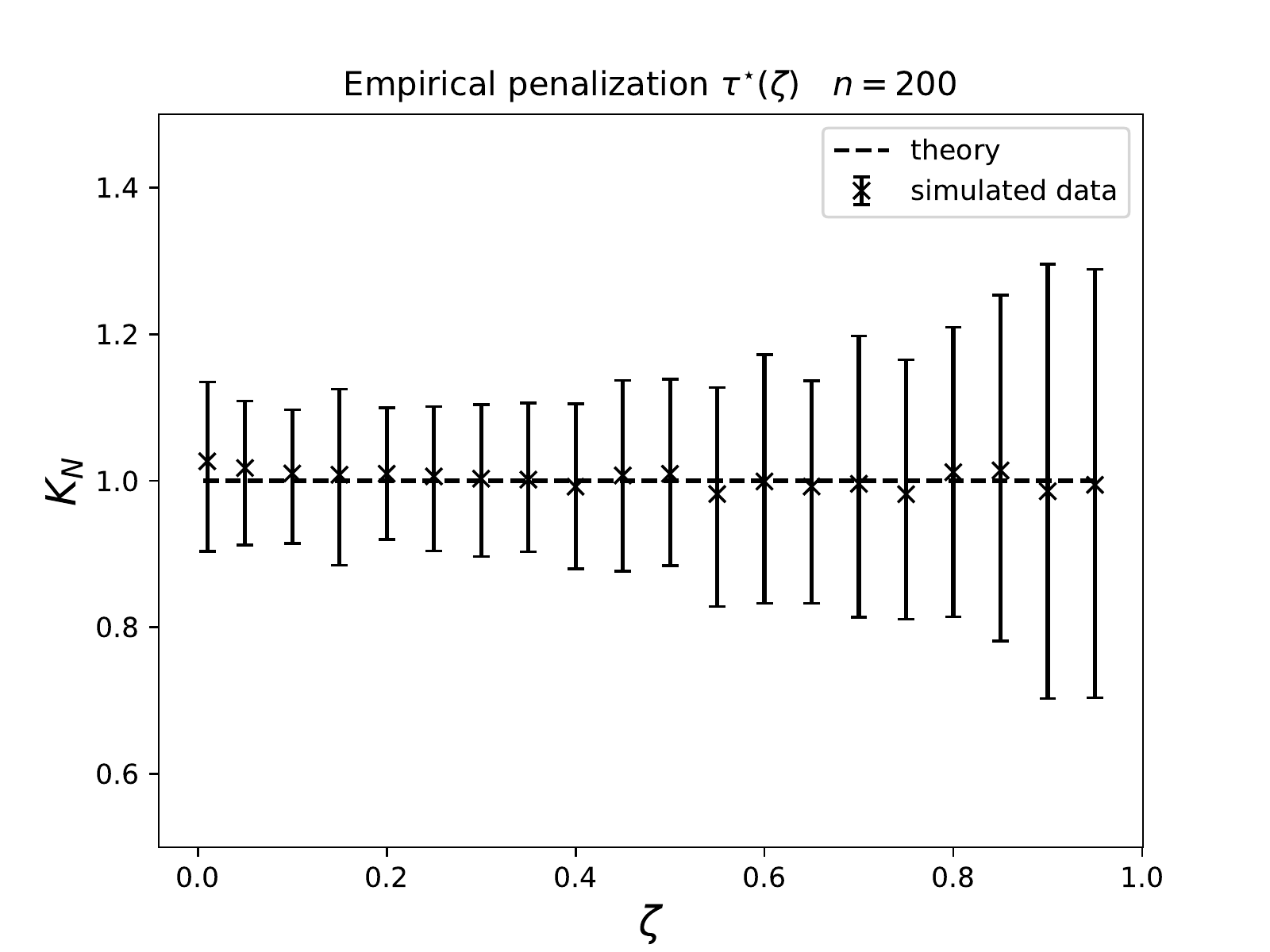}
  \caption{}
  \label{fig:2b}
\end{subfigure}
\medskip 
\begin{subfigure}{.47\textwidth}
  \includegraphics[width=\linewidth]{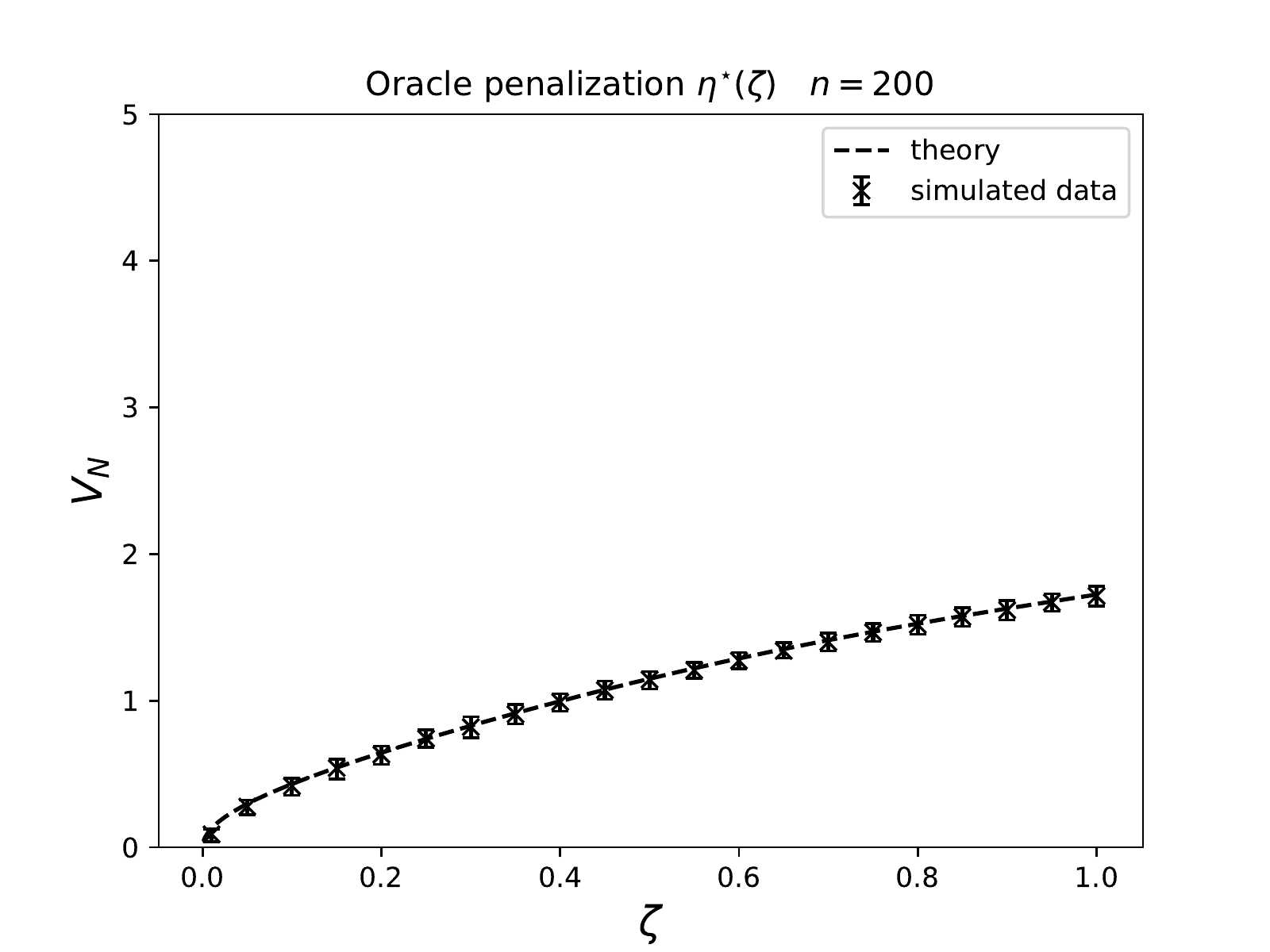}
  \caption{}
  \label{fig:2c}
\end{subfigure}%
\hfill
\begin{subfigure}{.47\textwidth}
  \includegraphics[width=\linewidth]{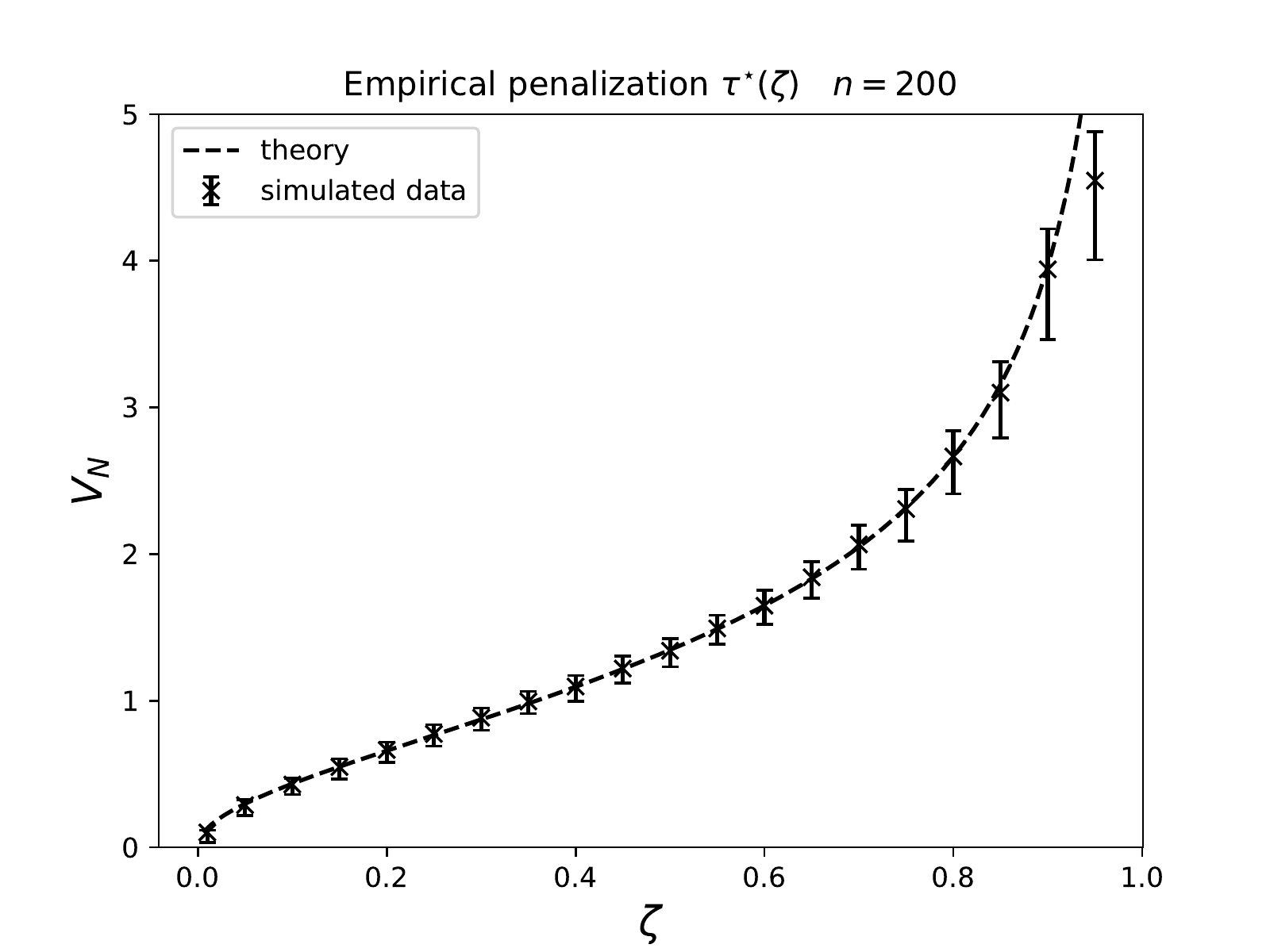}
  \caption{}
  \label{fig:2d}
\end{subfigure}
\caption{ Simulation data for the Logit regression model along the zero bias line (for $n$=200, $S=1$), with optimal oracle (left) or empirical (right) covariant penalization. The penalization parameters now depend on  $\zeta$ in such a way that the asymptotic bias is zero, i.e.  $w_{\star}/S=1.0$. This is confirmed in Figures  \ref{fig:2a} and \ref{fig:2b}. The estimator fluctuations are shown in Figures \ref{fig:2c} and \ref{fig:2d}. The error bars represent the first and third sample quartiles. Dashed lines give the theoretical predictions (solutions of the RS equations).}
\label{fig:zero_bias_lines}
\end{figure}

The special values of $\eta'$ and $\tau'$ that would make the estimator $\hat{\bbeta}_n$ strictly unbiased in both direction and amplitude, which we will denote by $\eta^\star$ and $\tau^\star$, can be computed as a function of $\zeta$ and $S$ by solving the RS equations, subject to the condition $w/S=\omega/S(1-\tau'\zeta\mu^2)^{-1}=1$.
This constraint gives 
\begin{equation}
    \zeta = \mathbb{E}_{T,Q,Z_0}\big[x(T-\tanh(SZ_0))\big]
\end{equation}
Equation (\ref{logit_RS2}) can then be used to compute  $\tau'$ or $\eta'$. For $\eta'=0$ one obtains
\begin{equation}
\label{tau_eq_RS}
    \tau^\star = \frac{\zeta- \chi}{1-\chi}\; \frac{1}{\mu^2\zeta}
~~~~~{\rm with}~~~~~
    \chi = \mathbb{E}_{T,Q,Z_0}\Big[\frac{u^2}{u^2 + \cosh^2(x_{\star})}\Big]
\end{equation}
Whilst for $\tau'=0$ one finds the prescription
\begin{equation}
\label{eta_eq_RS}
    \eta^\star = (\zeta- \chi)/\mu^2\zeta
\end{equation}
Equations (\ref{tau_eq_RS},\ref{eta_eq_RS}) are added to the RS equations in place of (\ref{logit_RS3}), upon which the combined system is once more solved by fixed point iteration.
The surfaces $\tau^\star(S,\zeta)$ (empirical covariant regularization) or $\eta^\star(S,\zeta)$ (oracle covariant regularization) define the values of the regularization parameters that should be used to obtain an unbiased estimator (that is, given $S$). 
\begin{figure}[t]
\begin{subfigure}{.5\textwidth}
  \includegraphics[width=\linewidth]{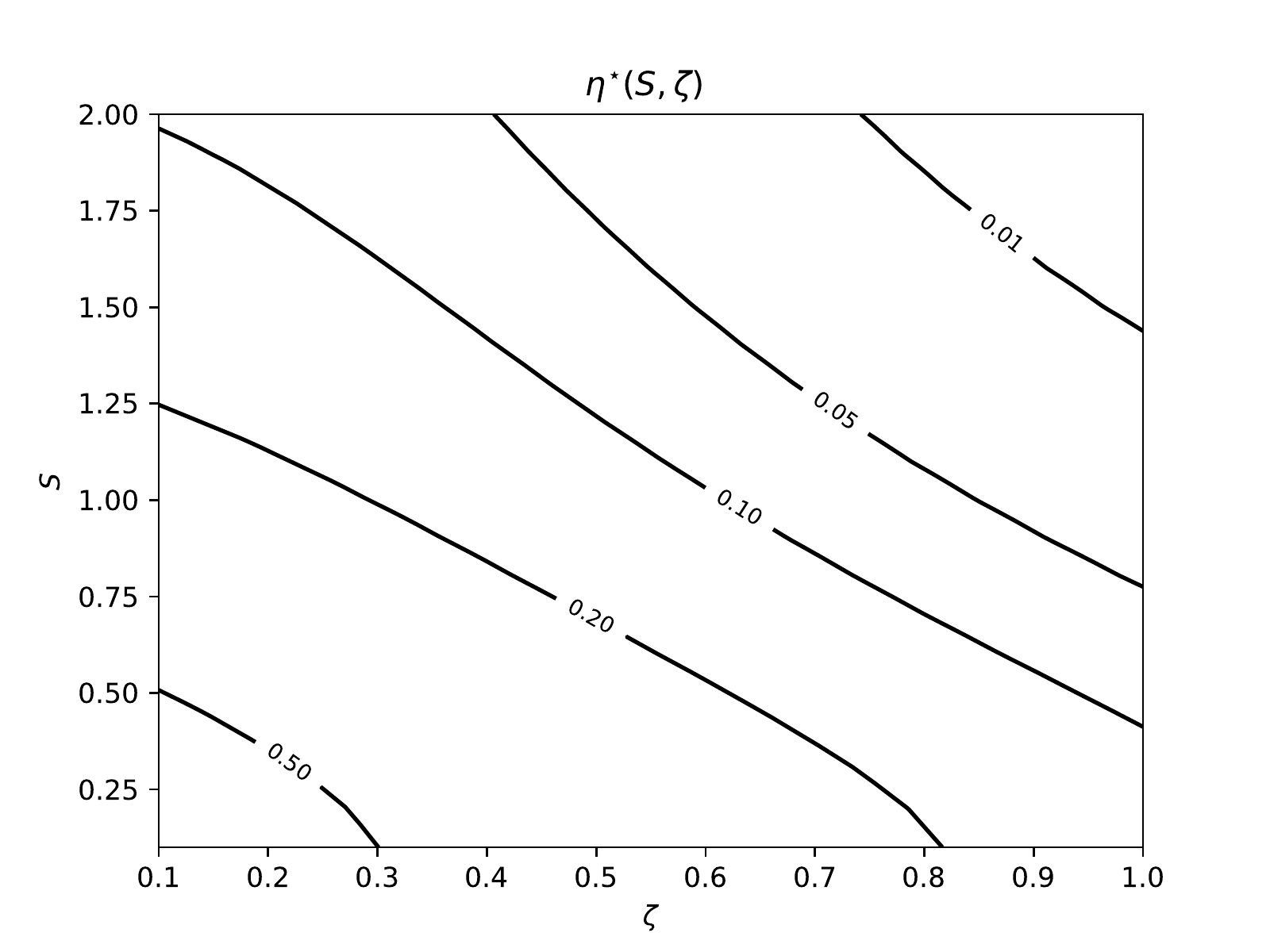}
  \caption{}
  \label{fig:3a}
\end{subfigure}%
\hfill
\begin{subfigure}{.5\textwidth}
  \includegraphics[width=\linewidth]{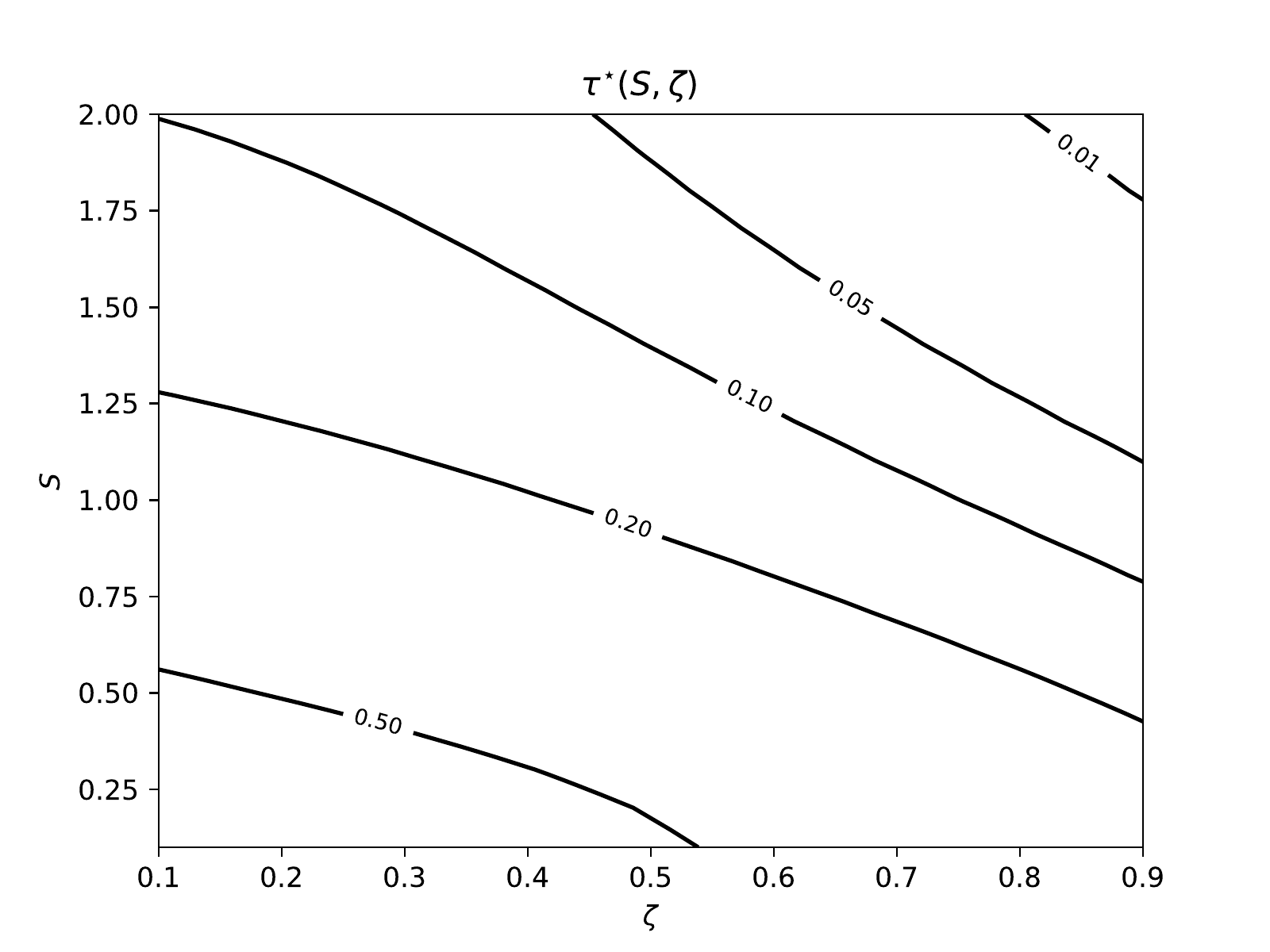}
  \caption{}
  \label{fig:3b}
\end{subfigure}
\caption{Contour plots of the surfaces $\tau^\star(S,\zeta)$ and $\eta^\star(S,\zeta)$. These are the regularizion parameter values of the covariant penalizations that make the estimator $\hat{\bbeta}_n$ unbiased. Left: oracle penalization. Right: empirical penalization. Both values are seen to decrease with increasing  $S$ and $\zeta$.}
\label{fig:unbiased_hp}
\end{figure}
In Figure \ref{fig:zero_bias_lines} we show the result of carrying out regression following the latter protocols and compare the values computed from simulations with the theoretical curves. We observe  satisfactory agreement, with the simulation averages almost identical to the theoretical predictions. The improvements over  ML estimator are clear: the estimator is unbiased (for known $S$), and has a finite variance. 
It is also clear that $V_n$, measuring the typical fluctuations around the average, gets larger as $\zeta$ increases, especially for the estimator obtained with the empirical penalization (Fig \ref{fig:2d}).
Nevertheless, it is remarkable and of practical relevance that, according to Figures  \ref{fig:2c} and \ref{fig:2d}, for sufficiently small values of $\zeta$ (e.g. $\zeta<0.4$) the sample to sample fluctuations  of the estimator obtained with empirical penalization is  comparable to the one obtained with oracle penalization. This implies that for such $\zeta$ values, in absence of $\mathbf{A}_0$ one can safely use empirical penalization.

It is interesting to inspect the dependence on the global parameters $S$ and $\zeta$ of the special values for $\tau^\star(S,\zeta)$ and $\eta^\star(S,\zeta)$ that lead to unbiased $\hat{\bbeta}_n$. These values are show as contour plots in the $(S, \zeta)$ plane in Figure \ref{fig:unbiased_hp}. We see that the de-biasing regularization parameters decrease with both $S$ and $\zeta$. Whilst the decrease with $S$ is intuitive (it is easier to estimate a powerful signal than a weak one, hence there is less need for correction), the  decrease with $\zeta$ is not. On the other hand, our theory predicts that $\hat{\bbeta}_n\cdot \mathbf{A}_0\hat{\bbeta}_n$ and $\sum_{i=1}^n \hat{\bbeta}_n\cdot\mathbf{X}_i\mathbf{X}_i\hat{\bbeta}_n/n$ both increase with $\zeta$. It is hence not unreasonable that the optimal penalization parameter gets smaller for larger $\zeta$, to prevent penalization from contributing the leading term in the objective function (which would force $\hat{\bbeta}_n$ excessively towards the origin). One would expect to find that $\eta^{\star}$ or $\tau^{\star}$ goes to zero with $\zeta$. This is actually not needed as the penalization parameters are multiplied by $p$ (see \ref{pll}), hence when $\zeta$ \say{goes} to zero, the same is true for $p = \zeta n$. In (\ref{pll}) the amount of penalization is then tuned by $\zeta \tau'$ (or $\zeta\eta'$) and  the penalty amounts to a small contribution in the objective function, no matter what finite value of $\tau'$ or $\eta'$, which are indeed found to be non vanishing, but finite. 
Notice finally,  that while it is possible to make the estimator $\hat{\bbeta}_n$ strictly {unbiased}, this may not always be meaningful. It is well known that in some cases a biased estimator with lower variance might outperform an unbiased one; for an unbiased estimator with excessive variance the signal may get lost in a sea of fluctuations. Once $v_\star$ and $\alpha$ are estimated, one has access to the variance of $\hat{\bbeta}_n$ (\ref{asymptotic_var_beta}). Alternative one can control the bias-variance trade-off manually via the solution of the RS equations.

\subsection{The Weibull proportional hazards model}

\begin{figure}[t]
\begin{subfigure}{.5\textwidth}
  \includegraphics[width=\linewidth]{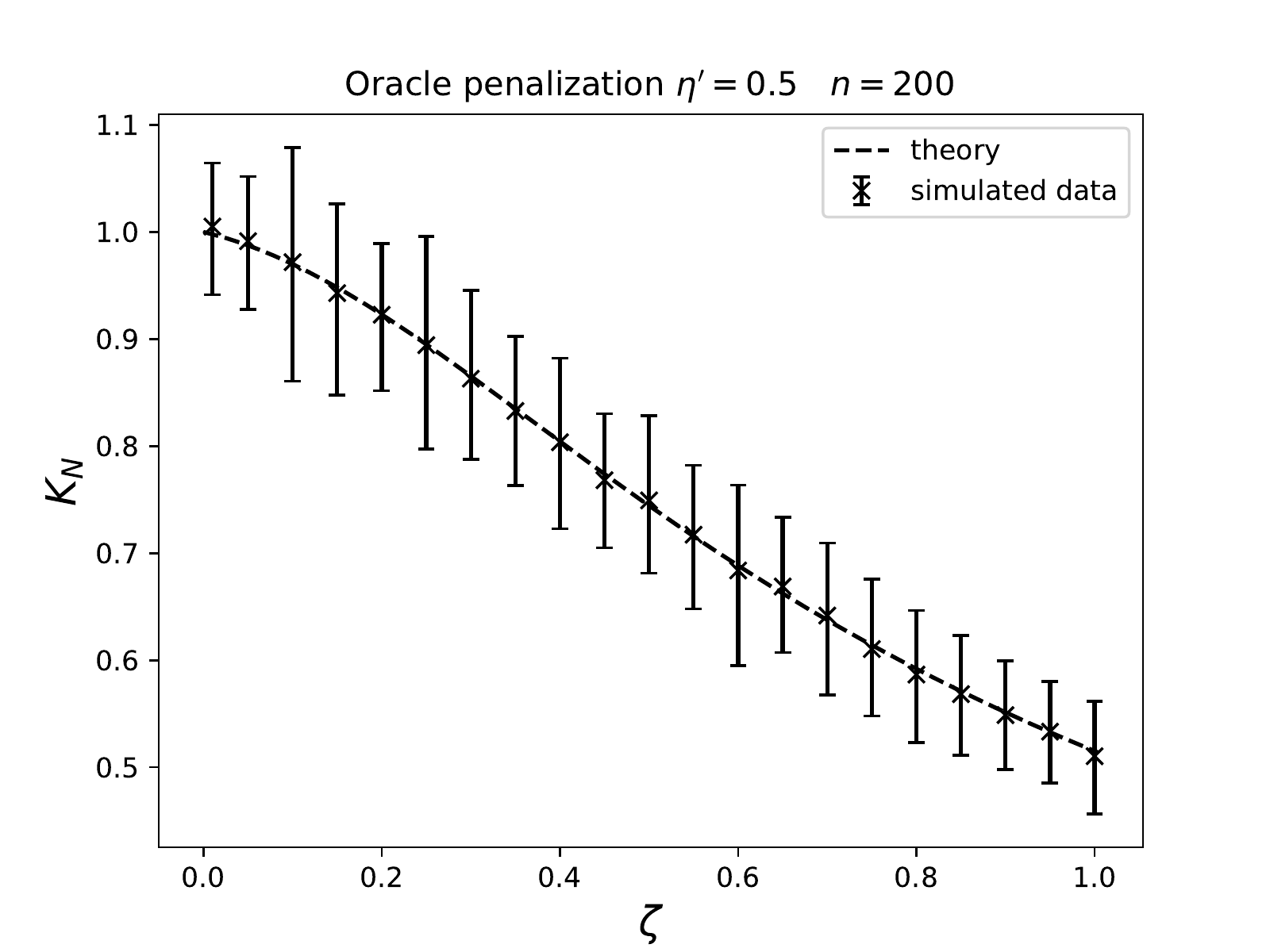}
  \caption{}
  \label{fig:5a}
\end{subfigure}%
\hfill
\begin{subfigure}{.5\textwidth}
  \includegraphics[width=\linewidth]{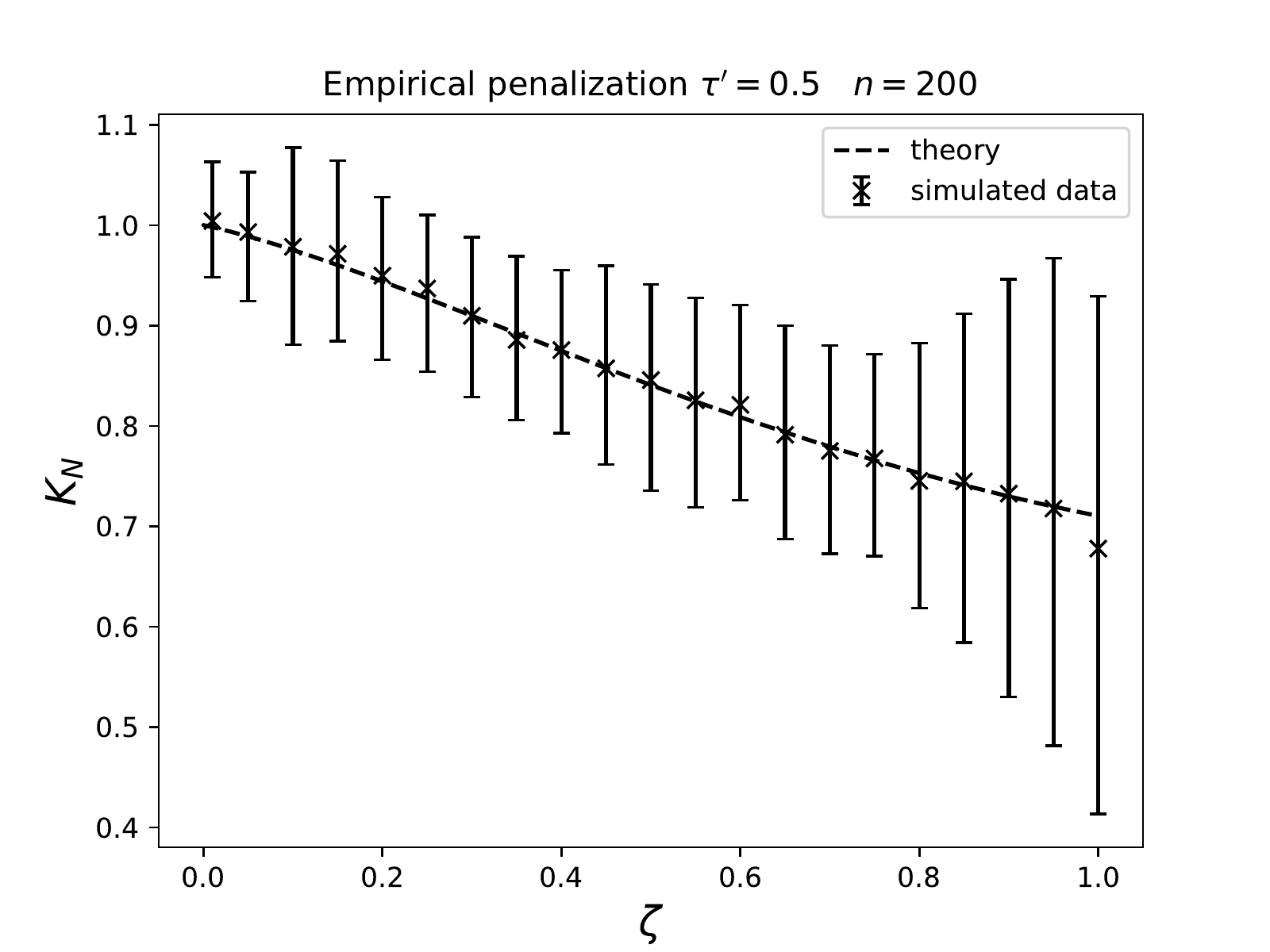}
  \caption{}
  \label{fig:5b}
\end{subfigure}
\medskip 
\begin{subfigure}{.5\textwidth}
  \includegraphics[width=\linewidth]{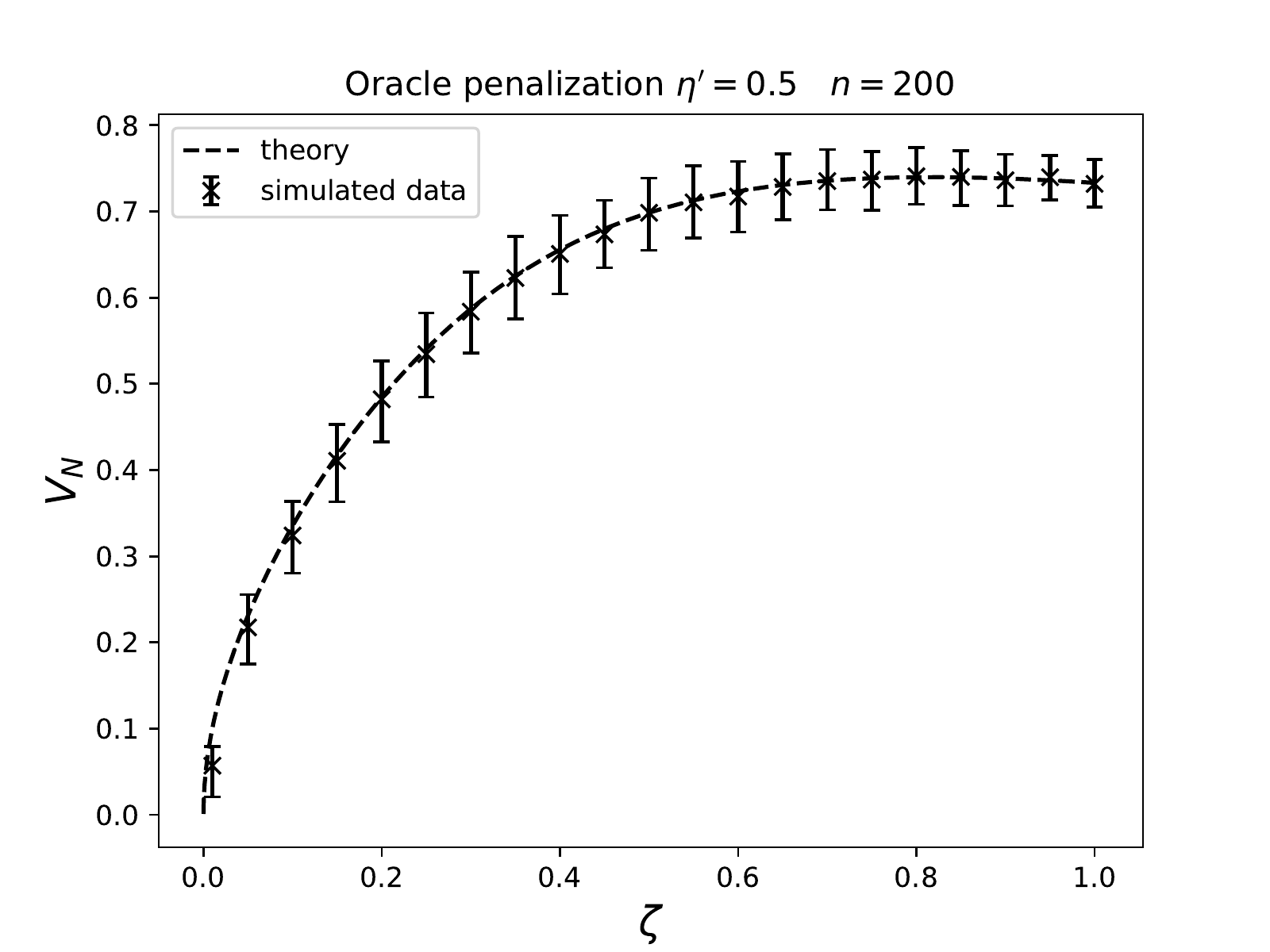}
  \caption{}
  \label{fig:5c}
\end{subfigure}%
\hfill
\begin{subfigure}{.5\textwidth}
  \includegraphics[width=\linewidth]{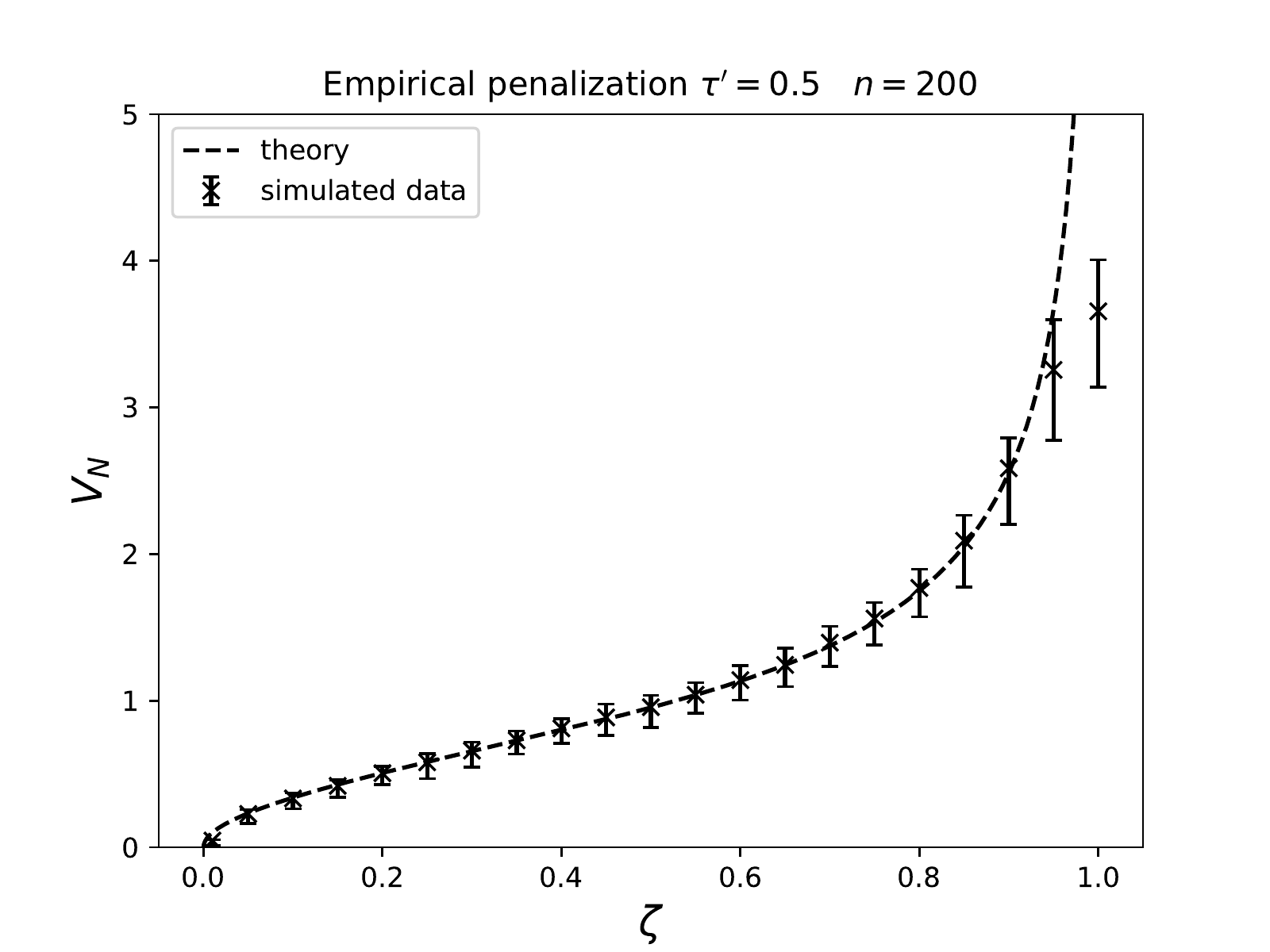}
  \caption{}
  \label{fig:5d}
\end{subfigure}
\caption{ Simulation data for the Weibull proportional hazards model (with $n$=200, $S=1$), with  oracle (left) or empirical (right) covariant penalization.  The error bars represent the first and third sample quartiles, while the dashed line is the curve obtained by solving the RS equations.
 Dashed lines give the theoretical predictions.
}
\label{fig:overlap_weib}
\end{figure}

In Proportional hazards models, the conditional density of the response $T\geq 0$, given the covariates, is of the form
\begin{equation}
    p(T|\mathbf{X}\cdot\bm{\beta},\bm{\sigma}) = h(T|\bm{\sigma})\rme^{\mathbf{X}\cdot\bm{\beta}-H(T|\bm{\sigma})\exp(\mathbf{X}\cdot\bm{\beta})}
\end{equation}
where $h(.|.)$ is the base hazard function, and 
$H(T|\bm{\sigma}) := \int_0^{T}\!\rmd t^\prime~ h(t^\prime|\bsigma)$
is the integrated base hazard rate.
A widely employed parametric survival model is the Weibull proportional hazards model, for which the integrated base hazard rate has the form\footnote{This is not the usual parametrization of the Weibull base hazard rate, which can be recovered by setting $\lambda= \exp(\phi)$ and $\rho= 1/\sigma$.}
\begin{equation}
    H(T|\phi,\sigma) = \rme^{\phi/\sigma} T^{1/\sigma}
\end{equation}
In \ref{appendix:weibull} we show that for this model the RS equations can be written in the form
\begin{eqnarray}
\label{weib1}
\hspace*{-20mm}
&& \frac{\nu^2\zeta}{(1-\tau'\zeta \mu^2)^2} =\\
\hspace*{-20mm}&&\hspace*{6mm}
        \mathbb{E}_{Z,Q,Z_0}\Big[\Big(\mu^2 - W\Big( Z^{\sigma_0/\sigma} \mu^2\rme^{\mu^2+(\omega -S\frac{\sigma_0}{\sigma})Z_0 + \nu Q + (\phi-\phi_0)/\sigma} \Big)-\frac{\tau'\zeta \mu^2}{1 - \tau'\zeta\mu^2}(\nu Q +\omega Z_0)\Big)^2\Big]\nonumber\\
\label{weib2}
\hspace*{-20mm}
&&    \frac{1}{1\!-\!\tau'\zeta\mu^2} \Big(1\!-\!\zeta\big(1\!-\!\frac{\mu^2\eta'}{1-\tau'\zeta\mu^2}\big)\Big) =
\\[-2mm]
\hspace*{-20mm}
        &&
        \hspace*{51mm}1-  \mathbb{E}_{Z,Q,Z_0}\Bigg[\frac{W\Big( Z^{\sigma_0/\sigma} \mu^2\rme^{\mu^2+(\omega -S\frac{\sigma_0}{\sigma})Z_0 + \nu Q + (\phi-\phi_0)/\sigma} \Big)}{1\!+\!W\Big( Z^{\sigma_0/\sigma} \mu^2\rme^{\mu^2+(\omega -S\frac{\sigma_0}{\sigma})Z_0 + \nu Q + (\phi-\phi_0)/\sigma} \Big)}\Bigg]\nonumber\\
\label{weib3}
\hspace*{-20mm}
      && \frac{\omega}{S} = \Big(1-\tau'\mu^2 - \frac{\eta'\mu^2}{1-\tau'\mu^2\zeta}\Big)\frac{\sigma_0}{\sigma}\\
\label{weib4}
\hspace*{-20mm}
&&       \mu^2 = \mathbb{E}_{Z,Q,Z_0}\Big[W\Big( Z^{\sigma_0/\sigma} \mu^2\rme^{\mu^2+(\omega -S\frac{\sigma_0}{\sigma})Z_0 + \nu Q + (\phi-\phi_0)/\sigma} \Big)\Big]\\
\label{weib5}
\hspace*{-20mm}
     &&\frac{\sigma}{\sigma_0} = \mathbb{E}_{Z,Q,Z_0}\Big[\Big(\log Z- SZ_0 \Big)\bigg(\frac{W\Big( Z^{\sigma_0/\sigma} \mu^2\rme^{\mu^2+(\omega -S\frac{\sigma_0}{\sigma})Z_0 + \nu Q + (\phi-\phi_0)/\sigma} \Big)}{\mu^2} -1 \bigg)\Big]
\end{eqnarray}
where $W(x)$ is Lambert's $W$-function, $Z \sim {\rm Exp}(1)$, $Q,Z_0\sim \mathcal{N}(0,1)$ with $Q\perp Z_0$ and $\mu,\omega,\nu$ are related to $u,w,v$ as in (\ref{transformation}). 
The above equations can be again solved by fixed point iteration. Similar to the  Logit model,  we see in Figure \ref{fig:overlap_weib} that the overlaps fluctuate closely around the theoretical prediction (i.e. the solution of the RS equations,  displayed as a dashed line). 
The estimator obtained with oracle penalization is again more biased with respect to the one obtained with the empirical penalization, but {fluctuates less} from sample to sample. These fluctuations of the two penalization formulae are of the same order of magnitude for sufficiently small $\zeta$.

To obtain {unbiased} estimators we impose the condition $w/S = \omega/S (1-\tau\zeta\mu^2)^{-1}=1$, which leads to an equation either for $\tau^\star$ or for $\eta^\star$, to be  substituted for (\ref{weib3}).
Setting $\eta'=0$ we obtain 
\begin{equation}
\label{tau_weib}
    \tau^\star = \frac{1-\sigma/\sigma_0}{1-\zeta\sigma/\sigma_0}
\end{equation}
while setting $\tau'=0$ we obtain 
\begin{equation}
\label{eta_weib}
    \eta^\star =\mu^{-2}(1-\sigma/\sigma_0)
\end{equation}
\begin{figure}[t]
\begin{subfigure}{.5\textwidth}
  \includegraphics[width=\linewidth]{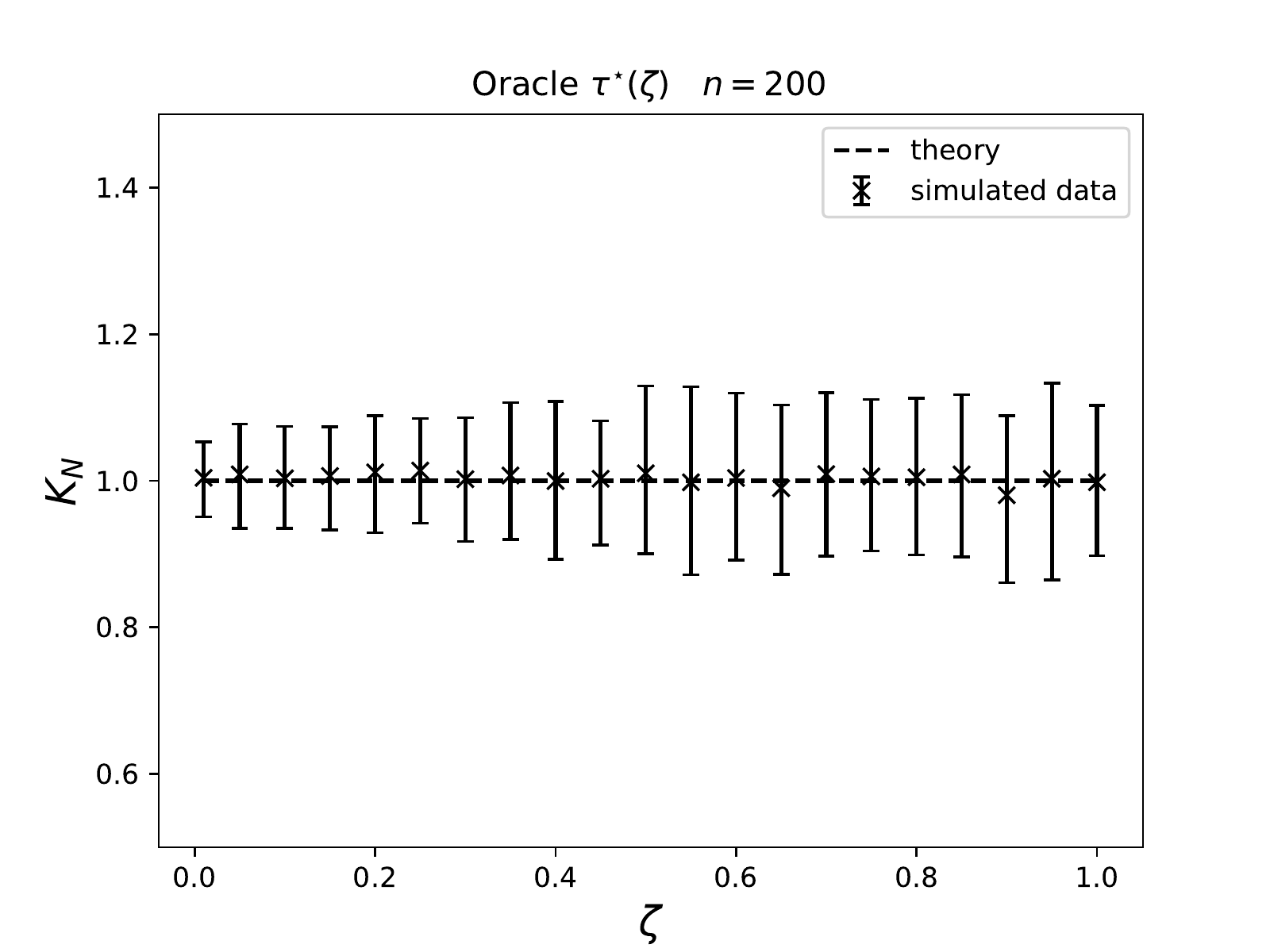}
  \caption{}
  \label{fig:6a}
\end{subfigure}%
\hfill
\begin{subfigure}{.5\textwidth}
  \includegraphics[width=\linewidth]{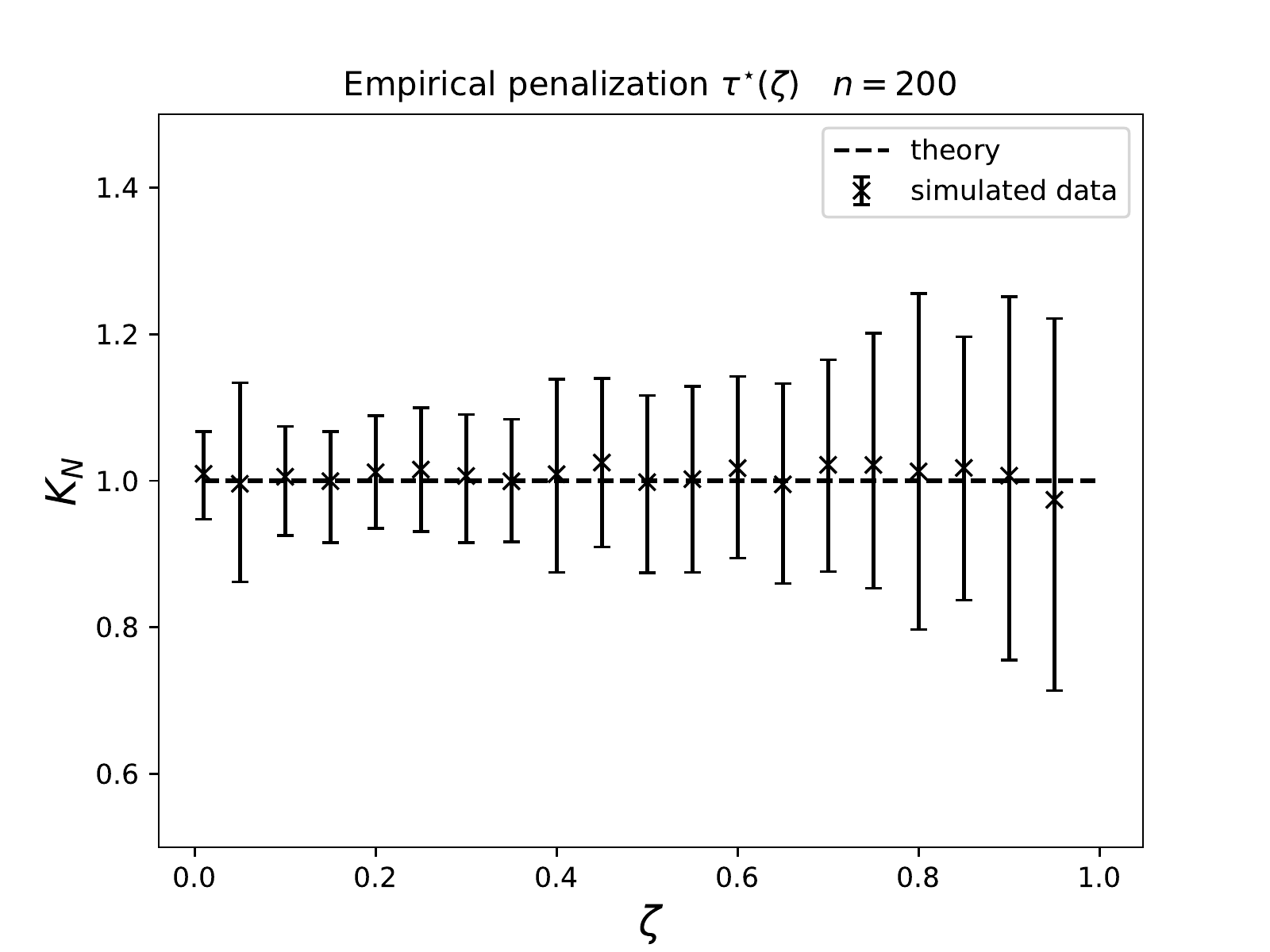}
  \caption{}
  \label{fig:6b}
\end{subfigure}
\medskip 
\begin{subfigure}{.5\textwidth}
  \includegraphics[width=\linewidth]{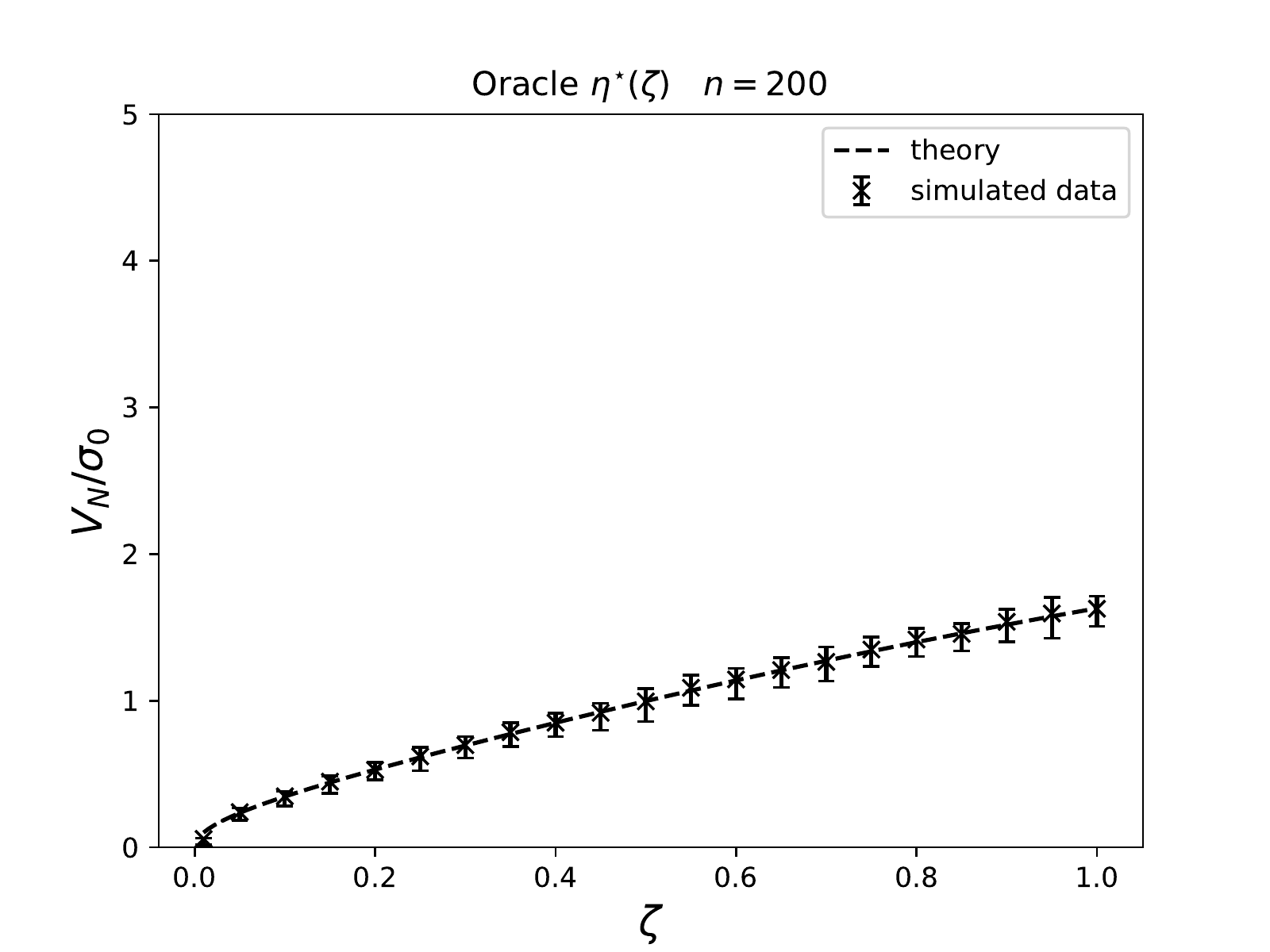}
  \caption{}
  \label{fig:6c}
\end{subfigure}%
\hfill
\begin{subfigure}{.5\textwidth}
  \includegraphics[width=\linewidth]{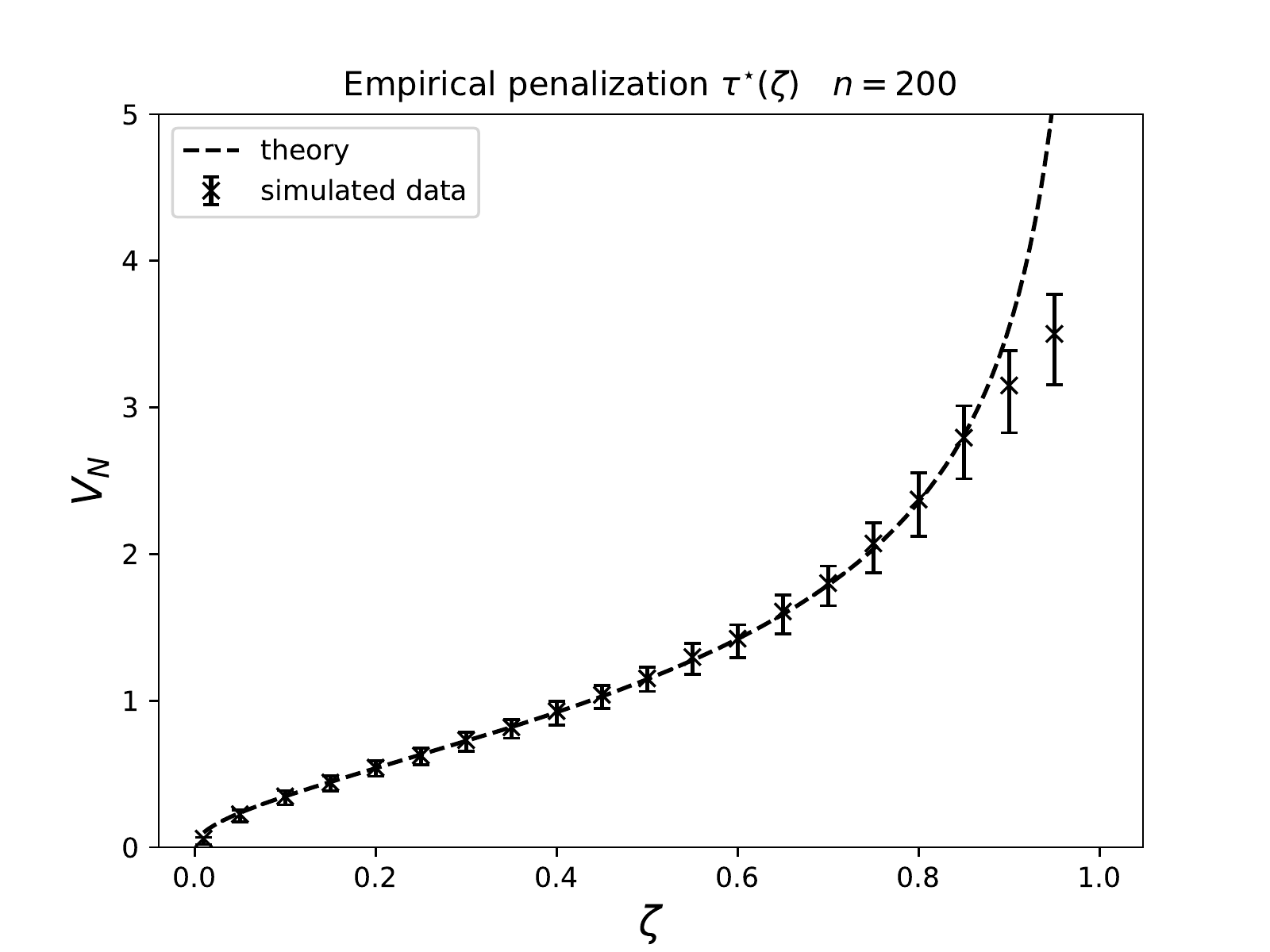}
  \caption{}
  \label{fig:6d}
\end{subfigure}
\caption{Simulation data for the Weibull proportional hazards model along the zero bias line (for $n$=200, $S=1$), with optimal oracle (left) or empirical (right) covariant penalization. The penalization parameters now depend on  $\zeta$ in such a way that the asymptotic bias is zero, i.e.  $w_{\star}/S=1.0$. This is confirmed in Figures  \ref{fig:6a} and \ref{fig:6b}. The estimator fluctuations are shown in Figures \ref{fig:6c} and \ref{fig:6d}. The error bars represent the first and third sample quartiles. Dashed lines give the theoretical predictions.
}
\label{fig:zero_bias_weib}
\end{figure}
We show the result of using the protocols  $\tau^\star(\zeta,S)$ and $\eta^\star(\zeta,S)$, at fixed $S=1$, in Figure \ref{fig:zero_bias_weib}. We see that $K_n$ now fluctuates around the value one, as desired, indicating unbiased inference. 
\begin{figure}[t]
\begin{subfigure}{.5\textwidth}
  \includegraphics[width=\linewidth]{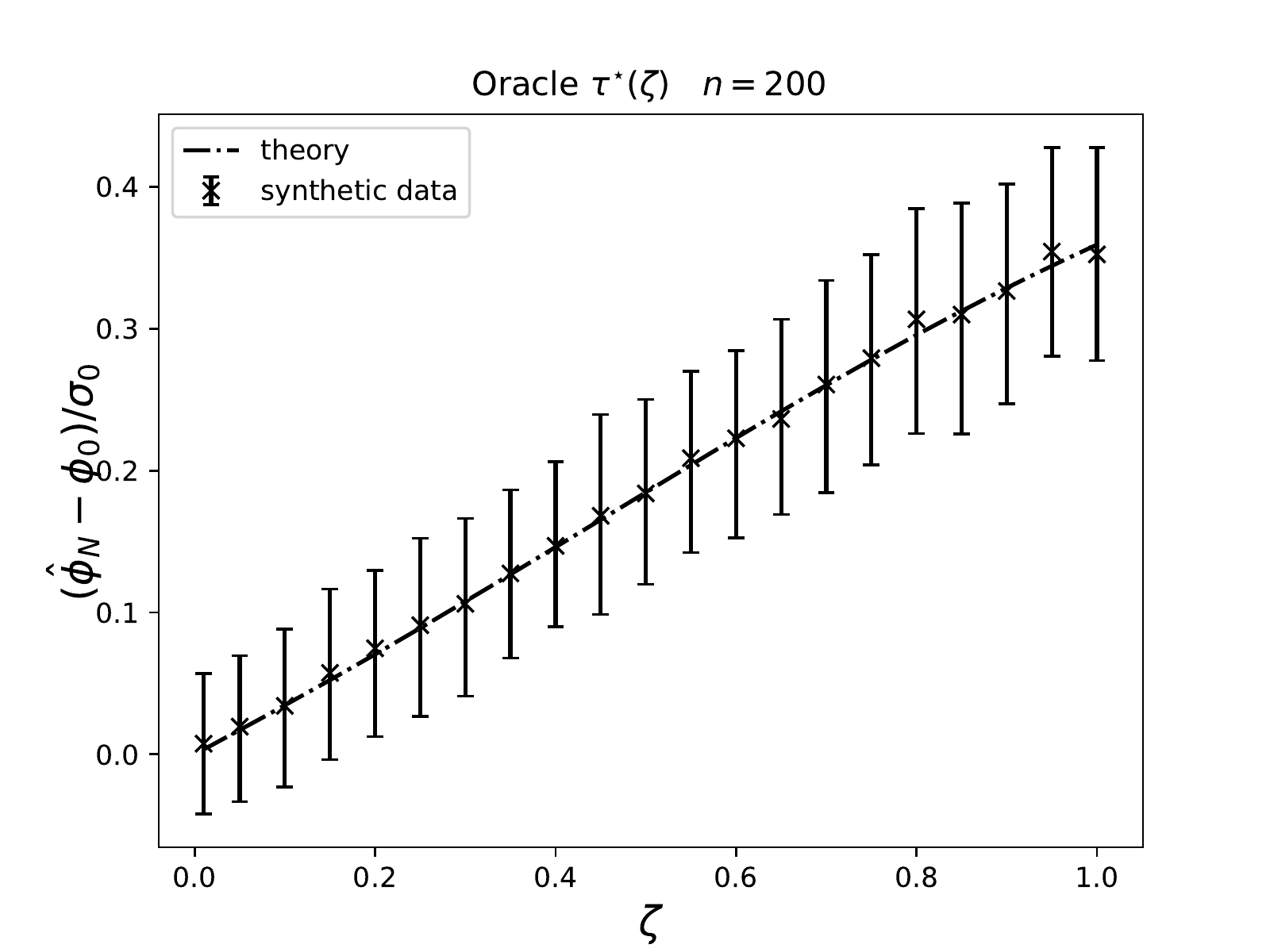}
  \caption{}
  \label{fig:7a}
\end{subfigure}%
\hfill
\begin{subfigure}{.5\textwidth}
  \includegraphics[width=\linewidth]{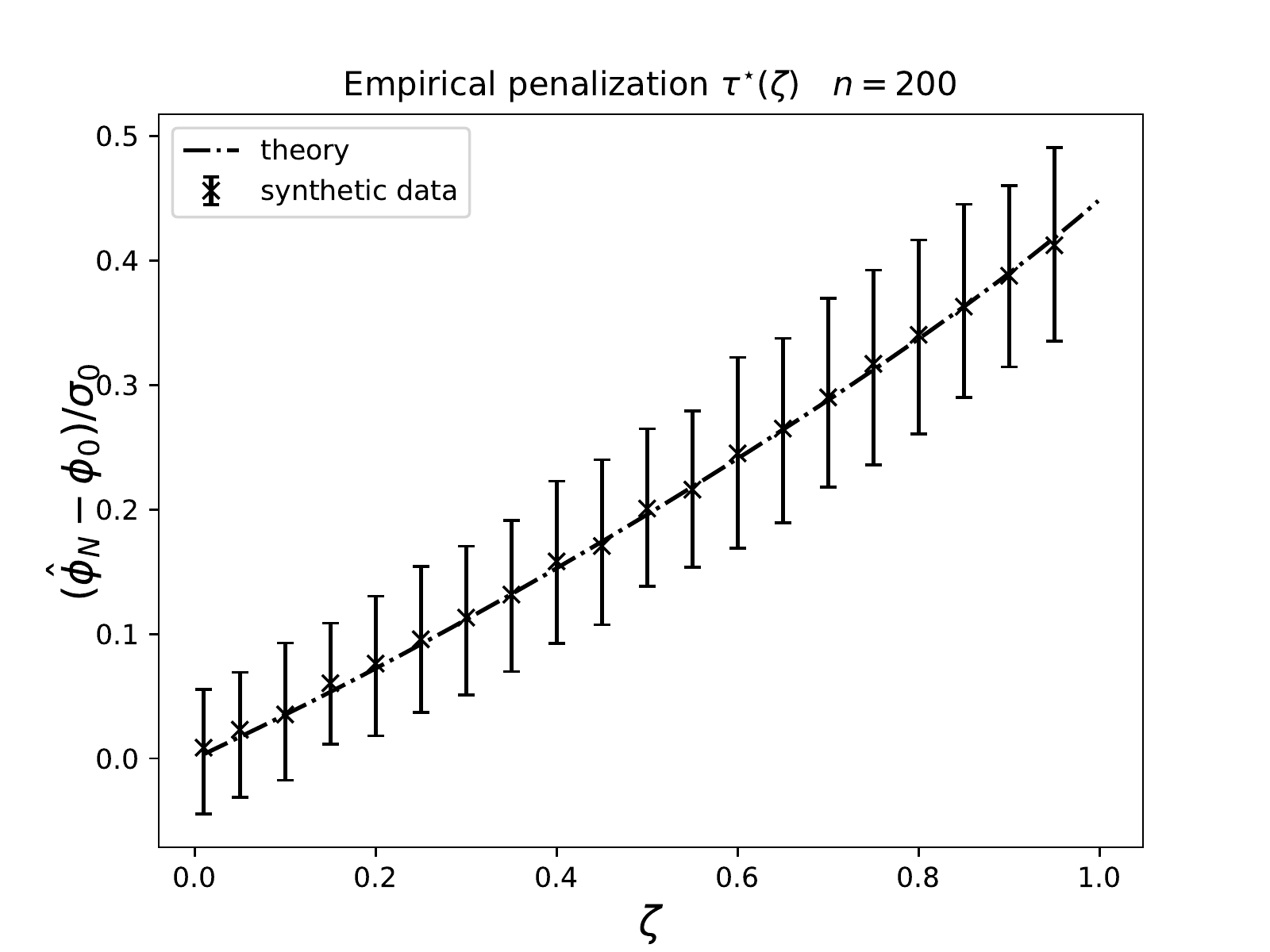}
  \caption{}
  \label{fig:7b}
\end{subfigure}
\medskip 
\begin{subfigure}{.5\textwidth}
  \includegraphics[width=\linewidth]{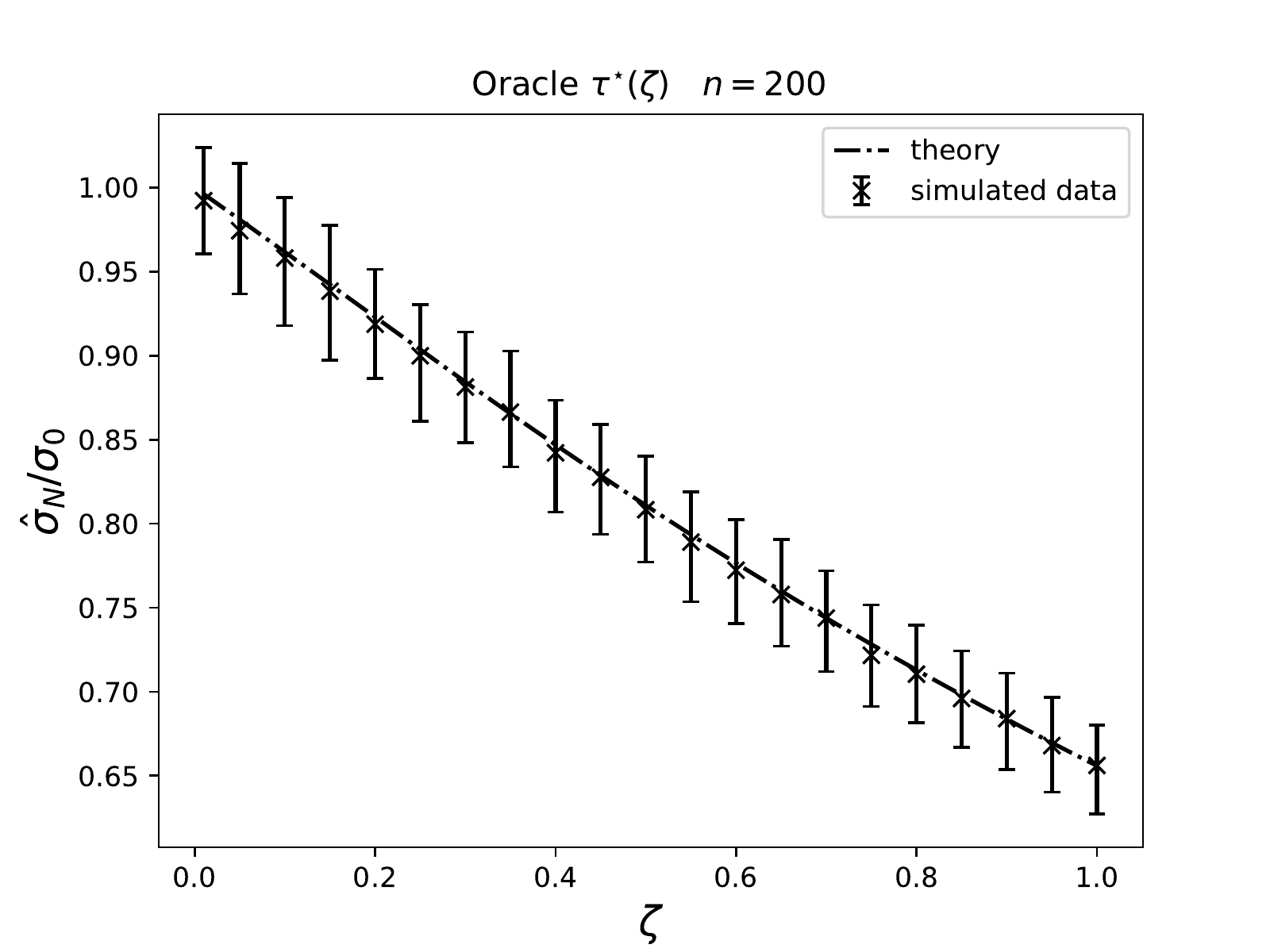}
  \caption{}
  \label{fig:7c}
\end{subfigure}%
\hfill
\begin{subfigure}{.5\textwidth}
  \includegraphics[width=\linewidth]{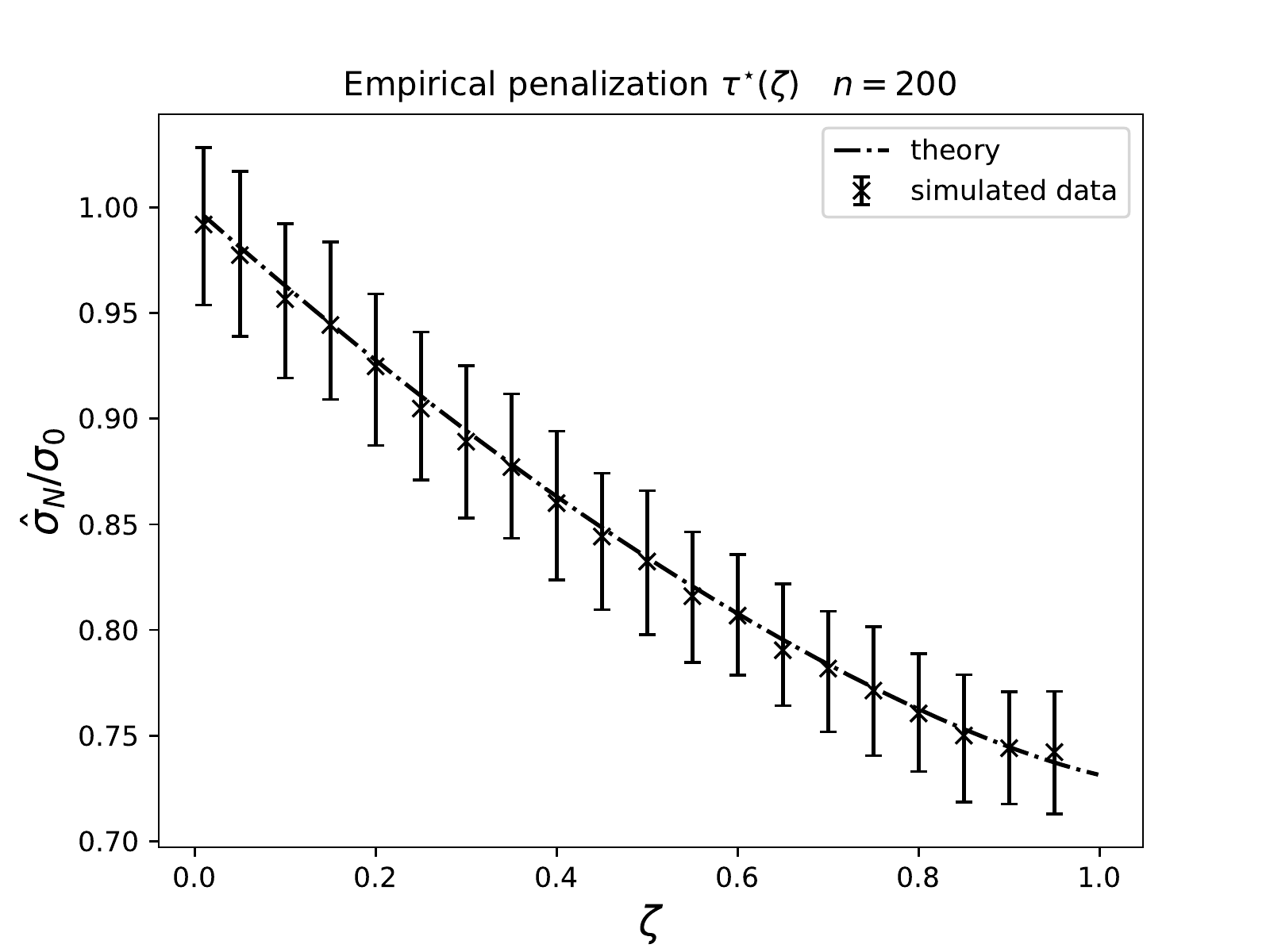}
  \caption{}
  \label{fig:7d}
\end{subfigure}
\caption{ Simulated data for the Weibull proportional hazards model ($n$=200, $S=1$). The error bars represents the first and third sample quartiles, while the dashed line is the curve obtained by solving the RS equations. The penalization parameters are obtained by solving the RS equations, upon requiring $\hat{\bbeta}_n$  be unbiased. It is evident that the nuisance parameters are biased and that the quantities displayed fluctuate closely around the solution of RS equations.}
\label{fig:nuis_weib}
\end{figure}

\begin{figure}[t]
\begin{subfigure}{.5\textwidth}
  \includegraphics[width=\linewidth]{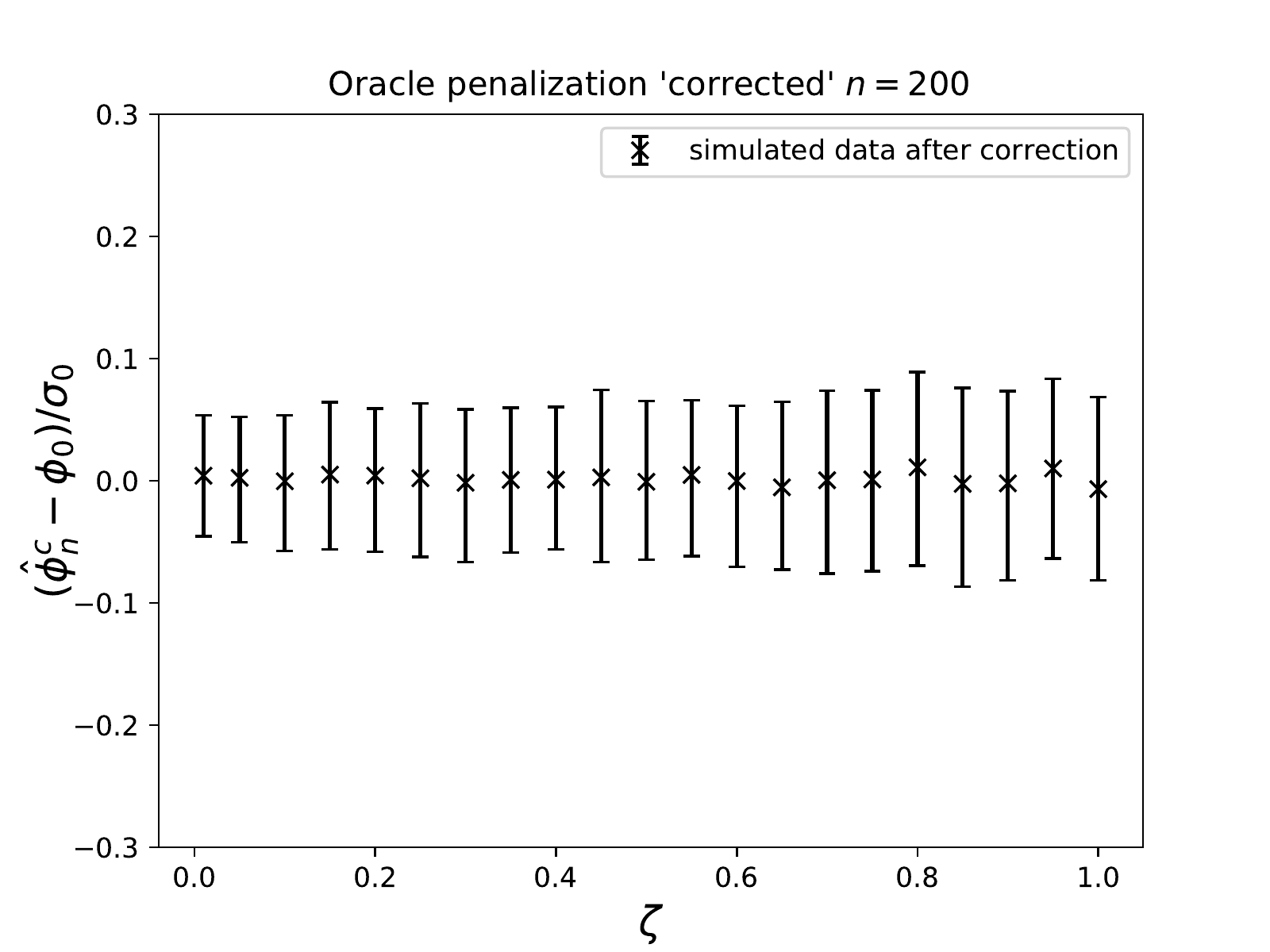}
  \caption{}
  \label{fig:8a}
\end{subfigure}%
\hfill
\begin{subfigure}{.5\textwidth}
  \includegraphics[width=\linewidth]{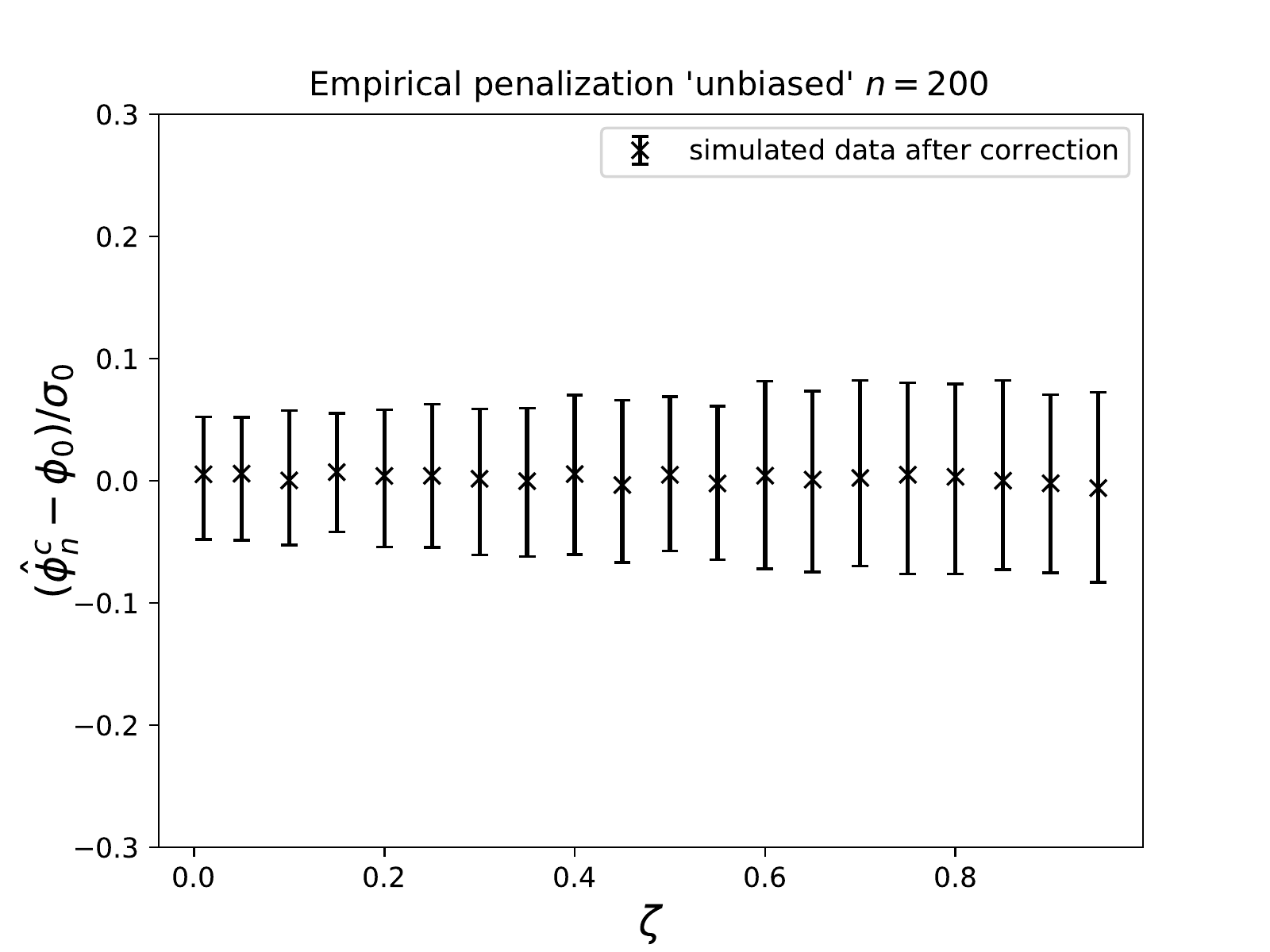}
  \caption{}
  \label{fig:8b}
\end{subfigure}
\medskip 
\begin{subfigure}{.5\textwidth}
  \includegraphics[width=\linewidth]{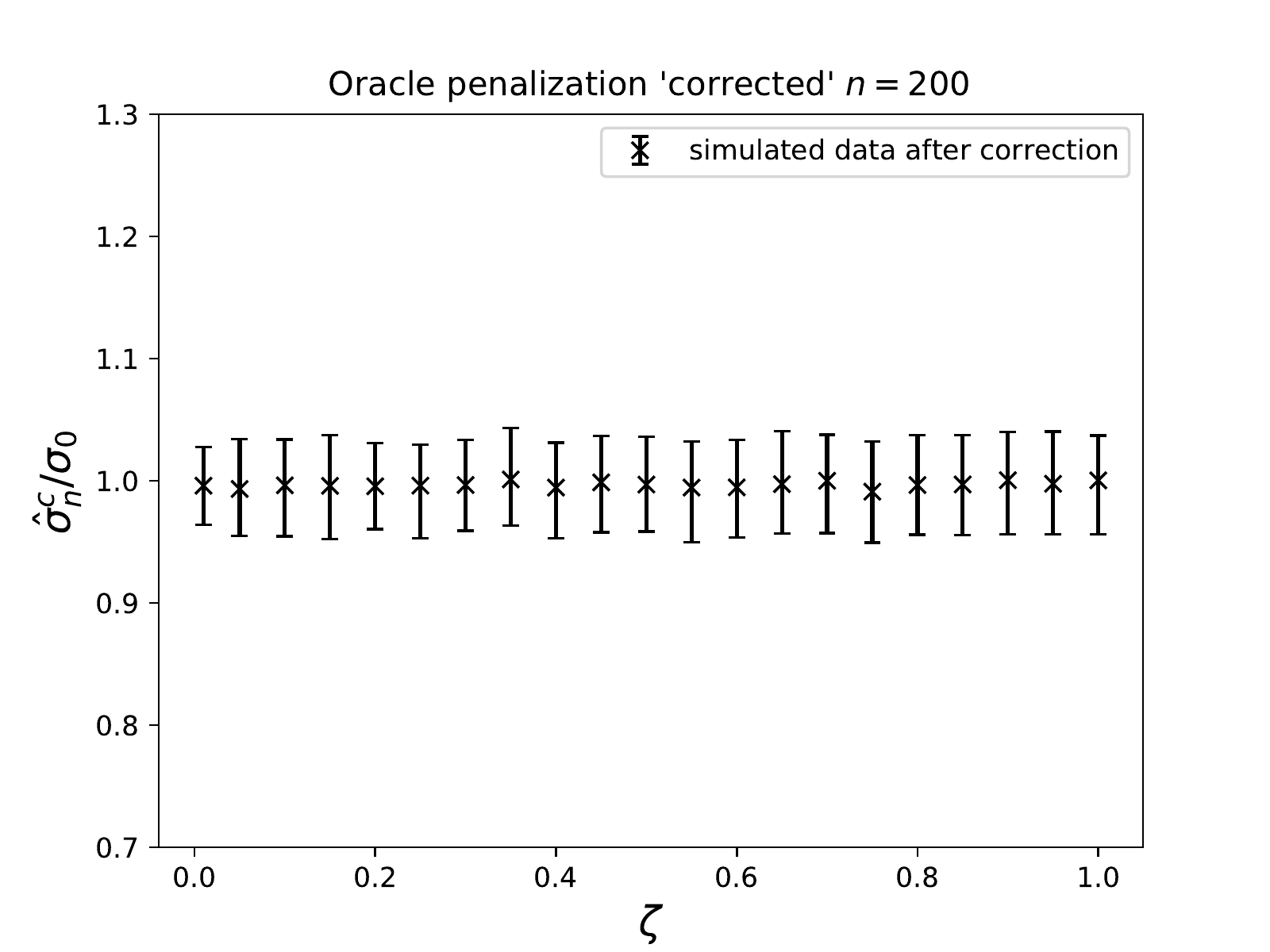}
  \caption{}
  \label{fig:8c}
\end{subfigure}%
\hfill
\begin{subfigure}{.5\textwidth}
  \includegraphics[width=\linewidth]{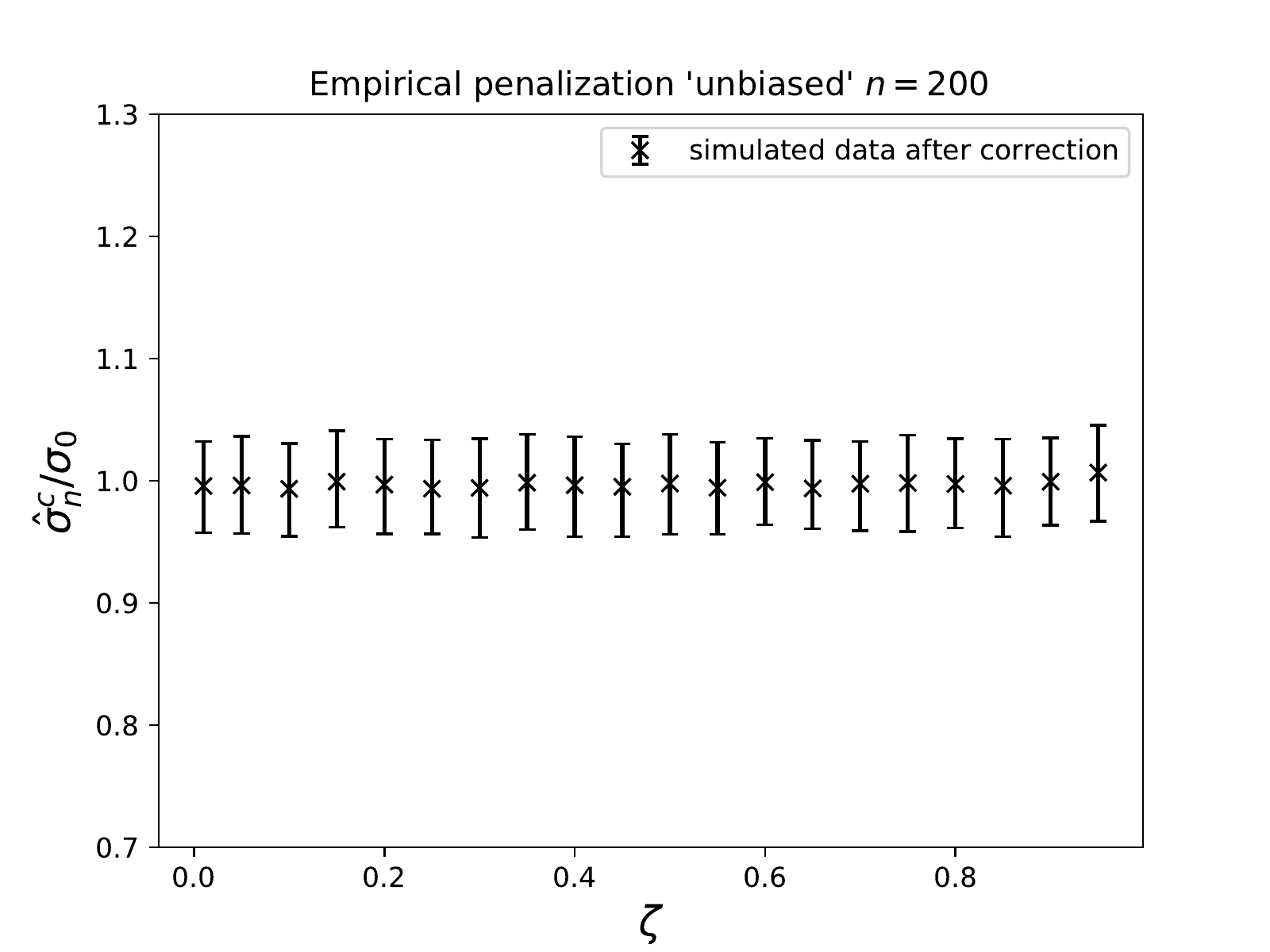}
  \caption{}
  \label{fig:8d}
\end{subfigure}
\caption{ Simulated data for the Weibull proportional hazards model ($n=200$, $S=1$) after application of regularization parameters computed with the RS equations upon demanding unbiased nuisance parameter estimators. Error bars represent the first and third sample quartiles, while the dashed line is the theoretical prediction. Figures \ref{fig:8a} and \ref{fig:8b} show that $(\hat{\phi}^c_n-\phi_0)/\hat{\sigma}_0$ fluctuates around zero, while Figures \ref{fig:8c} and \ref{fig:8d} show that $\hat{\sigma}_n^c/\sigma_0$ fluctuates around one. Hence the estimators $\hat{\phi}^c_n$ and $\hat{\sigma}_n^c$ are indeed {unbiased}.}
\label{fig:nobias_nuis_weib}
\end{figure}

In contrast to the Logit model, the Weibull model has  further nuisance parameters, namely $\phi$ and $\sigma$. It should be emphasised that, while $\hat{\bbeta}_n$ is {unbiased}, the nuisance parameters estimators are not, as one can see in Figure \ref{fig:nuis_weib}.  While in inference the focus is usually on the regression parameters, it might also be of interest to have unbiased estimators for the nuisance parameters; for instance in probabilistic prediction. We could in principle use our equations to predict for which $\tau'$ or $\eta'$  the nuisance parameters estimators will be unbiased, but a more efficient general alternative is to use the solution of the RS equation to compute a {de-biasing factor}, similar to the one employed in \cite{massa}.
In the present case the RS equations can be rewritten in terms of $(\phi-\phi_0)/\sigma_0$ and $\sigma/\sigma_0$. The latter combinations turn out to depend only on $\zeta$, $\tau'$ or $\eta'$, and $S$. Since the value of the estimator will be close to the solution of the RS equations, we may write
\begin{eqnarray}
    &&\frac{\hat{\phi}_n-\phi_0}{
    \sigma_0} = h(\zeta,S) + o(1),~~~~~~
    \frac{\hat{\sigma}_n}{\sigma_0} = g(\zeta,S) + o(1)
\end{eqnarray}
where $h(.,.)$ and $g(.,.)$ are expressions obtained by solving the RS equations. One subsequently solves these two equations for $\phi_0$ and $\sigma_0$, leading to new {unbiased} estimators $\hat{\phi}^c_n, \hat{\sigma}^c_n$. In Figure \ref{fig:nobias_nuis_weib} it can be observed that this procedure indeed leads to estimators for the nuisance parameters that are {unbiased}.

\section{Discussion and conclusion}
\label{section:conclusion}

Ridge regression, a popular inference approach in the high dimensional regime where $\zeta\equiv p/n=O(1)$, leads in generalized linear models typically to rotated estimators of the association parameters $\bbeta\in{\rm I\!R}^p$. This prompted us to search for  generalizations of the ridge penalty, for which the corresponding estimator will be unbiased. We focused on  two elliptical generalizations, of the form $\bbeta\cdot\mathbf{A}\bbeta$. An `oracle' penalization, defined by $\mathbf{A}=\mathbf{A}_0$, uses the population covariance matrix $\mathbf{A}_0$ of the covariates (information that is not always available). An `empirical' penalization uses instead the sample covariance matrix $\mathbf{A}=\sum_{i=1}^n \mathbf{X}_i\mathbf{X}_i/n$.  Both are covariant under linear transformations of the covariates. 
We derive an explicit expression for the distribution of the PML/MAP estimator $\hat{\bbeta}_n$, obtained by optimizing the corresponding objective functions, under the assumption of Gaussian covariates. This estimator is for both covariant priors shown to be typically aligned with the true parameter vector $\bbeta_0$, and its distribution depends solely on two scalar quantities $R^{(n)}_{0,1}$ and $R^{(n)}_{1,1}$, and on $\mathbf{A}_0$. The quantities $R^{(n)}_{0,1}$ and $R^{(n)}_{1,1}$ are finite $n$ equivalents of the replica symmetric (RS) order parameters of  \cite{massa,GLM,PH}. In the high dimensional regime they fluctuate weakly around deterministic values. Under the assumption that they are asymptotically  self-averaging, we show how their values can be computed from a small system of self-consistent equations, the RS equations. Deriving the RS equations alternatively via the cavity method clarifies the relevant underlying assumptions. The RS equations enable one to control the bias-variance trade-off; one can select the value of the regularization parameter by fixing the value of bias that is deemed to be acceptable for the problem at hand, and compute the RS curves to predict the associated variance. Our equations can also be used to give explicit expressions for covariant regularization parameters such that not only the direction but also the length of the estimator is correct, i.e. such that the estimators will be strictly unbiased. 

We tested our theoretical results against the parameter statistics computed from inferences with simulated data, and  found excellent agreement. We worked out the theory for two popular generalized linear models: teh Logit model for binary classification, and the Weibull proportional hazards model for survival analysis. For both applications,  the  solution of the RS equations correctly predicts the behaviour of the estimators obtained from simulated data, even for relatively small sample sizes.  
The PML/MAP estimators obtained with the oracle penalization achieve a lower variance than the ones computed using the empirical penalty, which can be explained in terms of the bias-variance trade-off.
When $\zeta$ is close to one, i.e. deep into the MLE overfitting regime, the sample covariance matrix may develop zero eigenvalues, leading to pathological estimators. Sufficiently below $\zeta=1$ zero eigenvalues will not occur, and the {empirical} penalization is an effective, practical and easy to implement method for achieving unbiased estimators. Oracle penalization obviously does not have this problem, and can  be used even for $\zeta>1$ (provided $\mathbf{A}_0$ has full rank). This also suggests an alternative third route, where one estimates the population covariance matrix based on previous data, with methods such as \cite{Ledoit}, to be used as a proxy for $\mathbf{A}_0$ in the oracle penalization (a possible subject of future investigation).

While the RS equations can be solved for any given value of $S=|\bbeta_0|$, in practice one does not know this value. For a significant class of models $S$ is in fact found to drop out of the RS equations, and for those models where this does not happen there are transparent approaches to estimate $S$. Still this is a subject hat warrants further research. 
It might also be interesting  to analyse what happens when the covariant regularizers are used in a model mis-match setting (see e.g. \cite{barbier}), which in principle can be addressed in the present formulation. 

Our results are derived under the assumption of Gaussian distributed covariates. The results concerning the first two moments of the estimator $\hat{\bbeta}_n$ should however hold more generally, as suggested by the work of e.g. 
\cite{el_karoui_rigorous,PH,sheik,GLM,massa}. The underlying reason is that even for non-Gaussian covariates, the quantities  $\mathbf{X}\cdot\bbeta$ may asymtotically have Gaussian statistics  due to the Central Limit Theorem (under certain conditions). See also \cite{goldt,gaussian_mixtures,loureiro}.  The universality of the asymptotic behaviour of the objective function has been proven for the machine learning setting in \cite{montanari22universality}. It would be interesting to study whether in arbitrary generalized linear models this universality would extend to the asymptotic properties of the estimators, at least concerning the first two moments. 

To conclude, we present covariant generalizations of ridge penalization in high-dimensional regression, such that the parameter estimators no longer exhibit the undesirable rotation that plagues conventional ridge regularization (especially for correlated data). For arbitrary generalized linear models the corresponding inference process is  described accurately by an RS theory (derived alternatively via the cavity method), even for modest sample sizes. This theory can be used e.g. to tune the bias-variance trade-off, or to predict which values of the covariant regularization parameters should be used in practice to obtain strictly unbiased estimators.

\appendix 

\section{Representation of the PML/MAP estimator}
\label{appendix:representation}

The PML/MAP estimator for the regression parameters is defined by the following optimization problem
\begin{eqnarray}
    \hat{\bbeta}^{\sim}_n &=&  \underset{\bbeta}{\arg\max}\bigg( \underset{\bsigma}{\max}\Big\{ \sum_{i=1}^n \log p (T_i|\bm{\mathcal{X}}_i\cdot\bm{\beta},\bm{\sigma}) -\frac{p}{2} \bbeta\cdot\Big( \tau'  \frac{1}{n}\sum_{i=1}^n \bm{\mathcal{X}}_i\bm{\mathcal{X}}_i +\eta' \bm{I}\Big)\bbeta \Big\}\bigg)
\end{eqnarray}
where
\begin{equation}
    T_i|\bm{\mathcal{X}}_i \sim p(T_i|\bm{\mathcal{X}}_i\cdot\bbeta^{\sim}_0,\bsigma_0)
\end{equation}
Now consider any rotation $\mathbf{R}_0$ in $\mathbb{R}^p$ around $\bbeta_0$. From property (\ref{property_ml}) in the main text and the fact that the conditional distribution of the response $T_i$ depends only on the projection $\bm{\mathcal{X}}_i\cdot \bbeta^{\sim}_0$, 
\begin{equation}
    T_i|\bm{\mathcal{X}}_i  = T_i|\bm{\mathcal{X}}_i\cdot \bbeta^{\sim}_0
\end{equation}
we obtain
\begin{equation}
\label{invariance}
    \hat{\bbeta}^{\sim}_n=\hat{\bbeta}_n (\{T_i,\bm{\mathcal{X}}_i\},\bbeta_0^{\sim}) = \mathbf{R}_0 \hat{\bbeta}_n (\{T_i,\mathbf{R}_0\bm{\mathcal{X}}_i\},\bbeta_0^{\sim})\overset{d}{=} \mathbf{R}_0 \hat{\bbeta}_n(\{T_i,\bm{\mathcal{X}}_i\},\bbeta_0^{\sim})
\end{equation}
where we have denoted equality in distribution with $\overset{d}{=}$. The first equality follows because $\bbeta_0^{\sim}= \mathbf{R}_0\bbeta_0^{\sim}$,  and the last equality follows because $\mathbf{O}\bm{\mathcal{X}}_i$ has the same distribution as $\bm{\mathcal{X}}_i$ for any rotation $\mathbf{O}$ since $\bm{\mathcal{X}}_i\sim\mathcal{N}(\bm{0},\bm{I})$ (hence also for any $\mathbf{R}_0$). 
To understand the implications of (\ref{invariance}), note that we can always decompose $\hat{\bbeta}^{\sim}_n$ into a component along $\bbeta^{\sim}_0$ and a component in the subspace of $\mathbb{R}^p$  orthogonal to $\bbeta^{\sim}_0$
\begin{equation}
    \hat{\bbeta}^{\sim}_n =  \hat{\bbeta}^{\sim}_{n,\|} +  \hat{\bbeta}^{\sim}_{n,\perp} \qquad \hat{\bbeta}^{\sim}_{n,\|}\cdot \hat{\bbeta}^{\sim}_{n,\perp}=0
\end{equation}
with 
\begin{equation}
    \hat{\bbeta}^{\sim}_{n,\|}:=\bbeta_0^{\sim}(\bbeta_0^{\sim}\cdot\hat{\bbeta}^{\sim}_n)/\|\bbeta_0^{\sim}\|^2 \qquad \hat{\bbeta}^{\sim}_{n,\perp} := \Big(\bm{I}-\bbeta_0^{\sim}\bbeta_0^{\sim}/\|\bbeta_0^{\sim}\|^2\Big)\hat{\bbeta}^{\sim}_n
\end{equation}
Then equation (\ref{invariance}) implies that 
\begin{equation}
    \mathbf{R}_0\hat{\bbeta}^{\sim}_{n,\perp} \overset{d}{=} \hat{\bbeta}^{\sim}_{n,\perp}
\end{equation}
because $\hat{\bbeta}^{\sim}_{n,\|}$ is not affected by any $\mathbf{R}_0$ as it is parallel to $\bbeta_0^{\sim}$. Hence all the values of $\hat{\bbeta}^{\sim}_{n,\perp}$ that have the same length, i.e.\ that lie on a sphere in the subspace of $\mathbb{R}^p$ orthogonal to $\bbeta_0^{\sim}$, must have the same probability density.
In turn this means that the length of $\hat{\bbeta}^{\sim}_{n,\perp}$ must be independent from its direction: conditional on the length, the direction is uniformly distributed over a sphere in the above mentioned subspace.
This means that, conditional on $\hat{\bbeta}^{\sim}_{n,\|}$, we have the following representation for $\hat{\bbeta}^{\sim}_{n,\perp}$
\begin{equation}
    \hat{\bbeta}^{\sim}_{n,\perp} = \|\hat{\bbeta}^{\sim}_{n,\perp}\|  \frac{\hat{\bbeta}^{\sim}_{n,\perp} }{\|\hat{\bbeta}^{\sim}_{n,\perp}\|} = \|\hat{\bbeta}^{\sim}_{n,\perp}\| \mathbf{U} 
\end{equation}
with $\mathbf{U}$ uniformly distributed on the unit sphere in the $p-1$ dimensional subspace $\mathcal{S}_{p-2}$ of $\mathbb{R}^p$ that is orthogonal to $\bbeta_0^{\sim}$.
Now note that 
\begin{equation}
    \|\hat{\bbeta}^{\sim}_{n,\perp}\| =  \sqrt{\|\hat{\bbeta}^{\sim}_n\|^2-\|\hat{\bbeta}^{\sim}_{n,\|}\|^2} = \sqrt{\|\hat{\bbeta}_n^{\sim}\|^2-(\bbeta_0^{\sim}\cdot\hat{\bbeta}_n^{\sim})^2/\|\bbeta_0^{\sim}\|^2}
\end{equation}
So finally we obtain the representation
\begin{equation}
    \hat{\bbeta}^{\sim}_n = \frac{\bbeta^{\sim}_0\cdot\hat{\bbeta}^{\sim}_n}{\|\bbeta^{\sim}_0\|^2}\bbeta^{\sim}_0 + \sqrt{\|\hat{\bbeta}^{\sim}_n\|^2-\frac{(\bbeta_0^{\sim}\cdot\hat{\bbeta}^{\sim}_n)^2}{\|\bbeta^{\sim}_0\|^2}}\mathbf{U}
\end{equation}

\section{The value of $\alpha^2$}
\label{appendix:alpha}

To compute the limit defining $\alpha^2$ in (\ref{alpha_def}), we note that 
\begin{eqnarray}
\label{exp_a_p}
     \mathbb{E}[A^2_p]&=& \mathbb{E}\Big[\mathbf{U}\cdot\mathbf{A}_0^{-1}\mathbf{U}\Big] = \mathbb{E}\Big[\frac{\mathbf{Z}\cdot\mathbf{P}\mathbf{A}_0^{-1}\mathbf{P}\mathbf{Z}}{\mathbf{Z}\cdot\mathbf{P}\mathbf{Z}}\Big] =\nonumber \\ &=&\frac{1}{p-1}\Big(\Tr (\mathbf{A}^{-1}_0) - \frac{\bbeta_0^{\sim}\cdot \mathbf{A}_0^{-1}\bbeta_0^{\sim}}{\|\bbeta_0^{\sim}\|^2}\Big) + o_P(1)
\end{eqnarray}
In (\ref{exp_a_p}), to obtain the third equality we used the fact that 
\begin{equation}
    \mathbf{U} \overset{d}{=} \frac{\mathbf{P}\mathbf{Z}}{\|\mathbf{P}\mathbf{Z}\|}  \qquad \mathbf{Z}\sim\mathcal{N}(\bm{0},\bm{I}), \quad \mathbf{P}:= \bm{I}-\frac{\bbeta^{\sim}_0\bbeta^{\sim}_0\cdot}{\|\bbeta^{\sim}_0\|^2}
\end{equation}
For the last equality we used the fact that quadratic forms like $\mathbf{Z}\cdot\mathbf{M}\mathbf{Z}/p$ concentrate around their expected value $\Tr(\mathbf{M})$ in the asymptotic high dimensional limit \cite{el_karoui1,Boucheron13}
\begin{equation}
    \mathbf{Z}\cdot\mathbf{M}\mathbf{Z}/p = \Tr(\mathbf{M}) + o_P(1)
\end{equation}
At this point we see that by definition of $\bbeta_0^{\sim}:= \mathbf{A}_0^{1/2}\bbeta_0$
\begin{equation}
    \frac{\bbeta_0^{\sim}\cdot \mathbf{A}_0^{-1}\bbeta_0^{\sim}}{\|\bbeta_0^{\sim}\|^2} = \frac{\bbeta_0\cdot \bbeta_0}{\bbeta_0\cdot \mathbf{A}_0\bbeta_0}
\end{equation}
and 
\begin{equation}
   \lambda_{\max}^{-1}(\mathbf{A}_0) \leq \frac{\bbeta_0\cdot \bbeta_0}{\bbeta_0\cdot \mathbf{A}_0\bbeta_0}\leq \lambda_{\min}^{-1}(\mathbf{A}_0)
\end{equation}
thus 
\begin{equation}
    \lambda_{\max}^{-1}(\mathbf{A}_0)/(p-1)\leq\Big|\alpha^2_p - \Tr(\mathbf{A}_0^{-1})\big/(p-1)\Big|\leq \lambda_{\min}^{-1}(\mathbf{A}_0)/(p-1)
\end{equation}
If $\lambda_{min}(\mathbf{A}_0)=O(p^{\alpha})$, with $\alpha>-1$, then $\lim_{p\rightarrow \infty}\lambda_{\min}^{-1}(\mathbf{A}_0)/(p-1)=0$.
Hence we conclude that
\begin{equation}
    \alpha^2 := \lim_{p\rightarrow \infty} \mathbb{E}[A_p^2] = \lim_{p\rightarrow \infty} \frac{1}{p}\Tr(\mathbf{A}_0^{-1})
\end{equation}

\section{Derivation of the RS equations with the Cavity method}
\label{appendix:derivation}

In this section we derive the Replica symmetric equations (\ref{RS1}, \ref{RS2}, \ref{RS3}, \ref{RS4}).
We first establish the asymptotic properties of the PML/MAP estimator of the regression parameters $\hat{\bbeta}_n$ for a fixed value of the nuisance parameters $\bsigma$. We subsequently show  that, given our assumptions and approximations, this is sufficient to establish also the asymptotic behaviour of $\hat{\bsigma}_n$.

\subsection{Analysis via statistical physics}

We exploit the analogy between optimization and statistical physics, whereby  variables over which an optimization is carried out are interpreted as degrees of freedom of {particles}, with a Hamiltonian that is minus the objective function to be optimized. The solution of the optimization problem then equals the ground state of the Hamiltonian. Here we aim to understand the properties of the PML/MAP estimator obtained with $\{T_i,\bm{\mathcal{X}}_i\}_{i=1}^n$ where $\bm{\mathcal{X}}_i\sim \mathcal{N}(\bm{0},\bm{I})$, $T_i|\bm{\mathcal{X}}_i \sim p(T_i|\bm{\mathcal{X}}_i\cdot\bbeta_0^{\sim},\bsigma_0)$ and $\bbeta_0^{\sim}:= \mathbf{A}_0^{1/2}\bbeta_0$. 
So we consider the Hamiltonian
\begin{equation}
    \mathcal{H}_{n}(\bm{\beta},\bm{\sigma}) := -l_n(\bbeta,\bsigma) = \frac{1}{2} \eta \|\bm{\beta}\|^2-\sum_{i=1}^n \rho (T_i|\bm{\mathcal{X}}_i\cdot\bm{\beta},\bm{\sigma}) 
\end{equation}
The statistical properties of the physical system at temperature $1/\gamma$ follow from the Gibbs measure
\begin{equation}
    \bbeta|\{T_i,\bm{\mathcal{X}}_i\}_{i=1}^n \sim p(\bm{\beta})=\rme^{-\gamma \mathcal{H}_{n}(\bm{\beta},\bm{\sigma})}/Z_{n}(\gamma)
\end{equation}
where the partition function ensures the normalization, i.e. 
\begin{equation}
    Z_{n}(\gamma)= \int\!\rmd \bbeta~ \rme^{-\gamma \mathcal{H}_{n}(\bm{\beta},\bm{\sigma})} 
\end{equation}
Under regularity conditions on the Hamiltonian function, and for a well behaved functions $f$, the Laplace argument implies that 
\begin{equation}
    f(\hat{\bbeta}^{\sim}_n) = \lim_{\gamma \rightarrow \infty} \int\!\rmd \bbeta~ f(\bbeta)\frac{\rme^{-\gamma \mathcal{H}_{n}(\bm{\beta},\bm{\sigma})}}{Z_{n}(\gamma)}  = \lim_{\gamma \rightarrow \infty} \langle f(\bbeta) \rangle
\end{equation}
where
\begin{equation}
    \hat{\bbeta}^{\sim}_n := \underset{\bbeta}{\arg \min} \Big\{ \mathcal{H}_n(\bm{\beta},\bm{\sigma})\Big\} = \underset{\bbeta}{\arg \max} \Big\{ l_n(\bm{\beta},\bm{\sigma})\Big\}
\end{equation}
We are interested in the typical properties of (functions of) the PML estimator, hence we study averages over all possible realizations of the data-set:
\begin{equation}
    \mathbb{E}_{T,\bm{\mathcal{X}}}\Big[f(\hat{\bbeta}^{\sim}_n)\Big] = \lim_{\gamma \rightarrow \infty} \mathbb{E}_{T,\bm{\mathcal{X}}}\Big[\langle f(\bbeta)\rangle \Big]
    \label{av_f_thermal}
\end{equation}
We switched the limit and the expectation over the data-set, which is valid under regularity conditions on $f$.
The right hand side of (\ref{av_f_thermal}) has the structure of the computation of a disorder-average observable $f$ in a spin glass model \cite{virasoro,parisi,Mezard_89}, with the role of disorder played by the data set $\{T_i,\bm{\mathcal{X}}_i\}_{i=1}^n$.
This enables the use of tools from the physics of disordered systems.

\subsection{Application of the Cavity method}

 We next use the cavity method to derive the asymptotic properties of the PML/MAP estimator.
Consider a sufficiently well behaved function $f$ that depends only on the linear predictor $\bm{\mathcal{X}}_i\cdot\bm{\beta}$, the response $T_i$ of the  $i$th observation and the (for now fixed) nuisance parameters $\bm{\sigma}$. By noting that
\begin{equation}
    \mathcal{H}_{n,p}(\bm{\beta},\bm{\sigma}) :=\mathcal{H}_{(i),p}(\bm{\beta},\bm{\sigma})- \rho (T_{i}|\bm{\mathcal{X}}_{i}\cdot\bm{\beta},\bm{\sigma}) 
\end{equation}
we can write
\begin{eqnarray}
   \langle f(T_i|\bm{\mathcal{X}}_i\cdot\bm{\beta},\bm{\sigma})\rangle &=& \frac{\int\!\rmd \bbeta~ f(T_i|\bm{\mathcal{X}}_i\cdot\bm{\beta},\bm{\sigma})\rme^{\gamma \rho(T_{i}|\bm{\mathcal{X}}_{i}\cdot\bm{\beta},\bm{\sigma})-\gamma \mathcal{H}_{(i)}(\bm{\beta},\bm{\sigma})} }{\int\!\rmd \bbeta~ \rme^{\gamma \rho(T_{i}|\bm{\mathcal{X}}_{i}\cdot\bm{\beta},\bm{\sigma})-\gamma \mathcal{H}_{(i)}(\bm{\beta},\bm{\sigma})}}
   \nonumber\\
   &=&\frac{\big\langle f(T_i|\bm{\mathcal{X}}_i\cdot\bm{\beta},\bm{\sigma})\exp\Big\{\gamma \rho (T_{i}|\bm{\mathcal{X}}_{i}\cdot\bm{\beta},\bm{\sigma})\Big\}\big\rangle_{(i)}}{\big\langle\exp\Big\{\gamma \rho (T_{i}|\bm{\mathcal{X}}_{i}\cdot\bm{\beta},\bm{\sigma})\Big\}\big\rangle_{(i)}}
   \label{cavity_exp}
\end{eqnarray}
Here $\langle\rangle_{(i)}$ refers to an average under the Gibbs measure with Hamiltonian $\mathcal{H}_{(i),p}(\bm{\beta},\bm{\sigma})$. 
The computation of (\ref{cavity_exp}) is difficult because we need the distribution of the linear predictor $Y_i :=\bm{\mathcal{X}}_i\cdot\bm{\beta}$ under the cavity measure.  The first two moments of  $Y_i$ are 
\begin{eqnarray}
    \langle Y_i\rangle = \bm{\mathcal{X}}_i\cdot\langle\bm{\beta}\rangle_{(i)},
    ~~~~~~~~
    \langle Y_i^2\rangle = \bm{\mathcal{X}}_i\cdot\langle\bm{\beta}\bm{\beta}\rangle_{(i)}\bm{\mathcal{X}}_i
\end{eqnarray}
Now standard results in the theory of concentration of measure \cite{el_karoui1,Boucheron13} tell us that 
\begin{eqnarray}
    \fl \mathcal{U}_{(i)}^2 &=& \langle Y_i^2\rangle-\langle Y_i\rangle^2 = \bm{\mathcal{X}}_i\cdot\Big(\langle\bm{\beta}\bm{\beta}\rangle_{(i)}-\langle\bm{\beta}\rangle_{(i)}\langle\bm{\beta}\rangle_{(i)}\Big)\bm{\mathcal{X}}_i \nonumber\\
    \fl &=& \Tr\Big(\langle\bm{\beta}\bm{\beta}\rangle_{(i)}-\langle\bm{\beta}\rangle_{(i)}\langle\bm{\beta}\rangle_{(i)}\Big) +o_P(1) =\langle\bm{\beta}\cdot\bm{\beta}\rangle_{(i)}-\langle\bm{\beta}\rangle_{(i)}\cdot\langle\bm{\beta}\rangle_{(i)} + o_P(1)
\end{eqnarray}
We note that $\mathcal{U}_{(i)}^2$ is exactly the nonlinear susceptibility under the cavity Gibbs measure. We assume that the cavity distribution of $Y_i$ is  Gaussian, so  
\begin{equation}
    Y_i = \bm{\mathcal{X}}_i\cdot\langle\bm{\beta}\rangle_{(i)} + \mathcal{U}_{(i)} Z_i, \qquad Z_i\sim\mathcal{N}(0,1)
\end{equation}
with $Z_i \perp T_i,\bm{\mathcal{X}}_i$.
Remembering that $T_i|\bm{\mathcal{X}}_i\sim p(T_i|\bm{\mathcal{X}}_i^T\bbeta^{\sim}_0,\bsigma_0)$, so the distribution of the disorder depends on $\bm{\mathcal{X}}_i\cdot\bbeta^{\sim}_0$ which is correlated with $\bm{\mathcal{X}}_i\cdot\langle\bm{\beta}\rangle_{(i)}$,  it is possible to  separate $\bm{\mathcal{X}}_i\cdot\langle\bm{\beta}\rangle_{(i)}$ into two independent contributions
\begin{eqnarray}
    \bm{\mathcal{X}}_i\cdot\langle\bm{\beta}\rangle_{(i)} &=& \bm{\mathcal{X}}_i\cdot\Big(\bm{I}-\frac{\bbeta^{\sim}_0\bbeta^{\sim}_0}{\|\bbeta^{\sim}_0\|^2}\Big)\langle\bm{\beta}\rangle_{(i)} + \frac{\bm{\mathcal{X}}_i\cdot\bbeta^{\sim}_0}{\|\bbeta^{\sim}_0\|} \frac{\bbeta^{\sim}_0\cdot \langle\bm{\beta}\rangle_{(i)}}{\|\bbeta^{\sim}_0\|}  \nonumber\\
    &=& V_{(i)}Q_i + W_{(i)}Z_{0,i}, \qquad Q_i\perp Z_{0,i}, \quad Q_i,Z_{0,i} \sim \mathcal{N}(0,1)
\end{eqnarray}
involving the generalized overlaps under the leave $i$-th out measure 
\begin{eqnarray}
    W_{(i)} &:=& \frac{\bbeta^{\sim}_0\cdot \langle\bm{\beta}\rangle_{(i)}}{\|\bbeta^{\sim}_0\|}\\
    V_{(i)} &:=& \Big\|\Big(\bm{I}-\frac{\bbeta^{\sim}_0\bbeta^{\sim}_0}{\|\bbeta^{\sim}_0\|^2}\Big)\langle\bm{\beta}\rangle_{(i)}\|^2 = \|\langle\bm{\beta}\rangle_{(i)}\|^2 - W_{(i)}^2
\end{eqnarray}
These are independent of $T_i,\bm{\mathcal{X}}_i$, hence
\begin{equation}
    Y_i =    W_{(i)}Z_{0,i}+V_{(i)}Q_i + \mathcal{U}_{(i)} Z_i
\end{equation}
and 
\begin{equation}
   \fl \langle f(T_i|\bm{\mathcal{X}}_i\cdot\bm{\beta},\bm{\sigma})\rangle = \frac{\mathbf{E}_{Z_i}\Big[f(T_i| W_{(i)}Z_{0,i}+V_{(i)}Q_i + \mathcal{U}_{(i)} Z_i,\bm{\sigma})\rme^{\gamma \rho (T_{i}| W_{(i)}Z_{0,i}+V_{(i)}Q_i + U_{(i)} Z_i,\bm{\sigma})}\Big]}{\mathbf{E}_{Z_i}\Big[\rme^{\gamma \rho (T_{i}| W_{(i)}Z_{0,i}+V_{(i)}Q_i + \mathcal{U}_{(i)} Z_i,\bm{\sigma})}\Big]}
\end{equation}
This can be rewritten via a simple change of variable, $\xi:= W_{(i)}Z_{0,i}+V_{(i)}Q_i + \mathcal{U}_{(i)} Z_i$, as
\begin{eqnarray}
    \fl \langle f(T_i|\bm{\mathcal{X}}_i\cdot\bm{\beta},\bm{\sigma})\rangle &=& \mathbf{E}_{\xi_i}\Big[f(T_i,\xi_i,\bsigma)\Big]\nonumber\\
    \fl &=&\frac{\int\!\rmd \xi_i~ f(T_i, \xi,\bm{\sigma})\exp\Big\{-\frac{1}{2} \Big(\frac{\xi_i-W_{(i)}Z_{0,i}+V_{(i)}Q_i}{\mathcal{U}_{(i)}}\Big)^2 +\gamma \rho (T_{i}, \xi,\bm{\sigma})\Big\}  }{\int\!\rmd \xi_i~ \exp\Big\{-\frac{1}{2} \Big(\frac{\xi_i-W_{(i)}Z_{0,i}+V_{(i)}Q_i}{\mathcal{U}_{(i)}}\Big)^2 +\gamma \rho (T_{i}, \xi,\bm{\sigma})\Big\}}
   \label{cavity_exp_2}
\end{eqnarray}
The disorder is now {contained} only in the leave-one-out generalized overlaps $W_{(i)},V_{(i)}, \mathcal{U}_{(i)}$. 
It seems reasonable to approximate the overlaps computed under the full Gibbs measure by their leave-one-out approximations for each $i$, 
\begin{eqnarray}
    W_n &:=& \frac{\bbeta^{\sim}_0\cdot\langle\bbeta\rangle}{\|\bbeta^{\sim}_0\|^2}= W_{(i)} + o_P(1)\\
    V_{n} &:=& \Big\|\Big(\bm{I}-\frac{\bbeta^{\sim}_0\bbeta^{\sim}_0}{\|\bbeta^{\sim}_0\|^2}\Big)\langle\bm{\beta}\rangle\|^2= V_{(i)} + o_P(1)\\[2mm]
    \mathcal{U}_n &:=& \langle\bm{\beta}\cdot\bm{\beta}\rangle-\langle\bm{\beta}\rangle\cdot\langle\bm{\beta}\rangle = \mathcal{U}_{(i)} + o_P(1)
\end{eqnarray}
 because, in the limit of large $n$, disregarding one data point should not influence the inference results. This, in turn, means that one can alternatively use the approximation 
\begin{eqnarray}
     W_n &:=& \frac{1}{n}\sum_{i=1}^n W_{(i)} + o_P(1), ~~~~~~~
    V_{n} := \frac{1}{n}\sum_{i=1}^n  V_{(i)} + o_P(1)\\
     \mathcal{U}_n &:=& \frac{1}{n}\sum_{i=1}^n  \mathcal{U}_{(i)} + o_P(1)
\end{eqnarray}
Which shows that it is plausible that $W_n,V_n$ and $\mathcal{U}_n$ concentrate around deterministic values, or are self averaging in physical jargon, which can be stated as 
\begin{eqnarray}
\label{concentration_hyp}
    W_n := w+ o_P(1), ~~~~~~~
    V_{n} := v + o_P(1), ~~~~~~~
    \mathcal{U}_n :=\tilde{u} + o_P(1)
\end{eqnarray}
The remaining task is to determine the values $(w,v,\tilde{u})$, which according to the above arguments should satisfy a set of coupled self-consistency equations.

\subsection{Self-consistent equations for the overlaps}

Given the Gibbs measure 
\begin{equation}
    \bbeta|\{T_i,\bm{\mathcal{X}}_i\} \sim \frac{1}{Z_n}\rme^{\gamma\big(\sum_{i=1}^n \rho(T_i|\bm{\mathcal{X}}_i\cdot\bm{\beta},\bm{\sigma})-\frac{1}{2}\eta\| \bbeta\|^2\big)}
\end{equation}
it is easy to derive the following equation via integration bu parts in $\bbeta$:
\begin{equation}
    \gamma \eta  \langle\bm{\beta}\rangle =\gamma \sum_{i=1}^n \bm{\mathcal{X}}_i \langle\varphi(T_i|\bm{\mathcal{X}}_i\cdot\bm{\beta},\bm{\sigma})\rangle
\end{equation}
By taking the scalar product with $\bbeta^{\sim}_0$ and dividing by $\|\bbeta^{\sim}_0\|=S$, we get 
\begin{equation}
    \gamma \eta  \bbeta_0^{\sim}\cdot\langle\bm{\beta}\rangle =\gamma \sum_{i=1}^n \bbeta^{\sim}_0\cdot\bm{\mathcal{X}}_i \langle\varphi(T_i|\bm{\mathcal{X}}_i\cdot\bm{\beta},\bm{\sigma})\rangle
\end{equation}
Using the concentration/self averaging hypothesis (\ref{concentration_hyp}) and taking the expectation over the disorder (i.e the data-set), we then obtain a self consistent equation for $w$
\begin{eqnarray}
    w &:=& \mathbb{E}\Big[\bm{\beta}^{\sim}_0\cdot\langle \bm{\beta}\rangle\Big]/S = \frac{1}{\eta} \sum_{i=1}^n \mathbb{E}\Big[\langle\bm{\mathcal{X}}_i\cdot\bm{\beta}^{\sim}_0 \varphi(T_i|\bm{\mathcal{X}}_i\cdot\bm{\beta},\bm{\sigma})\rangle\Big]/S 
    \nonumber\\
    &=&\frac{n}{\eta}  \mathbb{E}\Big[Z_0\mathbb{E}_{\xi}[ \varphi(T|\xi,\bm{\sigma})]\Big]
\end{eqnarray}
with 
\begin{equation}
    \mathbb{E}_{\xi}\Big[ f(\xi)\Big] := \frac{\int\!{\rm d} \xi~ \rme^{-\frac{1}{2}\big(\frac{\xi-vQ-wZ_0}{\tilde{u}}\big)^2 + \gamma \rho(T|\xi,\bm{\sigma})} f(\xi)}{\int\!{\rm d} \xi~ \rme^{-\frac{1}{2}\big(\frac{\xi-vQ-wZ_0}{\tilde{u}}\big)^2 + \gamma \rho(T|\xi,\bm{\sigma})}}
\end{equation}
We see that 
\begin{equation}
    \gamma \mathbb{E}_{\xi}\Big[ \varphi(T|\xi,\bm{\sigma})\Big] =  \mathbb{E}_{\xi}\Big[\Big(\frac{\xi-vQ -wZ_0}{\tilde{u}^2}\Big)\Big]
\end{equation}
This implies 
\begin{equation}
    w = \frac{n}{\gamma \eta \tilde{u}^2}\bigg(\mathbb{E}\Big[Z_0\mathbb{E}_{\xi}[ \xi]\Big] -w\bigg)
\end{equation}
and thus
\begin{equation}
    w(1+\gamma\tilde{u}^2 \eta/n) = \mathbb{E}\Big[Z_0\mathbb{E}_{\xi}[ \xi]\Big]
\end{equation}
Similarly 
\begin{eqnarray}
    \fl \gamma \eta  \mathbb{E}\Big[\langle\bm{\beta}\cdot\bm{\beta}\rangle\Big]  &=& \gamma \eta (v^2+w^2+\tilde{u}^2) =\gamma \sum_{i=1}^n \mathbb{E}\Big[\langle\bm{\beta}\cdot\bm{\mathcal{X}}_i \varphi(T_i|\bm{\mathcal{X}}_i\cdot\bm{\beta},\bm{\sigma})\rangle\Big]  +  \gamma p \nonumber\\
    \fl &=&\gamma p + \gamma n \mathbb{E}\Big[\mathbb{E}_{\xi}[\xi \varphi(T|\xi,\bm{\sigma})]\Big]\nonumber\\
    \fl &=& \gamma (p-n)+n \mathbb{E}\Big[\mathbb{E}_{\xi}[\xi(\xi-vQ-wZ_0)]\Big]/\tilde{u}^2 
    \label{r11}
\end{eqnarray}
where we have used 
\begin{equation}
    \gamma \mathbb{E}\Big[\mathbb{E}_{\xi}[\xi \varphi(T|\xi,\bm{\sigma})]\Big] = \mathbb{E}\Big[\mathbb{E}_{\xi}[\xi(\xi-vQ-wZ_0)]\Big]/\tilde{u}^2 - 1
\end{equation}
Also
\begin{eqnarray}
  \fl \gamma \eta \mathbb{E}\Big[\langle\bm{\beta}\rangle\cdot \langle\bm{\beta}\rangle\Big] &=& \gamma \eta (v^2 + w^2)= \gamma \sum_{i=1}^n \mathbb{E}\Big[\langle\bm{\mathcal{X}}_i \cdot\bm{\beta}\rangle\langle\langle\varphi(T_i|\bm{\mathcal{X}}_i\cdot\bm{\beta},\bm{\sigma})\rangle\Big] \nonumber\\
   \fl &=&\gamma n \mathbb{E}\Big[\mathbb{E}_{\xi}[\xi] \mathbb{E}_{\xi}[\varphi(T|\xi,\bm{\sigma})]\Big]= \frac{n}{\tilde{u}^2} \mathbb{E}\Big[\mathbb{E}_{\xi}[\xi] \mathbb{E}_{\xi}[\xi-vQ -wZ_0]\Big]
    \label{r12}
\end{eqnarray}
Now take the difference of (\ref{r11}) and (\ref{r12}) to find
\begin{equation}
    \gamma \eta  \tilde{u}^2 = p-n + n\mathbb{E}\Big[\mathbb{E}_{\xi}[\xi^2] -\mathbb{E}_{\xi}[\xi]^2\Big]/\tilde{u}^2=p-n+ n\mathbb{E}\Big[Q\mathbb{E}_{\xi}[\xi]\Big]/v 
\end{equation}
dividing by $n$ and rearranging the terms gives
\begin{equation}
    v(1-\zeta) + v\gamma\tilde{u}^2\eta /n  =\mathbb{E}\Big[Q\mathbb{E}_{\xi}[\xi]\Big]
\end{equation}
We obtain another self consistent equation by substitution
\begin{eqnarray}
 \hspace*{-10mm}    \eta\tilde{u}^2\gamma (v^2+w^2)/n &=& \mathbb{E}\Big[\mathbb{E}_{\xi}[\xi]^2\Big] - v\mathbb{E}\Big[Q\mathbb{E}_{\xi}[\xi]\Big] - w \mathbb{E}\Big[Z_0\mathbb{E}_{\xi}[\xi]\Big] \nonumber \\
 \hspace*{-10mm}     &= &\mathbb{E}\Big[\mathbb{E}_{\xi}[\xi^2]\Big] - v^2\Big(1-\zeta + \frac{\gamma\tilde{u}^2\eta}{n}\Big) -w^2(1+\frac{\gamma\tilde{u}^2 \eta}{n})
\end{eqnarray}
Rearranging terms leads to
\begin{equation}
    \mathbb{E}\Big[\mathbb{E}_{\xi}[\xi]^2\Big] = v^2(1-\zeta) + w^2 + 2 \tilde{u}^2\gamma (v^2+w^2)\eta/n
\end{equation}
In summary we have obtained the following equations 
\begin{eqnarray}
     \mathbb{E}\Big[\mathbb{E}_{\xi}[\xi]^2\Big] &=& v^2(1-\zeta) + w^2 + 2 \tilde{u}^2\gamma (v^2+w^2)\eta'\zeta \label{sp1}\\
     v(1-\zeta) + v\gamma\tilde{u}^2\eta'\zeta &=&\mathbb{E}\Big[Q\mathbb{E}_{\xi}[\xi]\Big]\label{sp2}\\
    w(1+\gamma\tilde{u}^2 \eta'\zeta) &=& \mathbb{E}\Big[Z_0\mathbb{E}_{\xi}[ \xi]\Big]\label{sp3}
\end{eqnarray}
where we used the scaling $\eta = \eta' p$ as in the main text. 
These equations are valid under our stated assumptions and approximations, and apply so far to any inverse temperature $\gamma$. Their solution  gives us the values  $\tilde{u}_{\star}, v_{\star}, w_{\star}$ as functions of $\gamma$, $\bsigma$, $S$ and $\zeta$. 
The value of $\bsigma$ though, is generally unknown and must be estimated from the data. 

\subsection{Equations for the nuisance parameters}

Following a strategy similar to the one followed so far, we assume also the estimator $\hat{\bsigma}_n$ for the nuisance parameters to assume deterministic values in the asymptotic limit $n,p\rightarrow\infty \quad p/n=\zeta$. 
Under the assumption of concentration of the overlaps (\ref{concentration_hyp}), we expect that  the internal energy per data point
\begin{equation}
\label{F_n}
    F_n(\bsigma) = \Big\langle\frac{1}{n}\sum_{i=1}^n \rho (T_i|\bm{\mathcal{X}}_i\cdot\bbeta,\bm{\sigma}) - \frac{1}{2} \eta'\zeta \|\bm{\beta}\|^2\Big\rangle 
\end{equation}
will be self averaging. i.e. concentrate on a deterministic function $f(\bsigma)$:
\begin{equation}
    F_n(\bsigma) = f(\bsigma) + o_P(1)
\end{equation}
Using (\ref{cavity_exp_2}) for the cavity average, together with (\ref{concentration_hyp}) for the overlap values, gives  
\begin{equation}
    F_n(\bsigma) = \frac{1}{n}\sum_{i=1}^n \mathbb{E}_{\xi}\Big[\rho (T_i|
    \xi,\bm{\sigma})\Big] - \frac{1}{2} \eta (w^2+v^2+\tilde{u}^2) + o_P(1)
\end{equation}
By the law of large numbers we expect that also this quantity will concentrate on
\begin{equation}
\label{int_en}
    f(\bsigma) = \mathbb{E}_{T,Q,Z_0}\bigg[ \mathbb{E}_{\xi}\Big[\rho (T|
    \xi,\bm{\sigma})\Big]\bigg] -\frac{1}{2} \eta (w^2+v^2+\tilde{u}^2)
\end{equation}
Hence the estimator for the nuisance parameter $\hat{\bsigma}_n$, the value that maximizes $F_n$, will be {close} to the  deterministic value that maximizes $f$. The latter solves 
\begin{equation}
    \frac{\rmd}{\rmd \bsigma} f(\bsigma) = \frac{\partial f(\bsigma)}{\partial \bsigma}  + \frac{\partial f}{\partial w}  \frac{\partial w }{\partial \bsigma}+ \frac{\partial f}{\partial v}  \frac{\partial v}{\partial \bsigma}+ \frac{\partial f}{\partial \tilde{u}}  \frac{\partial \tilde{u} }{\partial \bsigma} =0
\end{equation}
It can be verified by direct calculation that $\frac{\partial}{\partial w} f$,$\frac{\partial}{\partial v} f $ and $\frac{\partial}{\partial \tilde{u}} f $ vanish at the solution of the self consistent equations for the overlaps (as expected).
Hence the  deterministic value that maximizes $f$ at inverse temperature $\gamma$, is the solution of  the following equation, which must be solved simultaneously with the self consistent equations for the overlaps:
\begin{equation}
\label{sp4}
    \mathbb{E}_{T,Q,Z_0}\Big[ \mathbb{E}_{\xi}\big[\nabla_{\bsigma}\rho (T|\xi,\bm{\sigma})\big]\Big]  ={\bf 0}
\end{equation}

\subsection{The zero temperature limit}

What remains is to take the limit $\gamma \rightarrow\infty$ in our equations to recover the overlaps appearing in the representation of the MAP/PML estimator $\tilde{u}_{\star},v_{\star},w_{\star}$ and the deterministic values of the nuisance parameters estimator $\bsigma_{\star}$. 
Since the distribution of the PML/MAP estimator depends on the zero temperature limit of various quantities  
\begin{eqnarray}
     \fl R^{(n)}_{0,1} &=&\mathbb{E}\Big[\frac{\bbeta^{\sim}_0\cdot\hat{\bbeta}^{\sim}_n}{\|\bbeta^{\sim}_0\|^2}\Big] + o_P(1),\qquad \mathbb{E}\Big[\frac{\bbeta^{\sim}_0\cdot\hat{\bbeta}^{\sim}_n}{\|\bbeta^{\sim}_0\|^2}\Big] = \lim_{\gamma\rightarrow \infty}  \mathbb{E}\Big[\frac{\bbeta^{\sim}_0\cdot\langle\bbeta\rangle}{\|\bbeta^{\sim}_0\|^2}\Big]= \lim_{\gamma\rightarrow \infty} \frac{w}{S}\\
     \fl R^{(n)}_{1,1} &=& \mathbb{E}\Big[ \hat{\bbeta}^{\sim}_n\cdot\hat{\bbeta}^{\sim}_n\Big] +o_P(1), \qquad \mathbb{E}\Big[ \hat{\bbeta}^{\sim}_n\cdot\hat{\bbeta}^{\sim}_n\Big]= \lim_{\gamma \rightarrow \infty} \mathbb{E}\Big[ \langle \bbeta\cdot\bbeta \rangle \Big] = \lim_{\gamma\rightarrow \infty} (v^2+w^2+\tilde{u}^2)
\end{eqnarray}
together with the value $\bsigma_{\star}$ that minimizes the internal energy (\ref{int_en}) at zero temperature, we must compute the limit $\gamma \rightarrow \infty$ of equations (\ref{sp1}, \ref{sp2}, \ref{sp3}, \ref{sp4}).\\
To do so, we consider the general expression
\begin{equation}
    \lim_{\gamma\rightarrow \infty} \mathbb{E}_{\xi}[f(\xi)] = \lim_{\gamma\rightarrow \infty} \frac{\int\! {\rm d} \xi~ \rme^{-\frac{1}{2}\big(\frac{\xi-vQ-wZ_0}{\tilde{u}}\big)^2 + \gamma \rho(T|\xi,\bm{\sigma})} f(\xi)}{\int\!{\rm d} \xi~ \rme^{-\frac{1}{2}\big(\frac{\xi-vQ-wZ_0}{\tilde{u}}\big)^2 + \gamma \rho(T|\xi,\bm{\sigma})} }
\end{equation}
We assume the scaling $\tilde{u}^2 = u^2/\gamma $ suggested by equations (\ref{sp1},\ref{sp2},\ref{sp3}), where $\tilde{u}^2$ appears always multiplied by $\gamma$. We can now proceed via the Laplace argument:
\begin{equation}
     \lim_{\gamma\rightarrow \infty} \frac{\int\!{\rm d} \xi~ \rme^{-\gamma\big(\frac{1}{2}\frac{(\xi-vQ-wZ_0)^2}{u^2} - \rho(T|\xi,\bm{\sigma})\big)} f(\xi) \ }{\int\!{\rm d} \xi~ \rme^{-\gamma\big(\frac{1}{2}\frac{(\xi-vQ-wZ_0)^2}{u^2} - \rho(T|\xi,\bm{\sigma})\big)} } = f(\xi_{\star})
    \label{laplace}
\end{equation}
with $\xi_{\star}$ denoting the proximal mapping of the function $\rho(T|.,\bsigma)$, defined as
\begin{equation}
    \xi_{\star} = \underset{\xi}{\arg\min} \Big\{\frac{1}{2}\frac{(\xi-vQ-wZ_0)^2}{u^2} - \rho(T|\xi,\bm{\sigma})\Big\}
\label{xi_star}
\end{equation}
(as always under appropriate regularity conditions on $f$).
It is now easy to take the limit in  equations (\ref{sp1}, \ref{sp2}, \ref{sp3}, \ref{sp4}), obtaining
\begin{eqnarray}
 \mathbb{E}\Big[\xi^2_{\star}\Big] &=& v^2(1-\zeta) + w^2 + 2 u^2 (v^2+w^2)\eta'\zeta\\
 v(1-\zeta) + v u^2\eta'\zeta &=&\mathbb{E}\Big[Q\xi_{\star}\Big]\\
 w(1+u^2 \eta'\zeta) &=& \mathbb{E}\Big[Z_0 \xi_{\star}\Big]\\
 \nabla_{\bsigma} \rho(T|\xi_{\star},\bsigma) &=&{\bf 0}
\end{eqnarray}
Some additional algebraic steps then result in equations (\ref{RS1}, \ref{RS2}, \ref{RS3}) in the main text.

\section{Logit regression model}
\label{appendix:logit}

\subsection{Working out the RS equations}

The Logit regression model is defined by 
\begin{equation}
    T|\mathbf{X} \sim \frac{\rme^{T(\mathbf{X}\cdot\bm{\beta})}}{2\cosh(\mathbf{X}\cdot\bm{\beta})}
\end{equation}
where $T\in \{-1,+1\}$. The log-likelihood is 
 $\log p (T|\mathbf{X}\cdot\bm{\beta})= T(\mathbf{X}\cdot\bm{\beta}) - \log [2\cosh(\mathbf{X}\cdot\bm{\beta})]$, so the 
 first RS equation (\ref{logit_RS1}) follows by substitution into inside formula (\ref{RS1}) of 
\begin{equation}
    \xi_{\star} = \nu Q+\omega Z_0 +\mu^2 T -\mu^2\tanh(x)
\end{equation}
This gives
\begin{equation}
    \frac{\nu^2 \zeta}{(1-\tau'\zeta \mu^2)^2} = \mu^4~\mathbb{E}_{Y,Q,Z_0}\Big[\Big( T-\tau  \frac{\nu  Q+\omega Z_0}{1-\tau \mu^2} -  \tanh(x_{\star})\Big)^2\Big]
\end{equation}
For equation (\ref{logit_RS2}) we need the derivative 
$\partial \xi_{\star}/\partial Q =\partial x_{\star}/\partial Q$, which  
is obtained by differentiating the equation that defines $x_{\star}$:
\begin{equation}
    \frac{\partial \xi_{\star}}{\partial Q} = \nu - \frac{\mu^2}{\cosh^2(x_{\star})}\frac{\partial \xi_{\star}}{\partial Q}
\end{equation}
with solution
\begin{equation}
    \frac{\partial \xi_{\star}}{\partial Q} =  \frac{\nu\cosh^2(x_{\star})}{\mu^2 + \cosh^2(x_{\star})}
\end{equation}
Substitution into (\ref{RS2}), followed by division by the common factor $\nu$, leads after some simple rewriting to
\begin{equation}
    \zeta\Big(1-\eta'\mu^2/(1-\tau'\zeta\mu^2)\Big) = \mathbb{E}_{T,Q,Z_0}\Big[\frac{u^2}{u^2 + \cosh^2(x_{\star})}\Big]
\end{equation}
Next we compute
\begin{equation}
    \frac{\partial }{\partial Z_0} \log p(T|SZ_0) = S\Big(T - \tanh(SZ_0)\Big)
\end{equation}
and insert the result into  
\begin{equation}
    \frac{\zeta\omega}{1-\tau'\zeta\mu^2} =\mathbb{E}_{T,Q,Z_0}\Big[ \xi_{\star}\frac{\partial}{\partial Z_0} \log p (T|SZ_0)\Big] 
\end{equation}
to give
\begin{equation}
    \zeta \omega/S \zeta = (1-\tau'\zeta\mu^2) \mathbb{E}_{T,Q,Z_0}\big[(x-\phi)\big(T-\tanh(SZ_0)\big)\big] 
\end{equation}
We note that 
\begin{equation}
    \mathbb{E}_{T,Q,Z_0}[T] = \mathbb{E}_{T,Z_0}[\tanh(SZ_0)]
\end{equation}
so 
\begin{equation}
    \zeta \omega/S \zeta = (1-\tau'\zeta\mu^2) \mathbb{E}_{T,Q,Z_0}\big[x\big(T-\tanh(SZ_0)\big)\big]
\end{equation}

\subsection{Numerical solution of the RS equations}

The numerical solution of the RS equations is obtained by means of fixed point iteration. The system (\ref{logit_RS1},\ref{logit_RS2},\ref{logit_RS3}) can be seen as a dynamical mapping for $\mathbf{x} = (\zeta,\nu,\omega)$:
\begin{equation}
    \mathbf{x}_t= \mathbf{f}(\mathbf{x}_{t-1};\mu^2,S)
\end{equation}
The solution of (\ref{logit_RS1},\ref{logit_RS2},\ref{logit_RS3})  corresponds to a fixed point of the dynamical system, which is found by iterating the mapping above until some convergence criterion is met. In our case we chose to set the tolerance at $\|\mathbf{x}_t-\mathbf{x}_{t-1}\|\leq 1\cdot 10^{-10}$.
We see that 
\begin{equation}
    \mathbb{E}_{T,Q,Z_0}\Big[g(T,Q,Z_0)\Big] = \mathbb{E}_{Q,Z_0}\Big[\sum_{t=-1,1}p(t|SZ_0)g(t,Q,Z_0)\Big]
\end{equation}
so we just need to compute two Gaussian integrals. We computed the mapping $\mathbf{f}$ via Gauss-Hermite quadratures, at order $40$ to achieve a reasonably high accuracy (since the functions appearing inside the integrals are not simple polynomials of a fixed degree).

\section{Weibull Proportional Hazards model}
\label{appendix:weibull}

\subsection{Working out the RS equations}

We first derive the RS equations for an arbitrary {parametric} Proportional Hazards model, for which the  conditional response density given the covariates is
\begin{equation}
\label{pdfT}
    p(t|\mathbf{X}\cdot\bm{\beta},\bm{\sigma}) = h(t|\bm{\sigma})\rme^{\mathbf{X}^T\bm{\beta}}\rme^{-H(t|\bm{\sigma})\rme^{\mathbf{X}^T\bm{\beta}}}
\end{equation}
with $t\geq 0$ and $H(t|\bm{\sigma})=\int_0^t\!\rmd t^\prime~h(t^\prime|\bm{\sigma})$. 
Equation (\ref{RS_eqs}) new reads
\begin{equation}
    \mu^{-2}(\xi_{\star}-\nu Q-\omega Z_0) = 1 - H(T|\bm{\sigma})\rme^{\xi_{\star}}
\end{equation}
with the following solution, involving Lambert's $W$-function \cite{lambert_function}:
\begin{equation}
\label{xi_weib}
    \xi_{\star} = \nu Q +\omega Z_0 + \mu^2 - W\big(H(T|\bm{\sigma})\mu^2\rme^{\mu^2 + \nu Q+\omega Z_0}\big)
\end{equation}
By substituting (\ref{xi_weib}) into (\ref{RS1}) and writing $(u,v,w)$ in terms of $(\mu,\nu,\omega)$ one obtains
\begin{equation}
    \label{w1}
    \fl \frac{\zeta\nu^2}{(1-\tau \mu^2)^2}= \mathbb{E}_{T,Q,Z_0}\Big[\Big(\mu^2 - W\Big(H(T|\bm{\sigma})\mu^2\rme^{\mu^2 + \nu Q+\omega Z_0}\Big) - \frac{\tau \mu^2}{1 - \tau \mu^2}(\nu Q +\omega Z_0)\Big)^2\Big]
\end{equation}
To obtain equation (\ref{RS2}) we need $\partial \xi_\star/\partial Q$. This is computed by the inverse function rule
\begin{equation}
    \frac{\partial\xi_\star}{\partial Q} = \nu \Bigg( 1- \frac{W\big(H(T|\bm{\sigma})\mu^2\rme^{\mu^2 + \nu Q+\omega Z_0}\big)}{1+W\big(H(T|\bm{\sigma})\mu^2\rme^{\mu^2 + \nu Q+\omega Z_0}\big)}\Bigg)
\end{equation}
After substituting and dividing for the common factor $\nu$ we obtain
\begin{equation}
    \label{w2}
    \fl \frac{1}{1-\tau'\zeta\mu^2} \bigg(1-\zeta\Big(1-\frac{\mu^2\eta'}{1-\tau'\zeta\mu^2}\Big)\bigg)= 1- \mathbb{E}_{T,Q,Z_0}\Big[ \frac{W\big(H(T|\bm{\sigma})\mu^2\rme^{\mu^2 + \nu Q+\omega Z_0}\big)}{1+W\big(H(T|\bm{\sigma})\mu^2\rme^{\mu^2 + \nu Q+\omega Z_0}\big)}\Big]  
\end{equation}
For the third RS equation (\ref{RS3}) we need to compute 
\begin{equation}
    \frac{\partial}{\partial Z_0} \log p(T|Z_0,\phi_0,
    \sigma_0) = S - S H(T|SZ_0,\phi_0,\sigma_0)\rme^{SZ_0}
\end{equation}
Substitution into  (\ref{RS3}) gives
\begin{equation}
    \label{w3}
    \fl \frac{\omega \zeta }{1-\tau\mu^2} =  S~\mathbb{E}_{T,Q,Z_0}\Big[ W\big(H(T|\bm{\sigma})\mu^2\rme^{\mu^2 + \nu Q+\omega Z_0}\big)\Big(H(T|SZ_0,\phi_0,\sigma_0)\rme^{SZ_0}-1\Big)\Big]
\end{equation}

Note that the distribution (\ref{pdfT}) of $T|Z_0$ has the property that upon  changing to the new variable 
\begin{equation}
\label{change_of_var}
    Z = H(T|\bm{\sigma}_0)\rme^{SZ_0} 
\end{equation}
one will have $Z \sim Exp(1)$, $T = H^{-1}(Z\rme^{-SZ_0}|\bm{\sigma}_0)$ and 
\begin{equation}
    H(T|\bm{\sigma}) = H(H^{-1}(Z\rme^{-SZ_0}|\bm{\sigma}_0)|\bm{\sigma})
\end{equation}
Hence the first three RS equations derived above (\ref{w1},\ref{w2},\ref{w3}) can be rewritten as 
\begin{eqnarray}
\hspace*{-20mm} 
\frac{\zeta \nu^2}{(1-\tau\zeta\mu^2)^2} &=& \mathbb{E}_{Z,Q,Z_0}\Big[\Big(\mu^2 - W\big(\mu^2 H(H^{-1}(Ze^{-SZ_0}|\bm{\sigma}_0)|\sigma)\rme ^{\mu^2+\nu Q + \omega Z_0}\big)\Big)^2\Big]\\
\hspace*{-20mm}     &&\hspace*{-26mm} \frac{1}{1-\tau'\zeta\mu^2}\Big(1-\zeta\Big(1-\frac{\mu^2\eta'}{1-\tau\zeta\mu^2}\Big)\Big)\nonumber \\
 \hspace*{-20mm}    &=& 1-\mathbb{E}_{T,Q,Z_0}\Bigg[\frac{W\big(\mu^2 H(H^{-1}(Ze^{-SZ_0}|\bm{\sigma}_0)|\sigma)\rme ^{\mu^2+\nu Q + \omega Z_0}\big)}{1+W\big(\mu^2 H(H^{-1}(Ze^{-SZ_0}|\bm{\sigma}_0)|\sigma)\rme ^{\mu^2+\nu Q + \omega Z_0}\big)}\Bigg]\\
  \hspace*{-20mm}   \frac{\omega}{1-\tau\zeta\mu^2} \zeta &=& S\mathbb{E}_{T,Q,Z_0}\Big[ W\big(\mu^2 H(H^{-1}(Ze^{-SZ_0}|\bm{\sigma}_0)|\sigma)\rme ^{\mu^2+\nu Q + \omega Z_0}\Big)\big)(Z-1)\Big]
\end{eqnarray}
\vspace*{3mm}

To work out the RS equations for the nuisance parameters (\ref{RS4}) we must choose an explicit parametric form for the base hazard rate.
Upon  assuming the integrated base hazard to have the Weibull parametrization, 
\begin{equation}
    H(T|\phi,\sigma) = T^{1/\sigma} \rme^{\phi/\sigma}
\end{equation}
we obtain
\begin{equation}
\label{weib_cov}
    H(H^{-1}(Ze^{-SZ_0}|\phi_0,\sigma_0)|\phi,\sigma) = Z^{\sigma_0/\sigma} \rme^{-SZ_0 \sigma_0/\sigma + (\phi-\phi_0)/\sigma} 
\end{equation}
Substituting (\ref{weib_cov}) into (\ref{w1}, \ref{w2}) then directly leads to equations (\ref{weib1}, \ref{weib2}) in the main text.
We simplify (\ref{w3}) via integration by parts, use  (\ref{w2}), and after some minor algebraic simplifications obtain (\ref{weib3}).
For the remaining equations one has to compute 
\begin{eqnarray}
    \fl \frac{\partial }{\partial \phi} \log p(T|\xi_*,\phi,\sigma) &=& -\frac{1}{\sigma} - \frac{1}{\sigma\mu^2} W\big( Z^{\sigma_0/\sigma} \mu^2\rme^{\mu^2 +(\omega-S\frac{\sigma_0}{\sigma})Z_0 + \nu Q+(\phi-\phi_0)/\sigma}\big) \\
    \fl \frac{\partial}{ \partial \sigma} \log p(T|\xi_*,\phi,\sigma) &=& -\frac{1}{\sigma} - \frac{1}{\sigma} \Big(\frac{\phi-\phi_0}{\sigma}+\frac{\sigma_0}{\sigma} \log Z -\frac{\sigma_0}{\sigma} SZ_0\Big)\nonumber\\
    &\times& \bigg(\frac{1}{\mu^{2}}W\big( Z^{\sigma_0/\sigma} \mu^2\rme^{\mu^2 +(\omega-S\frac{\sigma_0}{\sigma})Z_0 + \nu Q+(\phi-\phi_0)/\sigma}\big)-1\bigg)
\end{eqnarray}
Upon taking the expectation and setting the previous two derivatives to zero we obtain the remaining RS equations (\ref{weib4}, \ref{weib5})
\begin{eqnarray}
     \fl \mu^2 &=& \mathbb{E}_{Z,Q,Z_0}\Big[W\Big(\mu^2\rme^{\mu^2} Z^{\sigma_0/\sigma} \rme^{(\omega-S\frac{\sigma_0}{\sigma})Z_0 + \nu Q} \rme^{(\phi-\phi_0)/\sigma}\Big)\Big]\\
     \fl \frac{\sigma}{\sigma_0} &=& \mathbb{E}_{Z,Q,Z_0}\Big[\Big(\log Z\!-\! SZ_0 \Big)\Big(\frac{1}{\mu^2}W\Big(\mu^2\rme^{\mu^2} Z^{\sigma_0/\sigma} \rme^{(\omega-S\frac{\sigma_0}{\sigma})Z_0 + \nu Q} \rme^{(\phi-\phi_0)/\sigma}\Big) -1 \Big)\Big]
\end{eqnarray}

\subsection{Numerical solution of the RS equations}

The numerical solution of the RS equations is again obtained via fixed point iteration. The system (\ref{weib1},\ref{weib2},\ref{weib3},\ref{weib4},\ref{weib5}) is interpreted as the fixed-point condition of a dynamical mapping of $\mathbf{x} = (\zeta,\nu,\omega,\mu,\sigma/\sigma_0)$: 
\begin{equation}
    \mathbf{x}_t= \mathbf{f}(\mathbf{x}_{t-1};(\phi-\phi_0)/\sigma,S)
\end{equation}
The mapping is iterated until some convergence criterion is met, here chosen to be   $\|\mathbf{x}_t-\mathbf{x}_{t-1}\|\leq 1\cdot 10^{-10}$.
We see that all the RS equations involve terms of the form
\begin{equation}
    \mathbb{E}_{Z,Q,Z_0}\Big[g(Z,Q,Z_0)\Big] 
\end{equation}
involving two Gaussian and one  exponential integral. We computed these via Gauss-Hermite and Gauss-Laguerre quadratures, to order $40$ to achieve sufficient accuracy.

\end{document}